\documentclass{amsart}
\usepackage{amssymb,latexsym}

\numberwithin{equation}{section}

\newtheorem{theorem}{Theorem}[section]
\newtheorem{proposition}[theorem]{Proposition}
\newtheorem{corollary}[theorem]{Corollary}
\newtheorem{remark}[theorem]{Remark}
\newtheorem{lemma}[theorem]{Lemma}
\newtheorem{example}[theorem]{Example}
\newtheorem{definition}[theorem]{Definition}

\def\proof{\smallskip\noindent {\bf Proof. }}
\def\endproof{\hfill$\square$\medskip}

\newcommand{\doublesubscript}[3]{
\displaystyle\mathop{\displaystyle #1_{#2}}_{#3}
}

\def\AA{\mathcal{A}}
\def\BB{\mathcal{B}}
\def\MM{\mathcal{M}}
\def\ZZ{\mathbb{Z}}
\def\CC{\mathbb{C}}
\def\RR{\mathbb{R}}
\def\LL{\mathcal{L}}
\def\l{\ell}
\def\wnot{w_\mathrm{o}}
\def\gg{\mathfrak{g}}

\def\hh{\mathfrak{h}}
\def\ii{\mathbf{i}}
\def\jj{\mathbf{j}}
\def\tt{\mathbf{t}}

\newcommand{\mat}[4]{\left(\!\!\begin{array}{cc}
#1 & #2 \\
#3 & #4 \\
\end{array}\!\!\right)}

\begin{document}

\title{Tensor product multiplicities, canonical bases and totally positive varieties}

\author{Arkady  Berenstein}
\address{Department of Mathematics, Harvard University,
Cambridge, MA 02138}
\email{arkadiy@math.harvard.edu}

\author{Andrei Zelevinsky}
\address{\noindent Department of Mathematics, Northeastern University,
  Boston, MA 02115}
\email{andrei@neu.edu}

\thanks{The research of both authors was supported in part
by NSF grants (A.B.: \#DMS-9970533; A.Z.: \#DMS-9625511 and
\#DMS-9971362.)}

\date{December 1, 1999}

\maketitle

\makeatletter
\renewcommand{\@evenhead}{\tiny \thepage \hfill  A.~BERENSTEIN and  A.~ZELEVINSKY \hfill}

\renewcommand{\@oddhead}{\tiny \hfill  TENSOR PRODUCT MULTIPLICITIES ...
 \hfill \thepage}
\makeatother

\tableofcontents

\newpage

\section{Introduction}
\label{sec:introduction}

This paper continues and to some extent concludes the project initiated twelve years ago
in \cite{bz88}.
The original goal of this project was to find explicit ``polyhedral" combinatorial
expressions for multiplicities in the tensor product of two simple
finite-dimensional modules over a complex semisimple Lie algebra $\gg$.
Here ``polyhedral" means that the multiplicity in question is to be expressed as the
number of lattice points in some convex polytope, or in more down-to-earth terms,
the number of integer solutions of a system of linear equations and inequalities.
The tensor product multiplicities are often called generalized Littlewood-Richardson
coefficients because for type $A_r$ they are given by the classical Littlewood-Richardson
rule (a polyhedral version of the rule was given in \cite{gz2},
and a much more symmetric version was given in \cite{bz92}).
Conjectural polyhedral expressions for these multiplicities in the
case of classical Lie algebras were given in \cite{bz88}.
A uniform combinatorial description of these multiplicities (even in a more general context
of Kac-Moody algebras) has been obtained by P.~Littelmann in \cite{lit1}; his answer
is in terms of the so-called path model.
It is however not an easy task to transform this description into
a polyhedral one; in particular, it is still not clear how to
deduce the conjectural expressions in \cite{bz88} from Littelmann's results.

In this paper, we explicitly construct a family of polyhedral expressions
for tensor product multiplicities for an arbitrary semisimple Lie algebra $\gg$.
To be more precise, we associate two such expressions
to every reduced word of $\wnot$, the longest element of the
Weyl group, and also produce two ``universal" expressions which we call
\emph{tropical Pl\"ucker models}.
It can be shown that these expressions include as special cases the
conjectural expressions in \cite{bz88}.
As another application, we obtain a family of polyhedral expressions
for the multiplicities that occur when one restricts a simple
finite-dimensional $\gg$-module to the Levi subalgebra of some
parabolic subalgebra in $\gg$.
Our answers use a new combinatorial concept of $\ii$-trails which resembles
Littelmann's paths but seems to be more tractable.

~From the beginning, another (and maybe more important) aim of our project was to develop
combinatorial understanding of ``good" bases in
finite-dimensional representations of $\gg$.
A preliminary concept of a ``good" or ``canonical" basis was introduced in
\cite{gz1,gz2} and independently in \cite{bac} but it was only given a firm mathematical
foundation in pioneering works of G.~Lusztig and M.~Kashiwara on canonical bases in quantum groups.
Thus the aim of our project can be more precisely formulated as follows:
understand the combinatorial structure of Lusztig's canonical bases or, equivalently of
Kashiwara's global bases.
Although Lusztig's and Kashiwara's approaches were shown by Lusztig to be equivalent
to each other, they lead to different combinatorial parametrizations of the canonical bases.
One of the main results of the present paper is an explicit description
of the relationship between these parametrizations.

In Lusztig's approach, every reduced word $\ii$ for $\wnot$ gives rise to a parametrization
(``of the PBW-type") of the canonical basis $\BB$ by the set
$\ZZ_{\geq 0}^m$ of all $m$-tuples of nonnegative
integers, where $m = \l (\wnot)$ is the number of positive roots.
Kashiwara's approach is closely related to another family of parametrizations
which we call string parametrizations; they were introduced and studied in \cite{bz93}.
String parametrizations are also associated to reduced words for $\wnot$; but this time
every such reduced word $\ii$ gives rise to a bijection between $\BB$
and the set of all lattice points
of some rational polyhedral convex cone $C_\ii$ in $\RR^m$.
In this paper we obtain an explicit description of these cones for an arbitrary semisimple
Lie algebra $\gg$ and any reduced word $\ii$.
Such a description was previously known in some special cases only: in \cite{bz93}, it was given
for a special reduced word $(1, 2,1, 3,2,1, \dots)$ for type $A$; in \cite{lit}
it was extended to a special choice of a reduced word for any semisimple $\gg$,
while in \cite{gp} it was given for any
reduced word for type $A$; a more general setting of Kac-Moody algebras was discussed
in \cite{nakzel97} but the results there were inconclusive.

Our approach to the above problems is based on a remarkable
observation by G.~Lusztig that combinatorics of the canonical basis
is closely related to geometry of the totally positive varieties.
We formulate this relationship in terms of two mutually inverse
transformations: ``tropicalization" and ``geometric lifting."
The starting point for this is an observation that different parametrizations of the
canonical basis are related to each other by piecewise-linear transformations that involve
only the operations of addition, subtraction, and taking the minimum of two integers.
As in \cite{bfz}, we shall represent such expressions as ``tropical"
subtraction-free rational expressions.
Recall from \cite{bfz} that a \emph{semifield} $K$ is a set equipped with
operations of addition, multiplication and division (but
no subtraction!) satisfying the usual field axioms.
Two main examples for us will be $\RR_{> 0}$, the set of positive real
numbers with usual operations, and the \emph{tropical semifield}
$\ZZ_{\rm trop}$, the set of integers with
multiplication and addition given by
$$a\, \odot \, b = a + b\ , \quad a \,\oplus \,b =
{\rm min}\, (a,b)\ . $$
We shall write $[Q]_{\rm trop}$ if we need to emphasize that
a subtraction-free rational expression $Q$ is understood in a
tropical sense.
We call $Q$ a \emph{geometric lifting} of a piecewise-linear expression
$[Q]_{\rm trop}$.
Note that a geometric lifting is \emph{not} unique; for example, if $Q$ is a Laurent
polynomial with nonnegative integer coefficients then $[Q]_{\rm trop}$
only depends on the Newton polytope of $Q$ (cf.~\cite[Proposition~4.1.1]{bfz}).

In this terminology, Lusztig's observation can be formulated as follows:
for $\gg$ of simply-laced type, the transition map between any two Lusztig parametrizations
of $\BB$ has a geometric lifting which describes the transition between
two natural parametrizations of the totally positive variety in the maximal unipotent
subgroup of the semisimple group $G$ corresponding to $\gg$.
In this paper we extend this result to a non-simply-laced case;
and we also find a similar geometric lifting for transition maps
between string parametrizations.
This geometric lifting allows us to deduce combinatorial properties of the canonical basis
from geometric properties of totally positive varieties.
To do this, we rely on (and further develop) the
``calculus of generalized minors"  and its applications
to the study of totally positive varieties in Schubert cells and double Bruhat cells
developed in \cite{bfz,bz2,fz,fzoscil}.
An intriguing feature of our results is that combinatorial parametrizations
of the canonical basis and related expressions for tensor product multiplicities for
a semisimple Lie algebra $\gg$ are expressed
in terms of geometry of totally positive varieties in
the Langlands dual semisimple group ${}^L G$.
It would be very interesting to give a conceptual explanation of this phenomenon.

The paper is organized as follows.
Section \ref{sec:LR} provides background on semisimple Lie algebras
and summarizes our main results on tensor product multiplicities
(with the exception of Pl\"ucker models which are given in
Theorems~\ref{th:LR-Plucker-Lusztig} and \ref{th:LR-Plucker-strings}) and
reduction multiplicities.
Section \ref{sec:can} provides background on quantum groups and
canonical bases, and summarizes our results on
parametrizations of canonical bases.
In Section \ref{sec:geometric} we provide needed background, and
present some new results on generalized minors, double Bruhat cells, and
totally positive varieties.
Section \ref{sec:lifting} presents our main results on geometric lifting and
tropicalization; among other things, we find geometric counterpart
of Kashiwara's crystal operators which play an important part in
our arguments.
The remaining sections contain the proofs of all the results in this paper.
The proofs are not very difficult because most of the
needed geometric technique was developed in the preceding papers \cite{bfz,bz2,fz}.
However we need a number of important modifications to the geometric
setup developed in the previous papers: for example, we replace
double Bruhat cells from \cite{fz} by reduced double Bruhat cells introduced in
Section \ref{sec:geometric} below.

\section{Tensor product multiplicities}
\label{sec:LR}

\subsection{Background on semisimple Lie algebras}
\label{sec:background}
We start by fixing the terminology and notation (mostly standard) related
to semisimple Lie algebras.
Throughout the paper, $\gg$ is a complex semisimple Lie algebra of rank~$r$
with Chevalley generators $e_i, \alpha_i^\vee$, and $f_i$ for $i=1, \ldots, r$.
The elements $\alpha_i^\vee$ are \emph{simple coroots} of $\gg$; they form a basis
of a Cartan subalgebra $\hh$ of $\gg$.
The \emph{simple roots} $\alpha_1, \dots, \alpha_r$ form a basis in the dual space
$\hh^*$ such that
$[h, e_i] = \alpha_i (h) e_i$, and $[h,f_i] = - \alpha_i (h) f_i$
for any $h \in \hh$ and $i \in [1,r]$.
The structure of $\gg$ is uniquely determined by the \emph{Cartan matrix}
$A = (a_{ij})$
given by $a_{ij} = \alpha_j (\alpha_i^\vee)$.

The \emph{weight lattice} $P$ of $\gg$ consists of all $\gamma \in \hh^*$ such that
$\gamma (\alpha_i^\vee) \in \ZZ$ for all $i$.
Thus $P$ has a $\ZZ$-basis $\omega_1, \dots, \omega_r$ of \emph{fundamental weights}
given by $\omega_j (\alpha_i^\vee) = \delta_{i,j}$.
A weight $\lambda \in P$ is \emph{dominant} if $\lambda (\alpha_i^\vee) \geq 0$ for
any $i \in [1,r]$; thus $\lambda = \sum_i l_i \omega_i$ with all $l_i$ nonnegative integers.
Let $V_\lambda$ denote the simple (finite-dimensional) $\gg$-module
with highest weight $\lambda$.
We denote by $c_{\lambda, \nu}^\mu$ the multiplicity of the simple module $V_\mu$
in the tensor product $V_\lambda \otimes V_\nu$.

Every finite-dimensional $\gg$-module $V$ is known to be $\hh$-diagonalizable,
and we denote by $V = \oplus_{\gamma \in P}  V(\gamma)$ its weight decomposition.
Let $P(V)$ denote the set of weights of $V$, i.e., the set of all weights $\gamma \in P$
such that $V(\gamma) \neq 0$.

The \emph{Langlands dual} complex semisimple Lie algebra ${}^L \gg$ corresponds to the
transpose Cartan matrix $A^T$.
Thus we can identify the Cartan subalgebra of ${}^L \gg$ with $\hh^*$, and simple roots (resp.
coroots) of ${}^L \gg$ with simple coroots (resp. roots) of $\gg$.
We denote the fundamental weights for ${}^L \gg$ by $\omega_1^\vee, \dots, \omega_r^\vee$;
they are elements of $\hh$ such that $\gamma (\omega_i^\vee)$ is the coefficient
of $\alpha_i$ in the expansion of $\gamma \in \hh^*$ in the basis of simple roots.

The Weyl groups of $\gg$ and $\gg^\vee$ are naturally identified with each other, and we denote
both groups by the same symbol $W$.
As an abstract group, $W$ is a finite Coxeter group generated by \emph{simple reflections}
$s_1, \dots, s_r$; it acts in $\hh^*$ and $\hh$ by
$$s_i (\gamma) = \gamma - \gamma (\alpha_i^\vee) \alpha_i, \,\,
s_i (h) = h - \alpha_i (h) \alpha_i^\vee$$
for $\gamma \in \hh^*$ and $h \in \hh$.

A \emph{reduced word} for $w \in W$ is a sequence of indices
$(i_1, \ldots, i_l)$ that satisfies $w = s_{i_1} \cdots s_{i_l}$
and has the shortest possible length~$l=\l(w)$.
The set of reduced words for~$w$ will be denoted by~$R(w)$.
As customary, $\wnot$ denotes the unique element of maximal length in~$W$.

\subsection{Polyhedral expressions for tensor product multiplicities}
\label{sec:multiplicities}

In this section, we present our first main results: two families of
combinatorial expressions for the tensor product multiplicities.

\begin{definition}
\label{def:trails}
{\rm Let $V$ be a finite-dimensional $\gg$-module,
$\gamma$ and $\delta$ two weights in $P(V)$, and
$\ii = (i_1, \dots, i_l)$ a sequence of indices from $[1,r] := \{1, \dots, r\}$.
An $\ii$-trail from $\gamma$ to $\delta$ in $V$ is a sequence of weights
$\pi = (\gamma = \gamma_0, \gamma_1, \dots , \gamma_l = \delta)$ such that:

\noindent (1) for $k = 1, \dots, l$, we have
$\gamma_{k-1} - \gamma_k = c_k \alpha_{i_k}$ for some nonnegative integer $c_k$;

\noindent (2) $e_{i_1}^{c_1} \cdots e_{i_l}^{c_l}$ is a non-zero linear map
from $V(\delta)$ to $V(\gamma)$.}
\end{definition}

Note that the numbers $c_k$ in Definition~\ref{def:trails} can be written as
\begin{equation}
\label{eq:ck}
c_k = c_k (\pi) = \frac{\gamma_{k-1} - \gamma_k}{2} (\alpha_{i_k}^\vee) \ .
\end{equation}

Our expressions for the tensor product multiplicity $c_{\lambda, \nu}^\mu$
for $\gg$ will involve $\ii$-trails in the fundamental modules $V_{\omega_i^\vee}$
of the Langlands dual Lie algebra ${}^L \gg$.
We denote by $m$ the length of $\wnot$.

\begin{theorem}
\label{th:LR Lusztig-trails}
Let $\lambda, \mu, \nu$ be three dominant weights for $\gg$.
For any reduced word $\ii = (i_1, \dots, i_m)  \in R(\wnot)$, the multiplicity
$c_{\lambda, \nu}^\mu$ is equal to the number of integer $m$-tuples $(t_1, \dots, t_m)$
satisfying the following conditions:

\smallskip

\noindent {\rm (1)} $t_k \geq 0$ for any $k = 1, \dots, m$;

\smallskip

\noindent {\rm (2)} $\sum_k t_k \cdot s_{i_1} \cdots s_{i_{k-1}} \alpha_{i_k}=
\lambda + \nu - \mu$;

\smallskip

\noindent {\rm (3)} $\sum_k c_k (\pi) t_k \geq  (s_i \lambda + \nu - \mu)(\omega_i^\vee)$
for any $i$ and any $\ii$-trail $\pi$ from $s_i \omega_i^\vee$ to $\wnot \omega_i^\vee$
in $V_{\omega_i^\vee}$;

\smallskip

\noindent {\rm (4)} $\sum_k c_{k}(\pi) t_k \geq  (\lambda + s_i \nu - \mu)(\omega_i^\vee)$
for any $i$ and any $\ii$-trail $\pi$ from $\omega_i^\vee$ to $\wnot s_i \omega_i^\vee$
in $V_{\omega_i^\vee}$.
\end{theorem}

To present our second family of polyhedral expressions, let us define,
for every $\ii$-trail $\pi = (\gamma_0, \dots, \gamma_l)$ in a $\gg$-module $V$,
and every $k = 1, \dots, l$ (cf.~(\ref{eq:ck})):
\begin{equation}
\label{eq:dk}
d_k = d_k (\pi) = \frac{\gamma_{k-1} + \gamma_k}{2} (\alpha_{i_k}^\vee) \ .
\end{equation}

\begin{theorem}
\label{th:LR string-trails}
Let $\lambda, \mu, \nu$ be three dominant weights for $\gg$.
For any reduced word $\ii = (i_1, \dots, i_m)  \in R(\wnot)$, the multiplicity
$c_{\lambda, \nu}^\mu$ is equal to the number of integer $m$-tuples $(t_1, \dots, t_m)$
satisfying the following conditions:

\smallskip

\noindent {\rm (1)} $\sum_k d_k (\pi) t_k \geq 0$ for any $i$ and any
$\ii$-trail $\pi$ from $\omega_i^\vee$ to $\wnot s_i \omega_i^\vee$
in $V_{\omega_i^\vee}$;

\smallskip

\noindent {\rm (2)} $\sum_k t_k \alpha_{i_k} = \lambda + \nu - \mu$;

\smallskip

\noindent {\rm (3)} $\sum_k d_{k} (\pi) t_k \geq - \lambda (\alpha_{i}^\vee)$
for any $i$ and any $\ii$-trail $\pi$ from $s_i \omega_i^\vee$ to
$\wnot \omega_i^\vee$ in $V_{\omega_i^\vee}$;

\smallskip

\noindent {\rm (4)} $t_k + \sum_{l > k} a_{i_k, i_l} t_l  \leq  \nu (\alpha_{i_k}^\vee)$
for any $k = 1, \dots, m$.
\end{theorem}

Explicit bijections between the $m$-tuples in Theorems~\ref{th:LR Lusztig-trails}
and \ref{th:LR string-trails} will be given in Theorem~\ref{th:Lusztig-string} below.

\subsection{Tensor product multiplicities for classical groups}
\label{sec:LR BCD}
In this section we present most concrete expressions for the multiplicity
$c_{\lambda, \nu}^\mu$ in the case when $\gg$ is one of the classical simple Lie algebras of
types $B,C$ and $D$.
We start with the types $B_r$ ($\gg = so_{2r+1}$) and
$C_r$ ($\gg = sp_{2r}$).
Choose a (non-standard) numeration of simple roots of
$\gg$ by the index set $[0, r-1]$ so that the roots
$\alpha_1, \dots, \alpha_{r-1}$ form a system of type $A_{r-1}$
in the standard numeration, and
$a_{01} = -2, \, a_{10} = -1$ for type $B_r$, and
$a_{01} = -1, \, a_{10} = -2$ for type $C_r$.
Let us denote $a = |a_{10}|$, i.e., $a = 1$ for $\gg = so_{2r+1}$,
and $a = 2$ for $\gg = sp_{2r}$.

\begin{theorem}
\label{th:LR string-trails B}
Let $\gg$ be a simple Lie algebra of type $B_r$ or $C_r$, and let
$\lambda, \mu, \nu$ be three dominant weights for $\gg$.
Then the multiplicity $c_{\lambda, \nu}^\mu$ is equal to the number of integer
tuples $(t_i^{(j)}: 0 \leq |i| \leq j < r)$
satisfying the following conditions:

\smallskip

\noindent {\rm (1)} $2 t_{-j}^{(j)} \geq \cdots \geq 2 t_{-1}^{(j)}\geq
a t_{0}^{(j)} \geq 2 t_{1}^{(j)}\geq \cdots \geq  2 t_{j}^{(j)}\geq 0$
for $0 \leq j < r$;

\smallskip

\noindent {\rm (2)} $\sum_{0 \leq |i| \leq j < r} t_i^{(j)} \alpha_{|i|} = \lambda + \nu - \mu$;

\smallskip

\noindent {\rm (3)} $\lambda (\alpha_{0}^\vee)
\geq t_0^{(0)}$, and $\lambda (\alpha_{j}^\vee) \geq \max \
\{t_j^{(j)}, a t_0^{(j)} - t_1^{(j-1)} - t_{-1}^{(j)},
t_1^{(j-1)} + t_{-1}^{(j)} - a t_0^{(j-1)}, \varphi_i^{(j)} (t) \,
(1 \leq i < j)\}$ for $0 \leq j < r$, where $\varphi_i^{(j)} (t) =
\max \{t_i^{(j)} + t_{-i}^{(j)} - t_{i+1}^{(j-1)} -
t_{-i-1}^{(j)}, t_{i+1}^{(j-1)} + t_{-i-1}^{(j)} - t_{-i}^{(j-1)} -
t_{i}^{(j-1)}, t_{\pm i}^{(j)} - t_{\pm i}^{(j-1)}\}$ (with the
convention that $t_{i}^{(j)} = 0$ unless $0 \leq |i| \leq j < r$);

\smallskip

\noindent {\rm (4)} $\nu (\alpha_{0}^\vee) \geq \max_{j \geq 0} \
(t_0^{(j)} + \frac{2}{a} \sum_{k > j} (a t_0^{(k)} - t_{-1}^{(k)} -
t_1^{(k-1)}))$,
$$\nu (\alpha_{1}^\vee) \geq \max_{j \geq 1} \ \max
(t_{-1}^{(j)} + \sum_{k > j} (2t_{-1}^{(k)} + 2t_1^{(k-1)} - a t_0^{(k-1)}
- t_{-2}^{(k)} - t_2^{(k-1)}),$$
$$t_{1}^{(j)} + \sum_{k > j} (2t_{-1}^{(k)} + 2t_1^{(k)} - a t_0^{(k)}
- t_{-2}^{(k)} - t_2^{(k-1)})) \ ,$$
$$\nu (\alpha_{i}^\vee) \geq \max_{j \geq i} \ \max
(t_{-i}^{(j)} + \sum_{k > j} (2t_{-i}^{(k)} + 2t_i^{(k-1)} - t_{-i+1}^{(k-1)}
- t_{i-1}^{(k-1)} - t_{-i-1}^{(k)} - t_{i+1}^{(k-1)}),$$
$$t_{i}^{(j)} + \sum_{k > j} (2t_{-i}^{(k)} + 2t_i^{(k)} -
t_{-i+1}^{(k)} - t_{i-1}^{(k)} - t_{-i-1}^{(k)} - t_{i+1}^{(k-1)}))$$
for $2 \leq i < r$.
\end{theorem}

Now consider the type $D_r$ ($\gg = so_{2r}$).
Choose a (non-standard) numeration of simple roots of
$\gg$ by the index set $\{-1\} \cup [1, r-1]$ so that the roots
$\alpha_1, \dots, \alpha_{r-1}$ form a system of type $A_{r-1}$
in the standard numeration, and
$a_{-1,2} = a_{2,-1} = -1$.

\begin{theorem}
\label{th:LR string-trails D}
Let $\gg$ be a simple Lie algebra of type $D_r$ and let
$\lambda, \mu, \nu$ be three dominant weights for $\gg$.
Then the multiplicity $c_{\lambda, \nu}^\mu$ is equal to the number of integer
tuples $(t_i^{(j)}: 1 \leq |i| \leq j < r)$
satisfying the following conditions:

\smallskip

\noindent {\rm (1)} $t_{-j}^{(j)} \geq \cdots \geq t_{-2}^{(j)}\geq
\max \ (t_{-1}^{(j)}, t_{1}^{(j)}) \geq \min \ (t_{-1}^{(j)}, t_{1}^{(j)})
\geq t_{2}^{(j)}\geq \cdots \geq  t_{j}^{(j)}\geq 0$
for $1 \leq j < r$;

\smallskip

\noindent {\rm (2)} $\sum_j (t_{-1}^{(j)} \alpha_{-1}  + t_{1}^{(j)}\alpha_1)
+ \sum_{2 \leq |i| \leq j < r} t_i^{(j)} \alpha_{|i|} = \lambda + \nu - \mu$;

\smallskip

\noindent {\rm (3)} $\lambda (\alpha_{\pm 1}^\vee) \geq t_{\pm 1}^{(1)}$,
and $\lambda (\alpha_{j}^\vee) \geq \max \
\{t_j^{(j)}, t_{\pm 1}^{(j)} - t_{\mp 1}^{(j-1)},
t_1^{(j)} + t_{-1}^{(j)} - t_{-2}^{(j)} - t_2^{(j-1)},
t_{-2}^{(j)} + t_{2}^{(j-1)} - t_1^{(j-1)} - t_{-1}^{(j-1)},
\varphi_i^{(j)} (t) \, (2 \leq i < j)\}$
for $2 \leq j < r$, where $\varphi_i^{(j)} (t)$
is the same as in Theorem~\ref{th:LR string-trails B};

\smallskip

\noindent {\rm (4)} $\nu (\alpha_{\pm 1}^\vee) \geq \max_{j \geq 1} \
(t_{\pm 1}^{(j)} + \sum_{k > j} (2 t_{\pm 1}^{(k)} - t_{-2}^{(k)} -
t_2^{(k-1)}))$,
$$\nu (\alpha_{i}^\vee) \geq \max_{j \geq i} \ \max
(t_{-i}^{(j)} + \sum_{k > j} (2t_{-i}^{(k)} + 2t_i^{(k-1)} - t_{-i+1}^{(k-1)}
- t_{i-1}^{(k-1)} - t_{-i-1}^{(k)} - t_{i+1}^{(k-1)}),$$
$$t_{i}^{(j)} + \sum_{k > j} (2t_{-i}^{(k)} + 2t_i^{(k)} -
t_{-i+1}^{(k)} - t_{i-1}^{(k)} - t_{-i-1}^{(k)} - t_{i+1}^{(k-1)}))$$
for $2 \leq i < r$.
\end{theorem}

\begin{remark}
{\rm We shall show that Theorems~\ref{th:LR string-trails B} and
\ref{th:LR string-trails D} can be obtained by specializing Theorem~\ref{th:LR string-trails}
for a specific reduced word $\ii$.
One can also show that specializing Theorem~\ref{th:LR Lusztig-trails}
for the same $\ii$ leads to the expressions for $c_{\lambda, \nu}^\mu$
that were conjectured in \cite{bz88}.}
\end{remark}

\subsection{Reduction multiplicities}
\label{sec:reduction}
For a subset $I \subset [1,r]$, let $\gg(I)$ denote the corresponding
\emph{Levi subalgebra} in $\gg$  generated by the Cartan subalgebra
$\hh$ and by all $e_i$ and $f_i$ for $i \in I$.
A weight $\beta \in P$ is dominant for $\gg (I)$ if
$\beta (\alpha_i^\vee) \geq 0$ for $i \in I$; for such a weight,
let $V_\beta^{(I)}$ denote the simple (finite-dimensional)
$\gg(I)$-module with highest weight $\beta$.
In this section, we compute the multiplicity of $V_\beta^{(I)}$
in the reduction to $\gg (I)$ of a simple $\gg$-module $V_\nu$.
This includes the weight multiplicities in $V_\nu$ as a special case
when $I = \emptyset$.

Let $\wnot (I)$ denote the longest element
of the parabolic subgroup in $W$ generated by all $s_i$ with $i \in I$.

\begin{theorem}
\label{th:reduction multiplicity Lusztig}
For any reduced word $\ii = (i_1, \dots, i_n)  \in R(\wnot (I)^{-1} \wnot)$,
the multiplicity of $V_\beta^{(I)}$ in the reduction of
$V_\nu$ to $\gg (I)$ is equal to the number of integer $n$-tuples $(t_1, \dots, t_n)$
satisfying the following conditions:

\smallskip

\noindent {\rm (1)} $t_k \geq 0$ for $k = 1, \dots, n$;

\smallskip

\noindent {\rm (2)} $\sum_k t_k \cdot s_{i_1} \cdots s_{i_{k-1}} \alpha_{i_k}=
\wnot (I) (\nu - \beta)$;

\smallskip

\noindent {\rm (3)} $\sum_k c_k (\pi) t_k \geq  (s_i \nu - \beta)(\omega_i^\vee)$
for any $i$ and any $\ii$-trail $\pi$ from $\wnot (I) \omega_i^\vee$ to $\wnot s_i \omega_i^\vee$
in $V_{\omega_i^\vee}$.
\end{theorem}

\begin{theorem}
\label{th:reduction multiplicity string}
For any reduced word $\ii = (i_1, \dots, i_n)  \in R(\wnot (I)^{-1} \wnot)$,
the multiplicity of $V_\beta^{(I)}$ in the reduction of
$V_\nu$ to $\gg (I)$ is equal to the number of integer $n$-tuples $(t_1, \dots, t_n)$
satisfying the following conditions:

\smallskip

\noindent {\rm (1)} $\sum_k d_k (\pi) t_k \geq 0$ for any $i$ and any
$\ii$-trail $\pi$ from $\wnot (I) \omega_i^\vee$ to $\wnot s_i \omega_i^\vee$
in $V_{\omega_i^\vee}$;

\smallskip

\noindent {\rm (2)} $\sum_k t_k \alpha_{i_k} = \nu - \beta$;

\smallskip

\noindent {\rm (3)} $t_k + \sum_{l > k} a_{i_k, i_l} t_l  \leq
\nu (\alpha_{i_k}^\vee)$ for $k = 1, \dots, n$.
\end{theorem}

Explicit bijections between the $n$-tuples in
Theorems~\ref{th:reduction multiplicity Lusztig}
and \ref{th:reduction multiplicity string} will be given in
Theorem~\ref{th:Lusztig-string} below.
We illustrate these theorems with the following example.

\begin{example}
\label{ex:pq reduction}
{\rm Let $\gg=sl_{r+1}$ be of type $A_r$ (with the standard numeration of simple roots).
Let $I=[1,r]\setminus \{p\}$ for some $p \in [1,r]$.
Let $q = r+1 - p$; then the algebra $\gg (I)$ is the intersection
of $sl_{r+1} = sl_{p+q}$ with the block-diagonal subalgebra
$gl_{p} \times gl_{p} \subset gl_{p+q}$.
Denote by $M_{p \times q}$ the set of all $p \times q$ matrices
$T=(t_{ij})$ with integer entries (we shall also use the
convention that $t_{ij} = 0$ unless $(i,j) \in [1,p] \times [1,q]$).
Theorems~\ref{th:reduction multiplicity Lusztig} and
\ref{th:reduction multiplicity string} specialize to the following
two expressions for the reduction multiplicity.

\begin{corollary}
\label{cor:pq reduction Lusztig}
The multiplicity of $V_\beta^{(I)}$ in the reduction of
$V_\nu$ to $\gg (I)$ is equal to the number of all $T \in M_{p \times q}$ satisfying
the following conditions:

\smallskip

\noindent {\rm (1)} $t_{ij} \geq 0$ for all $i$ and $j$;

\smallskip

\noindent {\rm (2)} $\sum_{i \leq l; j \leq p+q-l} t_{ij} = (\nu - \beta)(\omega_l^\vee)$
for any $l \in [1,r]$;

\smallskip

\noindent {\rm (3)} $\nu(\alpha_i^\vee) \geq \max_{j \in [1,q]}
(\sum_{k \geq 0} (t_{i,j+k}-t_{i+1,j+k+1}))$ for $1 \leq i < p$;
\, \,  $\nu(\alpha_p^\vee) \geq t_{pq}$; \,\, $\nu(\alpha_{p+q-j}^\vee) \geq \max_{i \in [1,p]}
(\sum_{k \geq 0} (t_{i+k,j}-t_{i+k+1,j+1}))$ for $1 \leq j < q$.
\end{corollary}

\begin{corollary}
\label{cor:pq reduction string}
The multiplicity of $V_\beta^{(I)}$ in the reduction of
$V_\nu$ to $\gg (I)$ is equal to the number of all $T \in M_{p \times q}$ satisfying
the following conditions:

\smallskip

\noindent {\rm (1)} $t_{ij} \geq \max (t_{i+1,j}, t_{i,j+1})$ for all $i$ and $j$
(that is, $T$ is a plane partition of shape $p \times q$);

\smallskip

\noindent {\rm (2)} $\sum_{j-i=l-p} t_{ij} = (\nu - \beta)(\omega_l^\vee)$
for any $l \in [1,r]$;

\smallskip

\noindent {\rm (3)} $\nu(\alpha_l^\vee) \geq \max_{j-i=l-p}
(t_{ij} + \sum_{k > 0} (2 t_{i+k,j+k}-t_{i+k-1,j+k} - -t_{i+k,j+k-1}))$
for any $l \in [1,r]$.
\end{corollary}
}
\end{example}

\begin{remark}
{\rm The polytope defined by conditions (1) - (3) in
Corollary~\ref{cor:pq reduction Lusztig} (with the additional
assumption that $\nu(\alpha_p^\vee) >> 0$ so that the inequality
$\nu(\alpha_p^\vee) \geq t_{pq}$ in (3) is skipped) appeared in a
different context in \cite{donin}. Comparing Corollary~\ref{cor:pq
reduction Lusztig} with \cite[Theorem~1]{donin}, we conclude
the following.
If $\nu(\alpha_p^\vee) >> 0$ then the reduction multiplicity in
Corollary~\ref{cor:pq reduction Lusztig} is equal to the inner product
of two skew Schur functions $s_{\nu^{(1)}/\beta^{(1)}}$ and
$s_{\nu^{(2)}/\beta^{(2)}}$, where
for each dominant $\gg(I)$-weight $\mu$ we define the partition
$\mu^{(1)}$ (resp. $\mu^{(2)}$) with $\leq p$ (resp. $\leq q$)
parts by setting $\mu^{(1)}_i = \mu (\omega_i^\vee -
\omega_{i-1}^\vee)$ (resp. $\mu^{(2)}_j = \mu (\omega_{p+q-j}^\vee -
\omega_{p+q-j+1}^\vee)$).
Note that a polyhedral expression for the inner product of any
skew Schur functions was first obtained in \cite{bz0}; the
expression given there is close to the one in
Corollary~\ref{cor:pq reduction string}.
A bijective correspondence between the tuples in
Corollaries~\ref{cor:pq reduction Lusztig} and
\ref{cor:pq reduction string} was constructed in an unpublished work of
one of the authors (A.B.) and Anatol Kirillov; they also observed
that their bijection can be interpreted as the
Robinson-Schensted-Knuth correspondence.}
\end{remark}

\section{Canonical bases and their parametrizations}
\label{sec:can}

\subsection{Background on canonical bases and their Lusztig parametrizations}
\label{sec:canonical bases}
Let us recall some basic facts about quantized universal enveloping algebras and
their canonical bases.
Unless otherwise stated, all the results in this section are due to G.~Lusztig and
can be found in \cite{lu}.
The quantized universal enveloping
algebra $U=U_q (\gg)$ associated to $\gg$ is defined as follows.
Fix positive integers $d_1, \dots, d_r$ such that $d_i a_{ij}=d_j a_{ji}$,
where $(a_{ij})$ is the Cartan matrix of $\gg$.
The algebra $U$ is a $\CC(q)$-algebra with unit generated by the elements
$E_i, K_i^{\pm 1}$, and $F_i$ for $i = 1, \dots, r$
subject to the relations
$$K_i K_j = K_j K_i, \,\, K_i E_j K_i^{-1} = q^{d_i a_{ij}}E_j, \,\,
K_iF_jK_i^{-1}=q^{-d_ia_{ij}}E_j \ ,$$
$$E_i F_j - F_j E_i= \delta_{ij} \frac{K_i-K_i^{-1}}{q^{d_i}-q^{-d_i}}$$
for all $i$ and $j$, and the \emph{quantum Serre relations}
$$\sum_{k+l=1-a_{ij}} (-1)^k  E_i^{(k)} E_j E_i^{(l)} = \sum_{k+l=1-a_{ij}}
(-1)^k F_i^{(k)}F_jF_i^{(l)} = 0$$
for $i\ne j$.
Here $E_i^{(k)}$ and $F_i^{(k)}$ stand for the \emph{divided powers}
defined by
$$E_i^{(k)}=\frac{1}{[1]_i [2]_i \cdots [k]_i}E_i^k, \,\,
F_i^{(k)}=\frac{1}{[1]_i [2]_i \cdots [k]_i} F_i^k \ ,$$
where $[l]_i=\frac{q^{d_i l}-q^{-d_i l}}{q^{d_i }-q^{-d_i }}$.
The algebra $U$ is graded by the root lattice of $\gg$ via
$${\rm deg}(K_i)=0, \,\, {\rm deg}(E_i)= -{\rm deg}(F_i)= \alpha_i \ .$$

To each $i = 1, \dots, r$, Lusztig associates an
algebra automorphism $T_i$ of $U$ uniquely determined by:
$$T_i (K_j)= K_j K_i^{- a_{ij}} \quad (j = 1, \dots, r) \ ,$$
$$T_i (E_i)= -K_i^{-1} F_i, \,\, T_i(F_i)= -E_i K_i \ ,$$
and, for all $j\neq i$,
$$T_i (E_j)= \sum_{k+l= -a_{ij}} (-1)^k q^{- d_i k} E_i^{(k)}E_j E_i^{(l)}
\ ,$$
$$T_i (F_j)=\sum_{k+l=-a_{ij}} (-1)^k q^{d_i k} F_i^{(l)}F_jF_i^{(k)} \,\,
\ .$$
(This automorphism was denoted by $T'_{i,-1}$ in \cite{lu}.)
The $T_i$ satisfy the braid relations
and so extend to an action of the braid group on $U$.

Let $U^+$ denote the subalgebra of $U$ generated by $E_1, \dots, E_r$.
We now recall Lusztig's definitions
of the PBW-type bases and the canonical basis in $U^+$.
For a reduced word $\ii=(i_1,\ldots,i_m) \in R(\wnot)$, and
an $m$-tuple $t = (t_1, \dots, t_m) \in \ZZ^m_{\ge 0}$, denote
$$p_\ii^{(t)} := E_{i_1}^{(t_1)} T_{i_1}(E_{i_2}^{(t_2)})
\cdots (T_{i_1} \cdots T_{i_{m-1}})(E_{i_m}^{(t_m)}) \ .$$
As shown in \cite{lu}, all these elements belong to $U^+$.
For a given $\ii$, the set of all $p_\ii^{(t)}$ with $t \in \ZZ^m_{\ge 0}$ is called
the \emph{PBW type basis} corresponding to $\ii$ and is denoted by
$\BB_\ii$.
This terminology is justified by the following proposition
proved in \cite[Corollary~40.2.2]{lu}.

\begin{proposition}
For every $\ii \in R(\wnot)$, the set $\BB_\ii$ is a $\CC(q)$-basis of $U_+$.
\end{proposition}

According to \cite[Proposition~8.2]{lu96}, the canonical basis $\BB$ of $U^+$ can
now be defined as follows.
Let $u \mapsto \overline u$ denote the
$\CC$-linear involutive algebra automorphism of $U^+$ such that
$\overline {q} = q^{-1}, \,\, \overline {E_i} = E_i$.

\begin{proposition}
\label{pr:characterization of BB}
For every $\ii \in R(\wnot)$ and $t \in  \ZZ_{\ge 0}^m$,
there is a unique element $b = b_\ii (t)$ of $U^+$ such that
$\overline {b} = b$, and $b - p_\ii^{(t)}$ is a linear combination
of the elements of $\BB_\ii$ with coefficients in $q^{-1} \ZZ[q^{-1}]$.
For any fixed $\ii$, the elements
$b_{\ii}(t)$ for all
$t \in  \ZZ_{\ge 0}^m$ constitute the canonical basis $\BB$.
\end{proposition}

In view of Proposition~\ref{pr:characterization of BB},
any $\ii \in R(\wnot)$ gives rise to a bijection
$b_\ii: \ZZ_{\ge 0}^m \to \BB$.
We refer to these bijections as \emph{Lusztig parametrizations} of $\BB$.
Let us summarize some of their properties.
To do this, we need some more notation.
Let $i \mapsto i^*$ denote the involution on $[1,r]$
defined by $\wnot (\alpha_i) = - \alpha_{i^*}$.
For every sequence $\ii = (i_1, \dots, i_m)$,
we denote by $\ii^*$ and $\ii^{\rm op}$ the sequences
\begin{equation}
\label{eq:* and op}
\ii^* = (i_1^*, \dots, i_m^*), \,\, \ii^{\rm op} = (i_m, \dots, i_1) \ ;
\end{equation}
clearly, both operations $\ii \mapsto \ii^*$ and $\ii \mapsto \ii^{\rm op}$
preserve the set of reduced words $R(\wnot)$.

\begin{proposition}
\label{pr:Lusztig parametrizations}
{\rm (i)} Any canonical basis vector $b_\ii (t) \in \BB$ is homogeneous of degree
$\sum_k t_k \cdot s_{i_1} \cdots s_{i_{k-1}} \alpha_{i_k}$.

\noindent {\rm (ii)} Every subspace of the form $E_i^n U^+$ in $U^+$ is spanned by
a subset of $\BB$.
Furthermore, let $l_i (b)$ denote the maximal integer $n$ such that $b \in E_i^n U^+$;
then, for any $\ii \in R(\wnot)$ which begins with $i_1 = i$, we have
$l_i (b_\ii (t_1, \dots, t_m)) = t_1$.

\noindent {\rm (iii)} The canonical basis $\BB$ is stable under
the involutive $\CC(q)$-linear algebra antiautomorphism $E \to E^\iota$ of
$U^+$ such that $E_i^\iota = E_i$ for all $i$.
Furthermore, we have $b_\ii (t)^\iota = b_{\ii^{* {\rm op}}} (t^{\rm op})$.
\end{proposition}

This proposition allows us to interpret the tensor product multiplicity $c_{\lambda, \nu}^\mu$
in terms of the canonical basis.
Indeed it is well-known (see e.g., \cite{bz1}) that $c_{\lambda, \nu}^\mu$ is equal to
the dimension of the homogeneous component of degree $\lambda + \nu - \mu$ in
$U^+/ \sum_i (E_i^{\lambda (\alpha_i^\vee) + 1} U^+ \ + \ U^+ E_i^{\nu (\alpha_i^\vee) + 1})$.
Thus Proposition~\ref{pr:Lusztig parametrizations} has the following corollary.

\begin{corollary}
\label{cor:multiplicity through BB}
The multiplicity $c_{\lambda, \nu}^\mu$ is equal to the number of vectors
$b \in \BB$ of degree $\lambda + \nu - \mu$ satisfying the following property:
if $b = b_\ii (t_1, \dots, t_m)$, and $\ii \in R(\wnot)$ begins with $i$ and ends with $j$ then
$t_1 \leq \lambda (\alpha_i^\vee)$ and $t_m \leq \nu (\alpha_{j^*}^\vee)$.
\end{corollary}

\subsection{Preliminaries on string parametrizations}
\label{sec:string parametrizations}
We now describe another way to parametrize $\BB$ with certain strings
of nonnegative integers, the so-called \emph{string parametrizations}
introduced in \cite{bz93,k93}.
As a general setup, consider any $U^+$-module $V$ such that
each $E_i$ acts on $V$ as a locally nilpotent operator,
i.e., for every non-zero $v \in V$ there exists some positive integer $n$
such that $E_i^n (v) = 0$.
Let $c_i (v)$ denote the maximal integer $n$ such that
$E_i^n (v) \neq 0$.
For any sequence of indices $\ii = (i_1, \dots, i_m)$,
the \emph{string} of $v$ in direction $\ii$ is a sequence
$c_\ii (v) = (t_1, \dots, t_m)$ of nonnegative integers
defined recursively as follows:
$$t_1 = c_{i_1} (v), t_2 = c_{i_2} (E_{i_1}^{t_1} (v)), \dots,
t_m = c_{i_m} (E_{i_{m-1}}^{t_{m-1}} \cdots E_{i_1}^{t_1} (v)) \ .$$

We apply this construction to a $U^+$-module $\AA$
defined as follows.
As a $\CC(q)$-vector space, $\AA$ is the restricted dual vector space
of $U^+$, i.e., the direct sum of dual spaces of all homogeneous
components of $U^+$ (recall that $U^+$ is graded by
the root lattice of $\gg$).
The action $(E,f) \mapsto E(f)$ of $U^+$ on $\AA$ is given by
$E(f) (u) = f(E^\iota u)$ for all $u \in U^+$
(see Proposition~\ref{pr:Lusztig parametrizations} (iii)).
Clearly, each $E_i$ acts in $\AA$ as a locally nilpotent operator,
so the corresponding strings are well defined.

Now we consider the \emph{dual canonical} basis $\BB^{\rm dual}$ in
$\AA$: its element $b^{\rm dual}$ corresponding to $b \in \BB$
is a linear form on $U^+$ such that $b^{\rm dual} (b') = \delta_{b,b'}$
for $b' \in \BB$.
The following proposition was essentially proved in \cite{lit}
(our results below will provide an independent proof).

\begin{proposition}
\label{pr:strings}
For any $\ii \in R(\wnot)$, the string parametrization $c_\ii$
is a bijection of $\BB^{\rm dual}$ onto the set of all lattice points $C_\ii(\ZZ)$
of some rational polyhedral convex cone $C_\ii$ in $\RR^m$.
\end{proposition}

We call the cone $C_\ii$ in Proposition~\ref{pr:strings} the
\emph{string cone} associated to $\ii \in R(\wnot)$
(it was called the cone of adapted strings in \cite{lit}).

Comparing the definition of strings with the definition of $l_i (b)$ in
Proposition~\ref{pr:Lusztig parametrizations} (ii), we see that
\begin{equation}
\label{eq:li = ci}
l_i (b) = c_i (b^{\rm dual})
\end{equation}
for any $i \in [1,r]$ and $b \in \BB$.
We shall also need some basic properties of $\BB$ related to Kashiwara's
crystal structure.
The following proposition is a consequence of results in \cite{k93} and
\cite{lu} (see also \cite{nakzel97}).

\begin{proposition}
\label{pr:crystal BB}
{\rm (i)} There exist embeddings $\tilde f_i: \BB \to \BB$ for $i = 1, \dots, r$
which satisfy the following property: if $b = b_\ii (t_1, \dots, t_m)$ and
$c_\ii (b^{\rm dual}) = (t'_1, \dots, t'_m)$ for some
$\ii = (i_1, \dots, i_m) \in R(\wnot)$ with $i_1 = i$ then
$\tilde f_i (b) = b_\ii (t_1 + 1, t_2, \dots, t_m)$, and
$c_\ii (\tilde f_i (b)^{\rm dual}) = (t'_1 + 1, t'_2, \dots, t'_m)$.

\noindent {\rm (ii)} Every $b \in \BB$ has the form
$b = \tilde f_{i_1} \cdots \tilde f_{i_N} (1)$
for some sequence of indices $i_1, \dots, i_N$ from $[1,r]$.
\end{proposition}

\subsection{New results on Lustzig and string parametrizations}
\label{sec:new stuff on can}
We start with an explicit formula for the relationship between Lusztig
and string parametrizations.

\begin{theorem}
\label{th:Lusztig-string}
{\rm (i)} Let $\ii$ and $\ii'$ be two reduced words of $\wnot$, and
let $t = c_{\ii} (b_{\ii'} (t')^{\rm dual})$ be the string in direction $\ii$
of the dual canonical basis vector with the Lusztig parameters $t'$ relative to
$\ii'$.
Then $t$ and $t'$ are related as follows: for any $k = 1, \dots, m$, we have
\begin{equation}
\label{eq:Lusztig-string}
t_k =
\min_{\pi_1} \ (\sum_{l=1}^m c_{l}(\pi_1) \cdot t'_l) -
\min_{\pi_2} \ (\sum_{l=1}^m c_{l}(\pi_2) \cdot t'_l) \ ,
\end{equation}
where $\pi_1$ (resp. $\pi_2$) runs over $\ii$-trails from
$s_{i_1} \cdots  s_{i_{k-1}} \omega_{i_k}^\vee$
(resp. from $s_{i_1} \cdots  s_{i_k} \omega_{i_k}^\vee$)
to $\wnot \omega_{i_k}^\vee$ in $V_{\omega_{i_k}^\vee}$.

\noindent {\rm (ii)} For any three dominant weights $\lambda$, $\mu$ and $\nu$,
the correspondence {\rm (\ref{eq:Lusztig-string})} restricts to a bijection
between the set of tuples $t$ in Theorem~{\rm \ref{th:LR string-trails}}
and the set of tuples $t'$ in Theorem~{\rm \ref{th:LR Lusztig-trails}} with
$\ii$ replaced by $\ii'$.
\end{theorem}

Using Theorem~\ref{th:Lusztig-string}, we shall obtain the following
explicit expression for the functions $l_i$ on $\BB$ in terms of Lusztig parameters
(see Proposition~\ref{pr:Lusztig parametrizations} (ii)).

\begin{theorem}
\label{th:li}
For every $i \in [1,r]$ and $\ii = (i_1, \dots, i_m) \in R(\wnot)$, we have
\begin{equation}
\label{eq:li}
l_i (b_\ii (t_1, \dots, t_m)) =
\max_{\pi} \sum_k (s_{i_1} \cdots s_{i_{k-1}} \alpha_{i_k} (\omega_i^\vee) - c_k (\pi)) t_k \ ,
\end{equation}
where $\pi$ runs over all $\ii$-trails from $s_i \omega_i^\vee$ to $\wnot \omega_i^\vee$
in $V_{\omega_i^\vee}$.
\end{theorem}

The inverse of the map in (\ref{eq:Lusztig-string}) which expresses Lusztig parameters of $b$ in
terms of strings of $b^{\rm dual}$ can be also computed explicitly
but this expression is more involved.
We shall only present the following two important special cases.

\begin{theorem}
\label{th:1st and last Lusztig-through-string}
Let $t$ and $t'$ have the same meaning as in Theorem~\ref{th:Lusztig-string} (i).
Then
\begin{equation}
\label{eq:1st Lusztig}
t'_1 = - \min_{\pi} \sum_k d_{k} (\pi) t_k \ ,
\end{equation}
where $\pi$ runs over all $\ii$-trails from $s_{i'_1} \omega_{i'_1}^\vee$ to
$\wnot \omega_{i'_1}^\vee$
in $V_{\omega_{i'_1}^\vee}$; and also
\begin{equation}
\label{eq:last Lusztig}
t'_m = \max_{k: i_k^* = i'_m} (t_k +  \sum_{l > k} a_{i'_m, i_l^*} t_l) \ .
\end{equation}
\end{theorem}

We conclude this section by an explicit description of the string cones
(see Proposition~\ref{pr:strings}).

\begin{theorem}
\label{th:string cones general}
For any reduced word $\ii \in R(\wnot)$, the string cone $C_\ii$  is the cone in
$\RR^m$ given by the inequalities {\rm (1)} in
Theorem~\ref{th:LR string-trails}, i.e., it consists of all real
$m$-tuples $(t_1, \dots, t_m)$ such that
$\sum_k d_k (\pi) t_k \geq 0$ for any $i$ and any
$\ii$-trail $\pi$ from $\omega_i^\vee$ to $\wnot s_i \omega_i^\vee$
in $V_{\omega_i^\vee}$.
\end{theorem}

\subsection{More on string cones}
\label{sec:string cones special}
In this section we describe some specializations of Theorem~\ref{th:string cones general}.
Our first result shows that, under some conditions on a reduced word $\ii \in R(\wnot)$,
the corresponding string cone $C_\ii$ splits into the direct product of smaller cones.
To formulate the result, we need some more notation.
First, for any sequence $\ii = (i_1, \dots, i_l)$ of indices from $[1,r]$, and any two elements
$u, v \in W$, let $C_\ii (u,v)$ denote the cone of all real
$l$-tuples $(t_1, \dots, t_l)$ such that
$\sum_k d_k (\pi) t_k \geq 0$ for any $i$ and any
$\ii$-trail $\pi$ from $u \omega_i^\vee$ to $v s_i \omega_i^\vee$
in $V_{\omega_i^\vee}$.
(In particular, Theorem~\ref{th:string cones general} claims that
$C_\ii = C_\ii (e, \wnot)$.)
Second, as in Section~\ref{sec:reduction}, for any subset $I \subset [1,r]$
let $\wnot (I)$ denote the longest element
of the parabolic subgroup in $W$ generated by all $s_i$ with $i \in I$.

\begin{theorem}
\label{th:string cone split}
Let $\emptyset = I_0 \subset I_1 \subset \cdots \subset I_p = [1,r]$ be any flag of subsets
in $[1,r]$.
Suppose $\ii \in R(\wnot)$ is the concatenation $(\ii^{(1)}, \dots, \ii^{(p)})$,
where $\ii^{(j)} \in R(\wnot (I_{j-1})^{-1} \wnot (I_j))$ for $j = 1, \dots, p$.
Then the string cone $C_\ii$ is the direct product of cones:
\begin{equation}
\label{eq:string cone split}
C_\ii =  C_{\ii^{(1)}} (e, \wnot(I_1)) \times C_{\ii^{(2)}} (\wnot (I_1)), \wnot(I_2))
\times \cdots \times C_{\ii^{(p)}} (\wnot(I_{p-1}), \wnot(I_p)) \ .
\end{equation}
\end{theorem}

Under some additional assumptions, it is possible to describe
the factors in (\ref{eq:string cone split}) much more explicitly.
Following \cite{st}, we call an element $w \in W$ \emph{fully commutative}
if any two reduced words for $w$ can be obtained from each other by a
series of switches $(i_k,i_{k+1}) \to (i_{k+1},i_k)$ such that $a_{i_k,i_{k+1}} = 0$.

\begin{theorem}
\label{th:fully commutative factors}
In the situation of Theorem~\ref{th:string cone split}, suppose that an index
$j \in [1,p]$ is such that $|I_j| = |I_{j-1}| + 1$, and the element
$\wnot (I_{j-1})^{-1} \wnot (I_j)$ is fully commutative of length $l$.
Then the corresponding cone $C_{\ii^{(j)}} (\wnot(I_{j-1}), \wnot(I_j)) \subset \RR^l$
is given by the following inequalities:

\smallskip

\noindent {\rm (1)} $t_l \geq 0$;

\smallskip

\noindent {\rm (2)} $t_{k(1)} \geq t_{k(2)}$ for any pair of indices $k(1) < k(2)$ in $[1,l]$
such that $\alpha_{i_{k(1)}}$ and $\alpha_{i_{k(2)}}$ generate a root
system of type $A_2$, and $i_k \notin \{i_{k(1)}, i_{k(2)}\}$
for $k(1) < k < k(2)$;

\smallskip

\noindent {\rm (3)} $|a_{ij}| t_{k(1)} \geq t_{k(2)} \geq |a_{ij}| t_{k(3)}$ for any indices
$k(1) < k(2) < k(3)$ in $[1,l]$ such that
$i_{k(1)} = i_{k(3)} = j \neq i = i_{k(2)}$,
$\alpha_{i}$ and $\alpha_{j}$ generate a root system of type $B_2$,
and $a_{i_k,j} = 0$ for all $k \neq k(2)$ between $k(1)$ and $k(3)$;

\smallskip

\noindent {\rm (4)} $2 |a_{ij}| t_{k(1)} \geq 2 t_{k(2)} \geq |a_{ij}| t_{k(3)} \geq 2 t_{k(4)}
\geq 2 |a_{ij}| t_{k(5)}$ for any indices
$k(1) < \cdots < k(5)$ in $[1,l]$ such that
$i_{k(1)} = i_{k(3)} = i_{k(5)}= j \neq i = i_{k(2)} = i_{k(4)}$, and
$\alpha_{i}$ and $\alpha_{j}$ generate a root system of type $G_2$.
\end{theorem}

Combining Theorems~\ref{th:string cone split} and \ref{th:fully commutative factors}, we
obtain the following refinement of the main result of \cite{lit}.

\begin{corollary}
\label{cor:fully commutative}
Suppose that $\emptyset = I_0 \subset I_1 \subset \cdots \subset I_r = [1,r]$ is a
flag of subsets in $[1,r]$ such that $|I_j| = j$, and
$\wnot (I_{j-1})^{-1} \wnot (I_j)$ is fully commutative for every $j \in [1,r]$.
Suppose $\ii \in R(\wnot)$ is the concatenation $(\ii^{(1)}, \dots, \ii^{(r)})$,
where $\ii^{(j)} \in R(\wnot (I_{j-1})^{-1} \wnot (I_j))$.
Then the string cone $C_\ii$ is the direct product of cones given by
{\rm (\ref{eq:string cone split})}, with every factor in the product
given by inequalities in Theorem~\ref{th:fully commutative factors}.
\end{corollary}

Theorem~\ref{th:string cones general} also implies a more explicit description
of \emph{all} string cones for the type $A_r$
(another description of these cones was found in \cite{gp}).
We need the following notation: for any $i \in [1,r]$,
let $u^{(i)}$ denote the minimal representative of the coset
$W_{\widehat i} s_i w_0$ in $W$, where $W_{\widehat i}$ is the maximal parabolic
subgroup in $W$ generated by all $s_j$ with $j \neq i$.

\begin{theorem}
\label{th:minuscule}
Let $\ii = (i_1, \dots, i_m) \in R(\wnot)$.
For any $i \in [1,r]$ and any subword $(i_{k(1)}, \dots, i_{k(p)})$ of $\ii$ which
is a reduced word for $u^{(i)}$, all the points
$(t_1, \dots, t_m)$ in the string cone $C_\ii$ satisfy the inequality
\begin{equation}
\label{eq:extremal inequalities}
\sum_{j=0}^p \sum_{k(j) < k < k(j+1)} (s_{i_{k(1)}} \cdots s_{i_{k(j)}} \alpha_{i_k})
(\omega_i^\vee) \cdot t_k \geq 0
\end{equation}
(with the convention that $k(0) = 0$ and $k(p+1) = m+1$).
Furthermore, if $\gg = sl_{r+1}$ then
$C_\ii$ is the set of all $t \in \RR^m$ satisfying the inequalities
{\rm (\ref{eq:extremal inequalities})}.
\end{theorem}

\begin{example}
\label{ex:string cones for A3}
{\rm Let us illustrate Theorem~\ref{th:minuscule} by an example when
$\gg = sl_4$.
We shall use the standard numeration of simple roots and corresponding simple reflections
for type $A_3$, so that $s_1$ and $s_3$ commute with each other.
In our case $m = \l(\wnot) = 6$, and the elements $u^{(i)}$ are given by:
$$u^{(1)} = s_1 s_2, \quad u^{(2)} = s_2 s_1 s_3 = s_2 s_3 s_1, \quad u^{(3)} = s_3 s_2 \ .$$

Let us consider a reduced word $\ii = (2,1,3,2,1,3)$ of $\wnot$.
Here is the list of subwords of $\ii$ which are reduced words for the elements $u^{(i)}$
(their entries are underlined), and the corresponding inequalities
of the form (\ref{eq:extremal inequalities}):
$$\begin{array}{rcl}
u^{(1)} & : & (2, \underline 1,3,\underline 2,1,3) \to t_6 \geq 0; \\[.1in]
u^{(2)} & : & (\underline 2, \underline 1,\underline 3,2,1,3) \to t_4 - t_5 - t_6 \geq 0,
\,\, (\underline 2, \underline 1, 3,2,1,\underline 3) \to t_3 - t_5 \geq 0, \\[.05in]
& & (\underline 2, 1,\underline 3,2,\underline 1, 3) \to t_2 - t_6 \geq 0,
\,\, (\underline 2, 1, 3,2,\underline 1,\underline 3) \to t_2 + t_3 - t_4 \geq 0, \\[.05in]
& & (2, 1, 3,\underline 2,\underline 1,\underline 3) \to t_1 \geq 0; \\[.1in]
u^{(3)} & : & (2, 1,\underline 3,\underline 2,1,3) \to t_5 \geq 0.
\end{array}$$
Therefore, $C_{\ii}$ is a cone in $\RR^6$ given by:
$$t_1 \geq 0, \,\, t_2 \geq t_6 \geq 0, \,\, t_3 \geq t_5 \geq 0, \,\,
t_2+t_3  \geq t_4 \geq t_5 + t_6 \ .$$}
\end{example}



\section{Reduced double Bruhat cells and totally positive varieties}
\label{sec:geometric}

\subsection{Background on semisimple groups}
\label{sec:geometric prelims}
Let $G$ be a simply connected complex semisimple Lie group with the Lie algebra $\gg$.
For $i \in [1,r]$, we denote by $x_i (t)$ and $y_i (t)$ the one-parameter subgroups
in $G$ given by
$$x_i (t) = \exp \ (t e_i), \,\, y_i (t) = \exp \ (t f_i)$$
(note that in \cite{fz} the notation $x_{\overline i} (t)$ was used instead of
$y_i (t)$).
Let $N$ (resp. $N_-$) be the maximal unipotent subgroup of $G$ generated by all
$x_i (t)$ (resp. $y_i (t)$).
Let $H$ be the maximal torus in $G$ with the Lie algebra $\hh$.
Let $B = HN$ and $B_- = HN_-$ be two
opposite Borel subgroups, so that $H=B_-\cap B$.
For every $i\in [1,r]$, let $\varphi_i: SL_2 \to G$ denote
the canonical embedding corresponding to the simple
root~$\alpha_i\,$; thus we have
$$x_i (t) = \varphi_i \mat{1}{t}{0}{1}, \,\, y_i (t) = \varphi_i \mat{1}{0}{t}{1} \ .$$
We also set
$$t^{\alpha_i^\vee} = \varphi_i \mat{t}{0}{0}{t^{-1}} \in H$$
for any $i$ and any $t \neq 0$.

The Weyl group $W$ of $\gg$ is naturally identified with ${\rm Norm}_G (H)/H$;
this identification sends each simple reflection
$s_i$ to the coset $\overline {s_i} H$, where the representative
$\overline {s_i} \in {\rm Norm}_G (H)$ is defined by
\[
\overline {s_i} = \varphi_i \mat{0}{-1}{1}{0} \, .
\]
The elements $\overline {s_i}$ satisfy the braid relations in~$W$;
thus the representative $\overline w$ can be unambiguously defined for any
$w \in W$ by requiring that
$\overline {uv} = \overline {u} \cdot \overline {v}$
whenever $\l (uv) = \l (u) + \l (v)$.

The weight lattice $P$ is identified with the group of
multiplicative characters of~$H$,
here written in the exponential notation:
a weight $\gamma\in P$ acts by $a \mapsto a^\gamma$.
Under this identification, the fundamental weights
$\omega_1, \ldots, \omega_r$ act in $H$ by
$(t^{\alpha_j^\vee})^{\omega_i} = t^{\delta_{ij}}$.
The action of~$W$ on $P$ can be now written as
$a^{w (\gamma)} = (w^{-1} a w)^\gamma$
for $w \in W$, $a \in H$, $\gamma \in P$.

We denote by $G_0=N_-HN$ the open subset of elements $x\in G$ that
have Gaussian decomposition; this (unique) decomposition will be written as
$x = [x]_- [x]_0 [x]_+ \,$.

Following G.~Lusztig, we define the variety $G_{\geq 0}$ of totally nonnegative elements in $G$
as the multiplicative monoid with unit generated by the elements
$t^{\alpha_i^\vee}$, $x_i (t)$, and
$y_i (t)$ for all $i$ and all $t > 0$.

\subsection{Preliminaries on generalized minors}
\label{sec:generalized minors}
We now recall some basic properties of \emph{generalized minors} introduced in \cite{fz}.
For $u,v \in W$ and $i \in [1,r]$, the \emph{generalized minor}
$\Delta_{u \omega_i, v \omega_i}$
is the regular function on $G$ whose restriction to the open set
${\overline {u}} G_0 {\overline {v}}^{-1}$ is given by
\begin{equation}
\label{eq:Delta-general}
\Delta_{u \omega_i, v \omega_i} (x) =
(\left[{\overline {u}}^{\ -1}
   x \overline v\right]_0)^{\omega_i} \ .
\end{equation}
As shown in \cite{fz}, $\Delta_{u \omega_i, v \omega_i}$ depends on
the weights $u \omega_i$ and $v \omega_i$ alone, not on the particular
choice of $u$ and~$v$.
In the special case $G=SL_n\,$, the generalized minors are nothing but
the ordinary minors of a matrix.

Although we do not need it in this paper, we would like to mention the following
characterization of the totally nonnegative variety obtained in \cite{fzoscil}:
an element $x \in G$ is totally nonnegative if and only if all generalized minors
take nonnegative real values at $x$.

Generalized minors have the following  properties:
(see \cite[(2.14), (2.25)]{fz}):
\begin{equation}
\label{eq:minor properties}
\Delta_{\gamma, \delta} (a_1 x a_2) = a_1^\gamma a_2^\delta
\Delta_{\delta, \gamma} (x) \,\, (a_1, a_2  \in H; \ x \in G),
\quad  \Delta_{\gamma,\delta} (x) =  \Delta_{\delta, \gamma} (x^T) \ ,
\end{equation}
where $x \mapsto x^T$ is the ``transpose" involutive antiautomorphism of
$G$ given by
\begin{equation}
\label{eq:T}
a^T = a \quad (a \in H) \ , \quad x_i (t)^T = y_i (t) \ ,
\quad y_i (t)^T = x_i (t) \ .
\end{equation}

Later we shall need the involutive antiautomorphism $\tau_{\wnot}$ of $G$
introduced in \cite[(2.56)]{fz}; it is defined by
\begin{equation}
\label{eq:tau-wnot}
\tau_{\wnot} (x) = \overline \wnot (x^{-1})^{\iota T} {\overline \wnot}^{\ -1} \ ,
\end{equation}
where $x \mapsto x^\iota$ is the involutive antiautomorphism of $G$
given by
\begin{equation}
\label{eq:iota}
a^\iota = a^{-1} \quad (a \in H) \ , \quad x_i (t)^\iota = x_i (t) \ ,
\quad y_{i} (t)^\iota = y_{i} (t) \ .
\end{equation}
By \cite[(2.25) and Lemma~2.25]{fz}, we have
\begin{equation}
\label{eq:minors-tau}
\Delta_{\gamma, \delta} (x) = \Delta_{-\delta, -\gamma} (x^\iota) =
\Delta_{\wnot \delta, \wnot \gamma} (\tau_{\wnot}(x))
\end{equation}
for any generalized minor $\Delta_{\gamma, \delta}$, and any $x \in G$.

Now we present some less obvious identities for generalized minors.
The following identity was obtained in \cite[Theorem~1.17]{fz}.

\begin{proposition}
\label{pr:minors-Dodgson}
Suppose $u,v \in W$ and $i \in [1,r]$
are such that $\l (us_i) = \l (u) + 1$ and $\l (vs_i) = \l (v) + 1$.
Then
\begin{eqnarray}
\begin{array}{l}
\label{eq:minors-Dodgson}
\Delta_{u \omega_i, v \omega_i} \Delta_{us_i \omega_i, v s_i \omega_i}
= \Delta_{us_i \omega_i, v \omega_i} \Delta_{u \omega_i, v s_i \omega_i}
+ \prod_{j \neq i} \Delta_{u \omega_j, v \omega_j}^{- a_{ji}} \ .
\end{array}
\end{eqnarray}
\end{proposition}

The next proposition presents some generalized Pl\"ucker relations;
they follow from~\cite[Corollary~6.6]{bz2} (see also \cite[Theorem~1.16]{fz}).



\begin{proposition}
\label{pr:minors-Plucker}
Let $u, v \in W$ and $i, j \in [1,r]$.

\noindent \emph{1.}
If $a_{ij}\!=\!a_{ji}\!=\!-1$
and $\l (vs_i s_j s_i) = \l (v) + 3$, then
\[
\Delta_{u \omega_i, vs_i \omega_i} \Delta_{u \omega_j, vs_j
  \omega_j} =
\Delta_{u \omega_i, v\omega_i} \Delta_{u \omega_j, vs_i s_j \omega_j} +
\Delta_{u \omega_i, vs_j s_i \omega_i} \Delta_{u \omega_j, v \omega_j} \ .
\]

\noindent \emph{2.}
If $a_{ij}\! =\! -2$, $a_{ji}\! =\! -1$, and
$\l (vs_i s_j s_i s_j) = \l (v) + 4$, then
\[\begin{array}{l}
\Delta_{u \omega_i, vs_i \omega_i} \Delta_{u \omega_i, vs_j s_i \omega_i}
\Delta_{u \omega_j, vs_j \omega_j} \\[.1in]
=\Delta^2_{u \omega_i, v s_j s_i \omega_i}
 \Delta_{u \omega_j, v \omega_j}
\\[.1in]
\qquad\qquad
+\Delta_{u \omega_i, v \omega_i}
(\Delta_{u \omega_i, v \omega_i}
 \Delta_{u \omega_j, vs _js_i s_j \omega_j}
+\Delta_{u \omega_i, v s_i s_j s_i \omega_i}
 \Delta_{u \omega_j, vs_j \omega_j})
\end{array}\]
and
\[\begin{array}{l}
\Delta_{u \omega_j, vs_i s_j\omega_j}
\Delta^2_{u \omega_i, vs_j s_i \omega_i}
\Delta_{u \omega_j, vs_j \omega_j}
\\[.1in]
=
\Delta_{u \omega_j, vs _j s_i s_j \omega_j}
\Delta^2_{u \omega_i, v s_j s_i \omega_i}
\Delta_{u \omega_j, v \omega_j}
\\[.1in]
\qquad\qquad\qquad\qquad
+(\Delta_{u \omega_i, v \omega_i}
\Delta_{u \omega_j, vs _j s_i s_j \omega_j}\! +\!
\Delta_{u \omega_i, v s_i s_j s_i \omega_i}
\Delta_{u \omega_j, vs_j \omega_j})^2 \ .
\end{array}\]

\noindent \emph{3.}
If $a_{ij}\! =\! -3$, $a_{ji}\! =\! -1$, and
$l(vs_i s_j s_i s_j s_i s_j)= l(v) + 6$, then
\[\begin{array}{l}
 \Delta_{u \omega_i, vs_i \omega_i} \Delta^2_{u \omega_i,vs_j s_i \omega_i}
\Delta_{u \omega_i,vs_j s_i s_j s_i\omega_i}
\Delta_{u \omega_j,vs_j \omega_j} \Delta_{u \omega_j,vs_j s_i s_j \omega_j} \\[.1in]
\qquad = \Delta^3_{u \omega_i,vs_j s_i \omega_i} \Delta_{u \omega_i,vs_j s_i s_j s_i\omega_i}
\Delta_{u \omega_j,v\omega_j} \Delta_{u \omega_j,vs_j s_i s_j \omega_j} \\[.1in]
 + \bigl(\Delta_{u \omega_i,v\omega_i} \Delta_{u \omega_j,vs_j s_i s_j \omega_j} +
\Delta_{u \omega_i,vs_j s_i s_j s_i\omega_i} \Delta_{u \omega_j,vs_j \omega_j}\bigr)^2
\Delta_{u \omega_i,v\omega_i} \Delta_{u \omega_i,vs_j s_i s_j s_i\omega_i} \\[.1in]
 + \bigl(\Delta_{u \omega_i,vs_j s_i \omega_i} \Delta_{u \omega_j,vs_j s_i s_j s_i s_j \omega_j} +
\Delta_{u \omega_i,vs_i s_j s_i s_j s_i\omega_i} \Delta_{u \omega_j,vs_j s_i s_j \omega_j}\bigr)
\\[.1in]
\qquad \cdot \Delta_{u \omega_i,v\omega_i} \Delta^2_{u \omega_i,vs_j s_i\omega_i}
\Delta_{u \omega_j,vs_j \omega_j} \ ;
\end{array}\]
\[\begin{array}{l}
 \Delta_{u \omega_i,vs_i s_j s_i \omega_i} \Delta^3_{u \omega_i,vs_j s_i \omega_i}
\Delta^2_{u \omega_i,vs_j s_i s_j s_i\omega_i}
\Delta_{u \omega_j,vs_j \omega_j} \Delta^2_{u \omega_j,vs_j s_i s_j \omega_j} \\[.1in]
 \qquad = \Delta^3_{u \omega_i,vs_j s_i \omega_i} \Delta^3_{u \omega_i,vs_j s_i s_j s_i\omega_i}
\Delta_{u \omega_j,v\omega_j} \Delta^2_{u \omega_j,vs_j s_i s_j \omega_j} \\[.1in]
 + \bigl(\Delta_{u \omega_i,v\omega_i} \Delta_{u \omega_j,vs_j s_i s_j \omega_j} +
\Delta_{u \omega_i,vs_j s_i s_j s_i\omega_i} \Delta_{u \omega_j,vs_j \omega_j}\bigr)^3
\Delta^3_{u \omega_i,vs_j s_i s_j s_i\omega_i} \\[.1in]
 + \bigl(\Delta_{u \omega_i,vs_j s_i \omega_i} \Delta_{u \omega_j,vs_j s_i s_j s_i s_j \omega_j} +
\Delta_{u \omega_i,vs_i s_j s_i s_j s_i\omega_i} \Delta_{u \omega_j,vs_j s_i s_j \omega_j}\bigr)^2
\\[.1in]
\qquad \cdot \Delta^4_{u \omega_i,vs_j s_i\omega_i} \Delta_{u \omega_j,vs_j \omega_j} \\[.1in]
 + \bigl(3 \Delta_{u \omega_i,v\omega_i} \Delta_{u \omega_i,vs_j s_i \omega_i}
\Delta_{u \omega_j,vs_j s_i s_j \omega_j} \Delta_{u \omega_j,vs_j s_i s_j s_i s_j \omega_j} \\[.1in]
\qquad \qquad +
2 \Delta_{u \omega_i,v\omega_i} \Delta_{u \omega_i,vs_i s_j s_i s_j s_i\omega_i}
\Delta^2_{u \omega_j,vs_j s_i s_j \omega_j} \\[.1in]
\qquad \qquad + 2 \Delta_{u \omega_i,vs_j s_i s_j s_i \omega_i}
\Delta_{u \omega_i,vs_i s_j s_i s_j s_i \omega_i}
\Delta_{u \omega_j,vs_j \omega_j} \Delta_{u \omega_j,vs_j s_i s_j \omega_j} \\[.1in]
 + 2 \Delta_{u \omega_i,vs_j s_i \omega_i} \Delta_{u \omega_i,vs_j s_i s_j s_i \omega_i}
\Delta_{u \omega_j,vs_j \omega_j} \Delta_{u \omega_j,vs_j s_i s_j s_i s_j \omega_j} \bigr)
\\[.1in]
\qquad \cdot \Delta^2_{u \omega_i,vs_j s_i \omega_i} \Delta^2_{u \omega_i,vs_j s_i s_j s_i \omega_i}
\Delta_{u \omega_j,vs_j \omega_j} \ ;
\end{array}\]
\[\begin{array}{l}
\Delta_{u \omega_j,vs_i s_j s_i s_j \omega_j} \Delta^3_{u \omega_i,vs_j s_i \omega_i}
\Delta^3_{u \omega_i,vs_j s_i s_j s_i\omega_i}
\Delta_{u \omega_j,vs_j \omega_j} \Delta^2_{u \omega_j,vs_j s_i s_j \omega_j}  \\[.1in]
= \Delta^3_{u \omega_i,vs_j s_i \omega_i} \Delta^3_{u \omega_i,vs_j s_i s_j s_i\omega_i}
\Delta_{u \omega_j,v\omega_j} \Delta^2_{u \omega_j,vs_j s_i s_j \omega_j}
\Delta_{u \omega_j,vs_j s_i s_j s_i s_j \omega_j} \\[.1in]
 + \bigl(\Delta_{u \omega_i,v\omega_i} \Delta_{u \omega_j,vs_j s_i s_j \omega_j} +
\Delta_{u \omega_i,vs_j s_i s_j s_i\omega_i} \Delta_{u \omega_j,vs_j \omega_j}\bigr)^3 \\[.1in]
\qquad \cdot \Delta^3_{u \omega_i,vs_j s_i s_j s_i\omega_i}
\Delta_{u \omega_j,vs_j s_i s_j s_i s_j \omega_j}
\\[.1in]
+ \bigl(\Delta_{u \omega_i,vs_j s_i \omega_i} \Delta_{u \omega_j,vs_j s_i s_j s_i s_j \omega_j} +
\Delta_{u \omega_i,vs_i s_j s_i s_j s_i\omega_i} \Delta_{u \omega_j,vs_j s_i s_j \omega_j}\bigr)^3
\\[.1in]
\qquad \cdot \Delta^3_{u \omega_i,vs_j s_i\omega_i} \Delta_{u \omega_j,vs_j \omega_j} \\[.1in]
+ \bigl(3 \Delta_{u \omega_i,v\omega_i} \Delta_{u \omega_i,vs_j s_i \omega_i}
\Delta_{u \omega_j,vs_j s_i s_j \omega_j} \Delta_{u \omega_j,vs_j s_i s_j s_i s_j \omega_j} \\[.1in]
\qquad \qquad + 3 \Delta_{u \omega_i,v\omega_i} \Delta_{u \omega_i,vs_i s_j s_i s_j s_i\omega_i}
\Delta^2_{u \omega_j,vs_j s_i s_j \omega_j} \\[.1in]
\qquad \qquad
+ 3 \Delta_{u \omega_i,vs_j s_i s_j s_i \omega_i}
\Delta_{u \omega_i,vs_i s_j s_i s_j s_i \omega_i}
\Delta_{u \omega_j,vs_j \omega_j} \Delta_{u \omega_j,vs_j s_i s_j \omega_j} \\[.1in]
+ 2 \Delta_{u \omega_i,vs_j s_i \omega_i} \Delta_{u \omega_i,vs_j s_i s_j s_i \omega_i}
\Delta_{u \omega_j,vs_j \omega_j} \Delta_{u \omega_j,vs_j s_i s_j s_i s_j \omega_j} \bigr) \\[.1in]
\qquad \cdot \Delta^2_{u \omega_i,vs_j s_i \omega_i}
\Delta^2_{u \omega_i,vs_j s_i s_j s_i \omega_i}
\Delta_{u \omega_j,vs_j \omega_j} \Delta_{u \omega_j,vs_j s_i s_j s_i s_j \omega_j} \ ;
\end{array}\]
\[\begin{array}{l}
\Delta_{u \omega_j,vs_i s_j \omega_j} \Delta^6_{u \omega_i,vs_j s_i \omega_i}
\Delta^3_{u \omega_i,vs_j s_i s_j s_i\omega_i}
\Delta^2_{u \omega_j,vs_j \omega_j} \Delta^3_{u \omega_j,vs_j s_i s_j \omega_j} \\[.1in]
\qquad = \Delta^3_{u \omega_i,vs_j s_i \omega_i} \Delta^2_{u \omega_i,vs_j s_i s_j s_i \omega_i}
\Delta_{u \omega_j,v\omega_j} \Delta^2_{u \omega_j,vs_j s_i s_j \omega_j} \\[.1in]
\cdot \bigl\{\Delta^3_{u \omega_i,vs_j s_i \omega_i}
\Delta_{u \omega_i,vs_j s_i s_j s_i\omega_i}
\Delta_{u \omega_j,v\omega_j} \Delta^2_{u \omega_j,vs_j s_i s_j \omega_j} \\[.1in]
\qquad  + 2 \bigl(\Delta_{u \omega_i,v\omega_i} \Delta_{u \omega_j,vs_j s_i s_j \omega_j} +
\Delta_{u \omega_i,vs_j s_i s_j s_i\omega_i} \Delta_{u \omega_j,vs_j \omega_j}\bigr)^3
\Delta_{u \omega_i,vs_j s_i s_j s_i\omega_i} \\[.1in]
+ \bigl(3 \Delta_{u \omega_i,v\omega_i} \Delta_{u \omega_i,vs_j s_i \omega_i}
\Delta_{u \omega_j,vs_j s_i s_j \omega_j} \Delta_{u \omega_j,vs_j s_i s_j s_i s_j \omega_j} \\[.1in]
\qquad \qquad + 3 \Delta_{u \omega_i,v\omega_i} \Delta_{u \omega_i,vs_i s_j s_i s_j s_i\omega_i}
\Delta^2_{u \omega_j,vs_j s_i s_j \omega_j} \\[.1in]
\qquad \qquad + 3 \Delta_{u \omega_i,vs_j s_i s_j s_i \omega_i}
\Delta_{u \omega_i,vs_i s_j s_i s_j s_i \omega_i}
\Delta_{u \omega_j,vs_j \omega_j} \Delta_{u \omega_j,vs_j s_i s_j \omega_j} \\[.1in]
+ 2 \Delta_{u \omega_i,vs_j s_i \omega_i} \Delta_{u \omega_i,vs_j s_i s_j s_i \omega_i}
\Delta_{u \omega_j,vs_j \omega_j} \Delta_{u \omega_j,vs_j s_i s_j s_i s_j \omega_j} \bigr)
\\[.1in]
\qquad \cdot \Delta^2_{u \omega_i,vs_j s_i \omega_i} \Delta_{u \omega_j,vs_j \omega_j} \bigr \} \\[.1in]
+ \bigl \{ \bigl(\Delta_{u \omega_i,v\omega_i} \Delta_{u \omega_j,vs_j s_i s_j \omega_j} +
\Delta_{u \omega_i,vs_j s_i s_j s_i\omega_i} \Delta_{u \omega_j,vs_j \omega_j}\bigr)^2
\Delta_{u \omega_i,vs_j s_i s_j s_i\omega_i} \\[.1in]
+ \bigl(\Delta_{u \omega_i,vs_j s_i \omega_i} \Delta_{u \omega_j,vs_j s_i s_j s_i s_j \omega_j} +
\Delta_{u \omega_i,vs_i s_j s_i s_j s_i\omega_i} \Delta_{u \omega_j,vs_j s_i s_j \omega_j}\bigr)
\\[.1in]
\qquad \cdot \Delta^2_{u \omega_i,vs_j s_i \omega_i} \Delta_{u \omega_j,vs_j \omega_j}
\bigr \}^3 \ .
\end{array}\]
\end{proposition}

\subsection{Reduced double Bruhat cells}
\label{sec:reduced Bruhat}

The group $G$ has two \emph{Bruhat decompositions},
with respect to opposite Borel subgroups $B$ and $B_-\,$:
$$G = \bigcup_{u \in W} B u B = \bigcup_{v \in W} B_- v B_-  \ . $$
The \emph{double Bruhat cells}~$G^{u,v}$ are defined by
$G^{u,v} = B u B  \cap B_- v B_- \,$.
These varieties were introduced and studied in \cite{fz}.

In this paper we shall concentrate on the following subset $L^{u,v} \subset G^{u,v}$
which we call a \emph{reduced double Bruhat cell}:
\begin{equation}
\label{eq:reduced cell}
L^{u,v} = N \overline u N  \cap B_- v B_- \ .
\end{equation}
(In particular, if $u$ is the identity element $e \in W$ then $L^{e,v} = N \cap B_- v B_-$
is the variety $N^v$ studied in \cite{bz2}.)
We also set
$$L^{u,v}_{> 0} := L^{u,v} \cap G_{\geq 0} \ ,$$
and call $L^{u,v}_{> 0}$ the \emph{totally positive part} of $L^{u,v}$.

The equations defining $L^{u,v}$ inside $G^{u,v}$ look as follows.

\begin{proposition}
\label{pr:Luv equations}
An element $x \in G^{u,v}$ belongs to $L^{u,v}$ if and only if
$[{\overline u}^{\ -1} x]_0 = 1$, or equivalently if $\Delta_{u \omega_i, \omega_i} (x) = 1$
for any $i \in  [1,r]$.
\end{proposition}

The maximal torus $H$ acts freely on $G^{u,v}$ by left (or right) translations,
and $L^{u,v}$ is a section of this action.
Thus $L^{u,v}$ is naturally identified with $G^{u,v}/H$, and all properties of $G^{u,v}$
can be translated in a straightforward way into the corresponding properties of $L^{u,v}$.
In particular, Theorem~1.1 in \cite{fz} implies the following.

\begin{proposition}
\label{pr:Luv-affine}
The variety $L^{u,v}$ is biregularly isomorphic to a Zariski
open subset of an affine space of dimension $\l(u)+\l(v)$.
\end{proposition}

We now introduce a family of systems of local coordinates in $L^{u,v}$.
For any nonzero $t \in \CC$ and any $i\in [1,r]$, denote
\begin{equation}
\label{eq:xnegative}
x_{- i} (t) = y_i (t) t^{- \alpha_i^\vee}  = \varphi_i \mat{t^{-1}}{0}{1}{t} \, .
\end{equation}

A \emph{double reduced word} for a pair $u,v \in W$
is a reduced word for an element $(u,v)$
of the Coxeter group $W \times W$.
To avoid confusion, we will use the indices
$- 1, \ldots, - r$ for the simple reflections
in the first copy of $W$, and $1, \ldots, r$ for the second copy.
A double reduced word for $(u,v)$ is nothing but a shuffle of a reduced
word $\ii$ for $u$ written in the alphabet
$[- 1, - r]$ ( we will denote such a word by $- \ii$) and a reduced word $\ii'$
for $v$ written in the alphabet $[1,r]$.
We denote the set of double reduced words for $(u,v)$ by $R(u,v)$.

For any sequence $\ii= (i_1, \ldots, i_m)$ of indices
from the alphabet $[1,r] \cup [- 1, - r]$,
let us define the \emph{product map} $x_\ii: \CC_{\neq 0}^m \to G$ by
\begin{equation}
\label{eq:productmap}
x_\ii (t_1, \ldots, t_m) =  x_{i_1} (t_1) \cdots x_{i_m} (t_m) \, .
\end{equation}
The following proposition is a modification of \cite[Theorems 1.2, 1.3]{fz}.

\begin{proposition}
\label{pr:Luv coordinates}
For any $u,v\in W$ and $\ii= (i_1, \ldots, i_m)\in R(u,v)$,
the map $x_\ii$ is a biregular isomorphism between
$\CC_{\neq 0}^m$ and a Zariski open subset of~$L^{u,v}$, and
it restricts to a bijection between the set
$\RR_{> 0}^m$ of $m$-tuples of positive real numbers and the totally positive part
$L^{u,v}_{> 0}$ of $L^{u,v}$.
\end{proposition}

\subsection{Factorization problem for reduced double Bruhat cells}
\label{sec:factorization}
In this section, we address the following
\emph{factorization problem} for $L^{u,v}$: for any double reduced word
$\ii \in R(u,v)$, find explicit formulas for the inverse birational isomorphism
$x_\ii^{-1}$ between $L^{u,v}$ and $\CC_{\neq 0}^m$,
thus expressing the \emph{factorization parameters} $t_k$
in terms of the product~$x = x_\ii (t_1, \ldots, t_m) \in L^{u,v}$.

Recall from \cite[Section 1.5]{fz} that there is a biregular ``shift" isomorphism
$\zeta^{u,v}: G^{u,v} \to G^{u^{-1}, v^{-1}}$.
This isomorphism does not however send $L^{u,v}$ to $L^{u^{-1}, v^{-1}}$, so
we will use the following modification.

\begin{definition}
\label{def:twist}{\rm
For any $u,v\in W$, the twist map
$\psi^{u,v}$ is defined by (see (\ref{eq:iota}))}
\begin{equation}
\label{eq:psi-u,v-x}
\psi^{u,v}(x) =
[(\overline{v} x^{\iota})^{-1}]_+ \, \overline{v} \,
([\overline{u}^{\ -1} x]_+)^{\iota} \ .
\end{equation}
\end{definition}

\begin{theorem}
\label{th:psi-regularity}
The twist map $\psi^{u,v}$ is a biregular isomorphism
between $L^{u,v}$ and  $L^{v, u}$, and it restricts to a bijection between
$L^{u,v}_{> 0}$ and  $L^{v, u}_{>0}$.
The inverse isomorphism is $\psi^{v,u}$.
\end{theorem}

Now let us fix a pair $(u,v) \in W \times W$ and a double reduced word
$\ii = (i_1, \ldots, i_m)\in R(u,v)$.
Recall that $\ii$ is a shuffle of a reduced word for $u$ written in the alphabet
$[- 1, - r]$ and a reduced word for $v$ written in the alphabet $[1,r]$.
In particular, the length $m$ of $\ii$ is equal to $\l (u) + \l (v)$.
For $k\in [1,m]$, we denote
\[
k^- = \max\{l:l<k, |i_l| = |i_k|\} \ , \,\, k^+ = \min\{l:l>k, |i_l| = |i_k|\} \ ,
\]
so that $k^-$ (resp. $k^+$) is the previous (resp. next) occurrence of an index
$\pm i_k$ in $\ii$; if $k$ is the first (resp. last) occurrence of $\pm i_k$ in $\ii$ then we set
$k^- = 0$ (resp. $k^+ = m+1$).

For $k\in [1,m]$, denote
\begin{equation}
\label{eq:v_leq k}
u_{\geq k} = \doublesubscript{\prod}{l = m, \ldots, k}
{i_l < 0}
s_{-i_l} \ , \quad
v_{<k} = \doublesubscript{\prod}{l = 1, \ldots, k-1}
{i_l > 0} s_{i_l} \ .
\end{equation}
This notation means that in the first (resp.\ second) product in
(\ref{eq:v_leq k}), the index $l$ is decreasing (resp. increasing);
for example, if
$\ii= (-2, 1, -3, 3, 2, -1, -2, 1, -1)$,
then, say, $u_{\geq 7} = s_1 s_2 \,$ and $v_{<7} = s_1 s_3 s_2\,$.
Let us define a regular function $M_k = M_{k,\ii}$ on $L^{u,v}$ by
\begin{equation}
\label{eq:M-factors}
M_k (x) = M_{k,\ii} (x) =
\Delta_{v_{<k} \omega_{|i_k|}, u_{\geq k} \omega_{|i_k|}} (\psi^{u,v} (x)) \ ;
\end{equation}
by convention, we set $M_{m+1} (x) = 1$ (by Theorem~\ref{th:psi-regularity},
$\psi^{u,v} (x) \in L^{v,u}$ for $x \in L^{u,v}$, hence
$\Delta_{v \omega_i,  \omega_i} (\psi^{u,v} (x)) = 1$ for all $i$
(see Proposition~\ref{pr:Luv equations})).

Now we are ready to state our solution to the factorization problem.

\begin{theorem}
\label{th:t-through-x}
Let $\ii = (i_1, \ldots, i_m)$ be
a double reduced word for $(u,v)$, and suppose an element
$x \in L^{u,v}$ can be factored as
$x= x_{i_1} (t_1) \cdots x_{i_m} (t_m)$,
with all $t_k$ nonzero complex numbers.
Then the factorization parameters $t_k$
are determined by the following formulas:
if $i_k < 0$ then
\begin{equation}
\label{eq:negative t-through-x}
t_k =  M_k (x) / M_{k^+} (x) \ ;
\end{equation}
if $i_k > 0$ then
\begin{equation}
\label{eq:positive t-through-x}
t_k =  \frac{1}{M_k (x) M_{k^+} (x)} \prod_{l: \ l^- < k < l} M_l (x)^{- a_{|i_l|, i_k}} \ .
\end{equation}
\end{theorem}

The following two special cases will be of
particular importance for us: $(u,v) = (e, \wnot)$,
and $(u,v) = (\wnot, e)$.
In these cases, Definition~\ref{def:twist} and Theorem~\ref{th:psi-regularity}
can be simplified as follows.

\begin{corollary}
\label{cor:twist special}
The twist map $\psi^{\wnot,e}$ is a biregular isomorphism
between $L^{\wnot,e}$ and  $L^{e, \wnot}$ (and also a bijection
between $L^{\wnot,e}_{> 0}$ and  $L^{e, \wnot}_{>0}$) given by
\begin{equation}
\label{eq:psi wnot-e}
\psi^{\wnot,e}(x) = ([\overline{\wnot}^{\ -1} x]_+)^{\iota} \ .
\end{equation}
The inverse biregular isomorphism $\psi^{e,\wnot}$ is given by
\begin{equation}
\label{eq:psi e-wnot}
\psi^{e,\wnot}(x) = ([x \overline{\wnot}]_- \ [x \overline{\wnot}]_0)^{\iota} \ .
\end{equation}
\end{corollary}

The formulas (\ref{eq:negative t-through-x}) and (\ref{eq:positive t-through-x})
now take the following form.

\begin{corollary}
\label{cor:factors special}
Let $\ii = (i_1, \dots, i_m)$ be a reduced word for $\wnot$, and $t_1, \dots, t_m$
be non-zero complex numbers.
\begin{itemize}
\item[(i)] If $x = x_{- i_1} (t_1) \cdots x_{- i_m} (t_m)$ then the factorization parameters
$t_k$ are given by
$$t_k =  \frac{\Delta_{\omega_{i_k}, s_{i_m} \cdots s_{i_k} \omega_{i_k}} (\psi^{\wnot,e}(x))}
{\Delta_{\omega_{i_k}, s_{i_m} \cdots s_{i_{k+1}} \omega_{i_k}} (\psi^{\wnot,e}(x))} \ , $$
where $\psi^{\wnot,e}(x)$ is given by {\rm (\ref{eq:psi wnot-e})}.

\item[(ii)] If $x = x_{i_1} (t_1) \cdots x_{i_m} (t_m)$ then the factorization parameters
$t_k$ are given by
$$t_k = \frac{1}
{\Delta_{s_{i_1} \cdots s_{i_{k-1}} \omega_{i_k}, \omega_{i_k}} (\psi^{e,\wnot}(x))
\Delta_{s_{i_1} \cdots s_{i_{k}} \omega_{i_k}, \omega_{i_k}} (\psi^{e,\wnot}(x))} $$
$$\times \prod_{j \neq i_k}
\Delta_{s_{i_1} \cdots s_{i_{k-1}} \omega_{j}, \omega_{j}} (\psi^{e,\wnot}(x))^{- a_{j, i_k}}
\ , $$
where $\psi^{e,\wnot}(x)$ is given by {\rm (\ref{eq:psi e-wnot})}.
\end{itemize}
\end{corollary}

In general, Theorem~\ref{th:t-through-x} expresses the factorization parameters of an element
$x \in L^{u,v}$ as monomials in generalized minors of the twisted element $\psi^{u,v} (x)$.
However, there are two important special cases when taking the twist is unnecessary.

\begin{proposition}
\label{pr:t1 special}
Let $\ii = (i_1, \dots, i_m)$ be a reduced word for $\wnot$.
\begin{itemize}
\item[(i)] If $x = x_{- i_1} (t_1) \cdots x_{- i_m} (t_m)$ then the first and the last
factorization parameters of $x$ are given by
\begin{equation}
\label{eq:negative t1 special}
t_1 =  \frac{1}
{\Delta_{s_{i_1} \wnot \omega_{i_1^*}, \omega_{i_1^*}} (x)} \ , \,\,
t_m = \Delta_{\wnot \omega_{i_m}, s_{i_m} \omega_{i_m}} (x) \ ,
\end{equation}
where $\omega_{i^*} = - \wnot \omega_i$.

\item[(ii)] If $x = x_{i_1} (t_1) \cdots x_{i_m} (t_m)$ then the first and the last
factorization parameters of $x$ are given by
\begin{equation}
\label{eq:positive t1 special}
t_1 = \frac{\Delta_{\omega_{i_1}, \wnot \omega_{i_1}} (x)}
{\Delta_{s_{i_1} \omega_{i_1}, \wnot \omega_{i_1}} (x)} \ , \,\,
t_m = \frac{\Delta_{\omega_{i_m^*}, \wnot \omega_{i_m^*}} (x)}
{\Delta_{\omega_{i_m^*}, s_{i_m} \wnot \omega_{i_m^*}} (x)} \ .
\end{equation}
\end{itemize}
\end{proposition}

\section{Geometric lifting and tropicalization}
\label{sec:lifting}

\subsection{Transition maps}
\label{sec:transition}
Let $\ii$ and $\ii'$ be two reduced words of $\wnot$.
In view of Proposition~\ref{pr:characterization of BB},
there is a bijective \emph{transition map}
$$R_\ii^{\ii'} =
(b_{\ii'})^{-1} \circ b_\ii:\ZZ_{\ge 0}^m \to \ZZ_{\ge 0}^m $$
between the corresponding Lusztig parametrizations of the canonical basis $\BB$.
Similarly, by Proposition~\ref{pr:strings},
there is a bijective transition map
$$R_{-\ii}^{-\ii'} =
c_{\ii'} \circ (c_\ii)^{-1} :C_\ii (\ZZ) \to C_{\ii'} (\ZZ)$$
between the two string parametrizations of the dual canonical basis
(the reason for the notation $R_{-\ii}^{-\ii'}$ will become clear soon).

It turns out that each component of a tuple $R_\ii^{\ii'} (t)$ or
$R_{-\ii}^{-\ii'}(t)$ can be expressed through the components of $t$
as a ``tropical" subtraction-free rational expression (see Introduction).

\begin{example}
\label{ex:A2 transition}
{\rm Let $\gg = sl_3$ be of type $A_2$, and let $\ii = (1,2,1)$ and $\ii' = (2,1,2)$
be the two reduced words for $\wnot$.
The transition map $R_\ii^{\ii'}$ between two Lusztig parametrizations
was computed in \cite{lu}: the components of $t' = R_\ii^{\ii'}(t)$ are given by
$$t'_1 = t_2 + t_3 - \min \ (t_1, t_3), \,\,  t'_2 = \min \ (t_1, t_3), \,\,
t'_3 = t_1 + t_2 - \min \ (t_1, t_3) \ ,$$
which can also be written as
\begin{equation}
\label{eq:Lustzig A2 transition}
t'_1 = \left[\frac{t_2 t_3}{t_1 + t_3}\right]_{\rm trop}, \,\,
t'_2 = [t_1 + t_3]_{\rm trop}, \,\,
t'_3 = \left[\frac{t_1 t_2}{t_1 + t_3}\right]_{\rm trop} \ .
\end{equation}
The transition map $R_{-\ii}^{-\ii'}$ between two string parametrizations
was computed in \cite{bz1}: the components of $t' = R_{-\ii}^{-\ii'}(t)$ are given by
$$t'_1 = \max \ (t_3, t_2 - t_1), \,\,  t'_2 = t_1 + t_3, \,\,
t'_3 = \min \ (t_1, t_2 - t_3) \ ,$$
which can also be written as
\begin{equation}
\label{eq:string A2 transition}
t'_1 = \left[\frac{t_2 t_3}{t_1 t_3 + t_2}\right]_{\rm trop}, \,\,
t'_2 = [t_1 t_3]_{\rm trop}, \,\,
t'_3 = \left[\frac{t_1 t_3 + t_2}{t_3}\right]_{\rm trop} \ .
\end{equation}
}
\end{example}

We now generalize this example by giving a geometric lifting for each of the transition
maps $R_\ii^{\ii'}$ and $R_{-\ii}^{-\ii'}$.
To do this, we notice that by Proposition~\ref{pr:Luv coordinates},
for any two reduced words $\ii$ and $\ii'$ of a pair $(u,v)$ of elements of $W$,
there is a bijective transition map
$$\tilde R_\ii^{\ii'} =
(x_{\ii'})^{-1} \circ x_\ii:\RR_{> 0}^m \to \RR_{> 0}^m $$
that relates the corresponding parametrizations of the totally positive variety $L^{u,v}_{>0}$.
In particular, any two reduced words $\ii$ and $\ii'$ for $\wnot$ give rise to transition maps
$\tilde R_\ii^{\ii'}$ and $\tilde R_{- \ii}^{-\ii'}$ (for the varieties
$L^{e,\wnot}_{>0}$ and $L^{\wnot, e}_{>0}$, respectively).
We shall use the notation $(\tilde R_\ii^{\ii'})^\vee$ for the transition maps defined in the same
way but for the Langlands dual group ${}^L G$ instead of $G$.

\begin{theorem}
\label{th:geometric lifting}
{\rm (i)} For any $\ii, \ii' \in R(u,v)$, each component of
$\tilde R_\ii^{\ii'} (t)$ is a subtraction-free rational expression in the components
of $t$.

\noindent {\rm (ii)} For any $\ii, \ii' \in R(\wnot)$, each component of
$(\tilde R_\ii^{\ii'})^\vee (t)$ (resp. $(\tilde R_{-\ii}^{-\ii'})^\vee (t)$)
is a geometric lifting of the corresponding component of
$R_\ii^{\ii'} (t)$ (resp. of $R_{-\ii}^{-\ii'} (t)$).
\end{theorem}

As in \cite{bfz} and \cite[Section~4]{bz2}, Theorem~\ref{th:geometric lifting}(i)
allows us to extend the definition of the totally positive varieties $L^{u,v}_{> 0}$
from the ``ground semifield" $\RR_{>0}$ to an arbitrary semifield $K$ (see Introduction).
To do this, we define $\LL^{u,v}(K)$ as the set of all tuples
$\tt=(t^{\ii})_{\ii \in R(u,v)}$,
where each $t^{\ii}=(t^{\ii}_1,\dots,t^{\ii}_m)$
is a  ``vector'' in $K^m$ (with $m = \l(u) + \l(v)$),
and these vectors satisfy the
relations $t^{\ii'} = \tilde R_\ii^{\ii'} (t^\ii)$ for all $\ii, \ii' \in R(u,v)$.

\begin{example}
\label{ex:real Luv}
{\rm By Proposition~\ref{pr:Luv coordinates}, the map $\LL^{u,v}(\RR_{> 0}) \to G$
given by $\tt \mapsto x_{\ii}(t^\ii)$ is well-defined, and it is a bijection
between $\LL^{u,v}(\RR_{> 0})$ and $L^{u,v}_{> 0}$. }
\end{example}

\begin{example}
\label{ex:tropical Luv Lusztig}
{\rm By Proposition~\ref{pr:characterization of BB} and
Theorem~\ref{th:geometric lifting} (ii), the map
that sends every canonical basis vector $b \in \BB$ to a tuple $\tt$
with $t^\ii = b_{\ii}^{-1}(b)$ is a bijection
between $\BB$ and the set of all $\tt \in \LL^{e,\wnot}(\ZZ_{\rm trop})^\vee$
such that $t^\ii \in \ZZ_{\geq 0}^m$ for all $\ii \in R(\wnot)$
(as usual in this paper, $\LL^{e,\wnot}(K)^\vee$ stands for the set defined in the same way as
$\LL^{e,\wnot}(K)$ but for the Langlands dual group).}
\end{example}

\begin{example}
\label{ex:tropical Luv strings}
{\rm By Proposition~\ref{pr:strings} and
Theorem~\ref{th:geometric lifting} (ii), the map
that sends every dual canonical basis vector $b^* \in \BB^{\rm dual}$ to a tuple $\tt$
with $t^{-\ii} = c_{\ii}(b^*)$ is a bijection
between $\BB^{\rm dual}$ and the set of all $\tt \in \LL^{\wnot,e}(\ZZ_{\rm trop})^\vee$
such that $t^{-\ii} \in C_\ii (\ZZ)$ for all $\ii \in R(\wnot)$.}
\end{example}

Using Examples~\ref{ex:tropical Luv Lusztig} and \ref{ex:tropical Luv strings},
we shall identify $\BB$ (resp. $\BB^{\rm dual}$) with a part of $\LL^{e,\wnot}(\ZZ_{\rm trop})^\vee$
(resp. $\LL^{\wnot,e}(\ZZ_{\rm trop})^\vee$).
It is easy to see that the correspondence $b \mapsto b^{\rm dual}$ extends to a piecewise-linear map
$\LL^{e,\wnot}(\ZZ_{\rm trop})^\vee \to \LL^{\wnot,e}(\ZZ_{\rm trop})^\vee$.
Our next goal is to find a geometric lifting for this map.
To do this, we consider the following modification of the twist maps
$\psi^{\wnot,e}$ and $\psi^{e,\wnot}$ in Corollary~\ref{cor:twist special}.
Let us define the maps $\eta^{\wnot,e}$ and $\eta^{e,\wnot}$ by setting
(see (\ref{eq:tau-wnot}))
\begin{equation}
\label{eq:eta1}
\eta^{\wnot,e} = \tau_{\wnot} \circ \psi^{\wnot,e},
\quad  \eta^{e, \wnot} = \psi^{e,\wnot} \circ \tau_{\wnot} \ ;
\end{equation}
an easy calculation shows that these maps are given by
(see (\ref{eq:T}))
\begin{equation}
\label{eq:eta2}
\eta^{\wnot,e} (x) = [(\overline {\wnot} x^T)^{-1}]_+,
\quad  \eta^{e, \wnot} (x) = ([\overline {\wnot}^{\ -1} x^T]_0)^{-1}
([\overline {\wnot}^{\ -1} x^T]_-)^{\ -1} \ .
\end{equation}
Theorem~\ref{th:psi-regularity} immediately implies the following.

\begin{corollary}
\label{cor:eta-regularity}
The map $\eta^{\wnot,e}$ is a biregular isomorphism
between $L^{\wnot,e}$ and  $L^{e, \wnot}$, and it restricts to a bijection
between $L^{\wnot,e}_{> 0}$ and  $L^{e, \wnot}_{>0}$.
The inverse map is $\eta^{e,\wnot}$.
\end{corollary}

Now we are ready to state our geometric lifting result.

\begin{theorem}
\label{th:twist special}
The map $\eta^{\wnot,e}: (L^{\wnot,e}_{> 0})^\vee \to (L^{e,\wnot}_{> 0})^\vee$
is a geometric lifting of $b \mapsto b^{\rm dual}$.
In other words, for any $\ii, \ii' \in R(\wnot)$, we have
\begin{equation}
\label{eq:geometric dual}
c_{\ii} (b_{\ii'} (t')^{\rm dual}) =
[(x_{-\ii}^{\ -1} \circ \eta^{e,\wnot} \circ x_{\ii'})^\vee (t')]_{\rm trop}
\end{equation}
(as before, the superscript $^\vee$ means that the corresponding
varieties and maps are related to the group ${}^L G$).
\end{theorem}

Theorems~\ref{th:geometric lifting} and  \ref{th:twist special}
play crucial role in our proofs of the results in Sections~\ref{sec:LR} and
\ref{sec:can}; in fact, they allow us to deduce combinatorial properties
of the canonical basis from the properties of totally positive varieties.
The appearance of $\ii$-trails is explained by the following evaluation
of generalized minors.

\begin{theorem}
\label{th:minors-trails}
Let $\gamma$ and $\delta$ be two weights in the $W$-orbit of the same fundamental weight $\omega_i$
of $\gg$, and let $\ii = (i_1, \dots, i_m)$ be any sequence of indices from $[1,r]$.

\noindent {\rm (i)} $\Delta_{\gamma, \delta} (x_\ii (t_1, \dots, t_m))$
is a positive integer linear combination of the monomials $t_1^{c_1 (\pi)} \cdots t_m^{c_m (\pi)}$
for all $\ii$-trails $\pi$ from $\gamma$ to $\delta$ in $V_{\omega_i}$.

\smallskip

\noindent {\rm (ii)} $\Delta_{\gamma, \delta} (x_{-\ii} (t_1, \dots, t_m))$
is a positive integer linear combination of the monomials
$t_1^{d_1 (\pi)} \cdots t_m^{d_m (\pi)}$
for all $\ii$-trails $\pi$ from $-\gamma$ to $-\delta$ in $V_{\omega_{i*}}$.
\end{theorem}

Theorem~\ref{th:minors-trails} has an important ``tropical" corollary.

\begin{corollary}
\label{cor:minors-trails-tropical}
In the situation of Theorem~\ref{th:minors-trails}, we have
$$[\Delta_{\gamma, \delta} (x_\ii (t_1, \dots, t_m))]_{\rm trop} =
\min_\pi \ (\sum_{k=1}^m c_k (\pi) t_k) \ ,$$
$$[\Delta_{\gamma, \delta} (x_{-\ii} (t_1, \dots, t_m))]_{\rm trop} =
 \min_{\pi'} \ (\sum_{k=1}^m d_{k} (\pi') t_k) \ ,$$
where $\pi$ (resp. $\pi'$) runs over all $\ii$-trails
from $\gamma$ to $\delta$ in $V_{\omega_i}$
(resp. from $- \gamma$ to $- \delta$ in $V_{\omega_{i^*}}$).
\end{corollary}

\subsection{Geometric lifting of Kashiwara's crystals}
\label{sec:crystals}
Let us recall the crystal operators $\tilde f_i: \BB \to \BB$ introduced
in Proposition~\ref{pr:crystal BB}.
In this section we compute a geometric lifting of the twisted operators
$\tilde f_i^\iota : \BB \to \BB$
defined by $\tilde f_i^\iota (b) = (\tilde f_i (b^\iota))^\iota$.
In view of Proposition~\ref{pr:Lusztig parametrizations} (iii)),
for every positive integer $n$, the operator $(\tilde f_i^\iota)^n$
acts as follows on Lusztig parameters corresponding
to any reduced word $\ii' \in R(\wnot)$ with $i'_m = i^*$:
$$(\tilde f_i^\iota)^n (b_{\ii'}(t_1, \dots, t_{m-1}, t_m)) =
b_{\ii'}(t_1, \dots, t_{m-1}, t_m+n) \ .$$
We consider the following geometric counterpart of this operator:
for any $c > 0$ define the bijection
$F_i^{(c)}: L^{e,\wnot}_{> 0} \to L^{e,\wnot}_{> 0}$ by
\begin{equation}
\label{eq:geom crystal Lusztig}
F_i^{(c)} (x_{\ii'}(t'_1, \dots, t'_{m-1}, t'_m)) =
x_{\ii'}(t'_1, \dots, t'_{m-1}, c t'_m) \ ,
\end{equation}
where $\ii' \in R(\wnot)$ ends with $i'_m = i^*$.
If $\ii' \in R(\wnot)$ does not end with $i^*$ then we do not have
a nice formula for $F_i^{(c)} (x_{\ii'}(t'_1, \dots, t'_m))$.
Our next result shows that such a formula exists for the bijection
$L^{\wnot,e}_{> 0} \to L^{\wnot,e}_{> 0}$ obtained by transferring $F_i^{(c)}$
with the help of the bijection $\eta^{\wnot,e}: L^{\wnot,e}_{> 0} \to L^{e, \wnot}_{>0}$
(see Corollary~\ref{cor:eta-regularity}).

\begin{theorem}
\label{th:geom crystal}
Let $\ii$ be a reduced word for $R(\wnot)$, and let
$T_k = t_k^{-1} \prod_{l > k} t_l^{- a_{i_l, i_k}}$ for $k = 1, \dots, m$.
Then
$$(\eta^{e,\wnot} \circ F_i^{(c)} \circ \eta^{\wnot,e}) (x_{- \ii}(t_1, \dots, t_m)) =
x_{- \ii}(\tilde t_1, \dots, \tilde t_m) \ ,$$
where $\tilde t_k = t_k$ unless $i_k = i$, and
\begin{equation}
\label{eq:geom crystal}
\displaystyle{\tilde t_k=t_k \frac
{\sum \limits_{l<k:i_l=i} T_l + c \sum \limits_{l\geq k:i_l=i} T_l}
{\sum \limits_{l\leq k:i_l=i} T_l + c \sum \limits_{l> k:i_l=i} T_l}} \ ,
\end{equation}
whenever $i_k = i$.
\end{theorem}

In view of Theorem~\ref{th:twist special}, one obtains an explicit formula for the action of
$(\tilde f_i^\iota)^n$ on $\BB$ in terms of the string parameters
by tropicalizing (\ref{eq:geom crystal}) (and passing from $G$ to ${}^L G$ as usual).

\begin{corollary}
\label{cor:Kashiwara crystals}
For a reduced word $\ii \in R(\wnot)$, define the linear forms
$T^\vee_k:\ZZ^m\to \ZZ^m$ for $k=1, \dots, m$ by
$T_k^\vee(t)= -t_k - \sum_{l>k} a_{i_k,i_l} t_l$.
If $c_\ii (b^{\rm dual}) = (t_1, \dots, t_m)$ then
$c_\ii ((\tilde f_i^\iota)^n (b)^{\rm dual}) = (\tilde t_1, \dots, \tilde t_m)$,
where $\tilde t_k = t_k$ unless $i_k = i$, and
$$\tilde t_k = t_k+ \min \left(\min \limits_{l<k:i_l=i} T^\vee_l(t),
n + \min \limits_{l\geq k:i_l=i} T^\vee_l(t)\right)$$
$$- \min \left(\min\limits_{l\leq k:i_l=i} T^\vee_l(t),
n + \min\limits_{l> k:i_l=i} T^\vee_l(t)\right)$$
whenever $i_k = i$ (with the agreement that minimum over the empty set is $+ \infty$).
\end{corollary}

\begin{remark}
{\rm Theorem~\ref{th:geom crystal} is a starting point of a new
concept of \emph{geometric crystals} introduced and developed by one of the
authors (A.B.) in a joint work in progress with D.~Kazhdan.}
\end{remark}

\subsection{Pl\"ucker models}
\label{sec:Plucker models}
Following \cite[Section~4]{bz2}, we now consider the ``variety"
$\MM^{\wnot} (K)$  of all tuples $(M_{\omega_i,\gamma})$ of elements of the ground
semifield $K$ (for all $i \in [1,r]$, and $\gamma \in W \omega_i$)
satisfying the relations in Proposition~\ref{pr:minors-Plucker} (with $u = e$, and
each generalized minor $\Delta_{\omega_i,\gamma}$ replaced with $M_{\omega_i,\gamma}$).
We shall show that each of the varieties $\LL^{e,\wnot}(K)$ and $\LL^{\wnot,e}(K)$
is in a natural bijection with a part of $\MM^{\wnot} (K)$.
To define these bijections, we use Theorem~\ref{th:minors-trails}, 
which assures
that both $t \mapsto \Delta_{\gamma, \delta} (x_\ii (t_1, \dots, t_m))$ and
$t \mapsto \Delta_{\gamma, \delta} (x_{-\ii} (t_1, \dots, t_m))$ are well-defined
mappings $K^m \to K \cup \{0\}$ for any semifield $K$.

\begin{theorem}
\label{th:Plucker model for tp cells}
{\rm (i)} For every semifield $K$, the correspondence
$$p^+: \tt \mapsto (M_{\omega_i,\gamma} = \Delta_{\omega_i,\gamma}(x_\ii (t^\ii)))$$
(where $\ii \in R(\wnot)$) is an embedding of $\LL^{e, \wnot}(K)$ into $\MM^{\wnot} (K)$.
The image of $p^+$ consists of all tuples
$(M_{\omega_i,\gamma}) \in  \MM^{\wnot} (K)$ such that
$M_{\omega_i,\omega_i} = 1$ for all $i$.

\noindent {\rm (ii)} The map
$p^-: \tt \mapsto (M_{\omega_i,\gamma} = \Delta_{\gamma, \omega_i}(x_{-\ii} (t^\ii)))$
is an embedding of $\LL^{\wnot, e}(K)$ into $\MM^{\wnot} (K)$.
The image of $p^-$ consists of all tuples
$(M_{\omega_i,\gamma}) \in  \MM^{\wnot} (K)$ such that
$M_{\omega_i, \wnot \omega_i} = 1$ for all $i$.
\end{theorem}

We will refer to the maps $p^+$ and $p^-$ as \emph{Pl\"ucker models} of
$\LL^{e, \wnot}(K)$ and $\LL^{\wnot, e}(K)$.

Specializing Theorem~\ref{th:Plucker model for tp cells} to the
tropical semifield $K = \ZZ_{\rm trop}$, we obtain embeddings
$p^+_{\rm trop} : \BB \to \MM^{\wnot} (\ZZ_{\rm trop})^\vee$
and $p^-_{\rm trop} : \BB^{\rm dual} \to  \MM^{\wnot} (\ZZ_{\rm trop})^\vee$.
(As before, the variety $\MM^{\wnot} (K)^\vee$
corresponds to the Langlands dual group ${}^L G$.)
Our next task is to describe the images of these embeddings.
To do this, we notice that for every $\gamma \in W \omega_i \setminus \{\omega_i\}$,
there is a naturally defined function $M_{s_i \omega_i, \gamma}: \MM^{\wnot} (K) \to K$;
it can be obtained as a subtraction-free expression in the components
$M_{\omega_j, \delta}$ by iterating the identity (\ref{eq:minors-Dodgson})
(with $\Delta$ replaced by $M$) and using the boundary condition
$M_{s_i \omega_i, \omega_i} = 0$.
An explicit formula for $M_{s_i \omega_i, \gamma}$ can be given as
follows: let $\gamma = u \omega_i = s_{i_1} \cdots s_{i_l}
\omega_i$, where $(i_1, \dots, i_l)$ is a reduced word for $u \in
W$ such that $i_l = i$; then we have
\begin{equation}
\label{eq:subflag minors}
M_{s_i \omega_i, \gamma} = M_{\omega_i, \gamma} \doublesubscript{\sum}{k \leq l}{i_k = i}
\frac{1}{M_{\omega_i, s_{i_1} \cdots s_{i_k} \omega_i}
M_{\omega_i, s_{i_1} \cdots s_{i_{k-1}} \omega_i}} \prod_{j \neq i}
M_{\omega_j, s_{i_1} \cdots s_{i_k} \omega_j}^{- a_{ji}} \ .
\end{equation}

\begin{theorem}
\label{th:Plucker tropical}
{\rm (i)} The image of the embedding $p^+_{\rm trop} : \BB \to \MM^{\wnot} (\ZZ_{\rm trop})^\vee$
consists of all integer tuples
$(M_{\omega_i^\vee,\gamma}) \in  \MM^{\wnot} (\ZZ_{\rm trop})^\vee$ such that
$M_{\omega_i^\vee,\omega_i^\vee} = 0$ and $M_{\omega_i^\vee,s_i \omega_i^\vee} \geq 0$ for all $i$.

\smallskip

\noindent {\rm (ii)} The image of the embedding
$p^-_{\rm trop} : \BB^{\rm dual} \to  \MM^{\wnot} (\ZZ_{\rm trop})^\vee$
consists of all integer tuples
$(M_{\omega_i^\vee,\gamma}) \in  \MM^{\wnot} (\ZZ_{\rm trop})^\vee$ such that
$M_{\omega_i^\vee, \wnot \omega_i^\vee} = 0$ and
$M_{s_i \omega_i^\vee, \wnot \omega_i^\vee} \geq 0$ for all $i$.
\end{theorem}

The tropical Pl\"ucker models just constructed allow us to give
two ``universal"  polyhedral expressions for the tensor product multiplicities.

\begin{theorem}
\label{th:LR-Plucker-Lusztig}
For any three dominant weights $\lambda, \mu, \nu$ for $\gg$, the multiplicity
$c_{\lambda, \nu}^\mu$ is equal to the number of integer tuples
$(M_{\omega_i^\vee,\gamma}) \in \MM^{\wnot} (\ZZ_{\rm trop})^\vee$
satisfying the following conditions for any $i \in [1,r]$:

\smallskip

\noindent {\rm (0)} $M_{\omega_i^\vee, \omega_i^\vee} = 0$;

\smallskip

\noindent {\rm (1)} $M_{\omega_i^\vee, s_i \omega_i^\vee} \geq 0$;

\smallskip

\noindent {\rm (2)} $M_{\omega_i^\vee, \wnot \omega_i^\vee} =
(\lambda + \nu - \mu)(\omega_i^\vee)$;

\smallskip

\noindent {\rm (3)} $M_{s_i \omega_i^\vee, \wnot \omega_i^\vee}  \geq
(s_i \lambda + \nu - \mu)(\omega_i^\vee)$;

\smallskip

\noindent {\rm (4)} $M_{\omega_i^\vee, \wnot s_i \omega_i^\vee} \geq
(\lambda + s_i \nu - \mu)(\omega_i^\vee)$.
\end{theorem}

\begin{theorem}
\label{th:LR-Plucker-strings}
For any three dominant weights $\lambda, \mu, \nu$ for $\gg$, the multiplicity
$c_{\lambda, \nu}^\mu$ is equal to the number of integer tuples
$(M_{\omega_i^\vee,\gamma}) \in \MM^{\wnot} (\ZZ_{\rm trop})^\vee$
satisfying the following conditions for any $i \in [1,r]$:

\smallskip

\noindent {\rm (0)} $M_{\omega_i^\vee, \wnot \omega_i^\vee} = 0$;

\smallskip

\noindent {\rm (1)} $M_{s_i \omega_i^\vee, \wnot \omega_i^\vee} \geq 0$;

\smallskip

\noindent {\rm (2)} $M_{\omega_i^\vee, \omega_i^\vee} =
- (\lambda + \nu - \mu)(\omega_i^\vee)$;

\smallskip

\noindent {\rm (3)} $M_{\omega_i^\vee, \wnot s_i \omega_i^\vee} \geq
- \lambda (\alpha_{i^*}^\vee)$;

\smallskip

\noindent {\rm (4)} $M_{\omega_i^\vee, s_i \omega_i^\vee} \geq
- (s_i \lambda + \nu - s_i \mu)(\omega_i^\vee)$.
\end{theorem}

\section{Proofs of results in Section~\ref{sec:geometric}}
\label{sec:geometric proofs}

\noindent {\bf Proof of Proposition~\ref{pr:Luv equations}.}
Any element of the Bruhat cell $B u B$ can be written as $x = n_1 \overline u a n_2$ with
$a \in H$ and $n_1, n_2 \in N$.
Here the element $a$ is uniquely determined by $x$, and we have $a = [{\overline u}^{\ -1} x]_0$.
Thus the condition that $x \in N \overline u N$ is equivalent to $[{\overline u}^{\ -1} x]_0 = 1$.
This in turn is equivalent to the condition that
$\Delta_{u \omega_i, \omega_i} (x) = ([{\overline u}^{\ -1} x]_0)^{\omega_i} = 1$
for all $i \in  [1,r]$.
\endproof

\noindent {\bf Proof of Proposition~\ref{pr:Luv coordinates}.}
In view of \cite[Theorems 1.2, 1.3]{fz}, we only need to prove that
the image $x_\ii (\CC_{\neq 0}^m)$ is contained in $N \overline u N$.
Trivially, $x_i (t) \in N$ for any $i \in [1,r]$; using the commutation relation
\cite[(2.13)]{fz}, we also see that
$$x_{- i} (t) = y_i (t) t^{- \alpha_i^\vee} = x_i (t^{-1}) \overline {s_i} x_i (t)
\in N  \overline {s_i} N \ .$$
Using induction on $\l (u) + \l (v)$, it only remains to show that
$\overline {s_i} N \overline {w}  \subset N \overline {s_i w} N$
for any $w \in W$ such that $\l (s_i w) = \l (w) + 1$.
For any $n \in N$, we have
$$\overline {s_i} n \overline {s_i}^{\ -1} = n' y_i (t)$$
for some $n' \in N$ and $t \in \CC$.
Furthermore, the condition $\l (s_i w) = \l (w) + 1$ implies that
$$n'' := \overline {s_i w}^{\ -1} y_i (t) \overline {s_i w} \in N \ .$$
Therefore, we obtain
$$\overline {s_i}  n \overline {w} = n' y_i (t) \overline {s_i w} =
n' \overline {s_i w} n'' \in N \overline {s_i w} N \ ,$$
as required.
\endproof

\noindent {\bf Proof of Theorem~\ref{th:psi-regularity}.}
Comparing (\ref{eq:psi-u,v-x}) with \cite[Definition~1.5]{fz}, it is easy to show that
\begin{equation}
\label{eq:two twists}
\psi^{u,v} (x) = (\zeta^{u,v} (x))^T = \zeta^{u^{-1}, v^{-1}} (x^T)
\end{equation}
for any $x \in L^{u,v}$; here $\zeta^{u,v}$ is a biregular isomorphism
between $G^{u,v}$ and $G^{u^{-1}, v^{-1}}$ (see
\cite[Theorem~1.6]{fz}),
and $x \mapsto x^T$ is the ``transpose" antiautomorphism of $G$
given by (\ref{eq:T}).
It is also clear from (\ref{eq:psi-u,v-x}) that $\psi^{u,v} (x) \in N \overline v N$.
Therefore, $\psi^{u,v}$ sends $L^{u,v}$ to $L^{v,u}$, and
our theorem becomes a consequence of \cite[Theorems~1.6, 1.7]{fz}.
\endproof

\noindent {\bf Proof of Theorem~\ref{th:t-through-x}.}
In order to distinguish our present notation from that in \cite{fz}, let us denote
by $\tilde x_{\ii}$ the product map $H \times \CC_{\neq 0}^m \to G^{u,v}$ defined
in \cite[(1.3)]{fz}.
Using (\ref{eq:xnegative}) and commutation relations \cite[(2.5)]{fz}, we can rewrite
$x = x_\ii (t_1, \dots, t_m)$ as $x = \tilde x_{\ii} (a; t'_1, \dots, t'_m)$,
where
$$a = \prod_{k: i_k < 0} t_k^{- \alpha_{-i_k}^\vee} \ ,$$
and
\begin{equation}
\label{eq:t'k}
t'_k = t_k^{\varepsilon (i_k)} \cdot \doublesubscript{\prod}{l > k}
{i_l < 0} t_l^{\varepsilon (i_k) a_{|i_l|,|i_k|}} \ ;
\end{equation}
here $\varepsilon (i)$ denotes the sign of $i$, i.e.,
$\varepsilon (i) = \pm 1$ for $i \in \pm [1,r]$.

Now to prove our theorem, one only has to substitute into (\ref{eq:t'k}) the expressions for
the $t_l$ given by (\ref{eq:negative t-through-x}) and (\ref{eq:positive t-through-x}),
and to verify that the resulting expression for $t'_k$ agrees with \cite[(1.18)]{fz}.
This is done by a straigtforward check that we omit (notice that the function
$\Delta_{l, \ii} (x')$ in \cite[(1.18)]{fz} is our present $M_l (x)$,
and that $\Delta_{m+j, \ii} (x') = \Delta_{v \omega_j, \omega_j} (\psi^{u,v} (x)) = 1$
for any $j \in [1,r]$, since $\psi^{u,v} (x) \in L^{v,u}$ by
Theorem~\ref{th:psi-regularity}).
\endproof

\noindent {\bf Proof of Proposition~\ref{pr:t1 special}.}
To compute minors on the right-hand sides of (\ref{eq:negative t1 special})
and (\ref{eq:positive t1 special}), we shall use \cite[Lemma~6.4 (b)]{bz2}
which says that
\begin{equation}
\label{eq:highest minor}
\Delta_{\omega_i, w^{-1} \omega_i} (x_{j_1} (t_1) \cdots x_{j_l} (t_l))
= \prod_{k = 1}^l t_k^{\omega_i (s_{j_1} \cdots s_{j_{k-1}} \alpha_{j_k}^\vee)}
\end{equation}
for any $i \in [1,r]$, $w \in W$, and $\jj = (j_1, \dots, j_l) \in R(w)$.
This formula (with $j = (i_1, \dots, i_m)$) directly applies to the two minors in
the numerators in (\ref{eq:positive t1 special}).
To compute the second denominator in (\ref{eq:positive t1 special}), we notice that
it is equal to
$$\Delta_{\omega_{i_m^*}, s_{i_m} \wnot \omega_{i_m^*}} (x) =
\Delta_{\omega_{i_m^*}, s_{i_m} \wnot \omega_{i_m^*}} (x_{i_1} (t_1) \cdots x_{i_{m-1}} (t_{m-1}))
$$
(since $\overline {s_{i_m} \wnot}^{\ -1} x_{i_m} (t_m) \overline {s_{i_m} \wnot} \in N$),
and so it is given by (\ref{eq:highest minor}) for $\jj = (i_1, \dots, i_{m-1})$.
This implies the formula for $t_m$ in (\ref{eq:positive t1 special}).

To compute the first denominator in (\ref{eq:positive t1 special}), we
shall use
the antiautomorphism $\tau_{\wnot}$ of $G$ introduced in (\ref{eq:tau-wnot}).
It is easy to see that
\begin{equation}
\label{eq:tau on products}
\tau_{\wnot} (x_{i_1} (t_1) \cdots x_{i_m} (t_m)) = x_{i^*_m} (t_m) \cdots x_{i^*_1} (t_1)
\end{equation}
for any sequence $(i_1, \dots, i_m)$ of indices from $[1,r]$.
Using (\ref{eq:minors-tau}), we can now rewrite the first denominator in
(\ref{eq:positive t1 special}) as
$$\Delta_{s_{i_1} \omega_{i_1}, \wnot \omega_{i_1}} (x) =
\Delta_{\omega_{i_1}, \wnot s_{i_1} \omega_{i_1}} (x_{i^*_m} (t_m) \cdots x_{i^*_1} (t_1))$$
$$= \Delta_{\omega_{i_1}, \wnot s_{i_1} \omega_{i_1}} (x_{i^*_m} (t_m) \cdots x_{i^*_2} (t_2))
\ ,$$
and then compute it by using (\ref{eq:highest minor}) with $\jj = (i^*_m, \dots, i^*_{2})$.
This implies the formula for $t_1$ in (\ref{eq:positive t1 special}).

To prove (\ref{eq:negative t1 special}), we shall use the following lemma.

\begin{lemma}
\label{lem:negative minors}
For any sequence $(i_1, \dots, i_m)$  of indices from $[1,r]$,
and any generalized minor $\Delta_{\gamma, \delta}$, we have
\begin{equation}
\label{eq:negative-to-positive}
(x_{-i_1} (t_1) \cdots x_{-i_m} (t_m))^T =
t_1^{- \alpha_{i_1}^\vee} \cdots t_m^{- \alpha_{i_m}^\vee}
x_{i_m} (t'_m) \cdots x_{i_1} (t'_1) \ ,
\end{equation}
and
\begin{equation}
\label{eq:negative minors}
\Delta_{\gamma, \delta} (x_{-i_1} (t_1) \cdots x_{-i_m} (t_m)) =
\prod_{k=1}^m t_k^{- \delta (\alpha_{i_k}^\vee)} \cdot
\Delta_{\delta, \gamma} (x_{i_m} (t'_m) \cdots x_{i_1} (t'_1)) \ ,
\end{equation}
where the $t'_k$ are given by
\begin{equation}
\label{eq:t-to-t'}
t'_k = t_k \cdot \prod_{l < k} t_l^{a_{i_l, i_k}} \ .
\end{equation}
\end{lemma}

\proof
The equality (\ref{eq:negative-to-positive}) follows from
(\ref{eq:xnegative}) and commutation relations \cite[(2.5)]{fz}.
The equality (\ref{eq:negative minors}) is an immediate consequence
of (\ref{eq:negative-to-positive}) and (\ref{eq:minor properties}).
\endproof

Using Lemma~\ref{lem:negative minors}, we deduce the computation of the minors in
(\ref{eq:negative t1 special}) to the computation of ``transpose minors" evaluated
at the product $x_{i_m} (t'_m) \cdots x_{i_1} (t'_1)$.
The latter minors are computed in the same way as above,
and (\ref{eq:negative t1 special}) follows.
\endproof

\begin{remark}
{\rm Although we do not need it in this paper, we would like to
give an explicit formula for the inverse transformation in (\ref{eq:t-to-t'}):
\begin{equation}
\label{eq:t'-to-t}
t_k = t'_k \cdot \prod_{l < k}
{t'_l}^{\ - \alpha_{i_k} (s_{i_{k-1}} \cdots s_{i_{l+1}} \alpha_{i_l}^\vee)} \ .
\end{equation}
The fact that transformations in (\ref{eq:t-to-t'}) and (\ref{eq:t'-to-t})
are inverse of each other is proved by a straightforward calculation.}
\end{remark}

\section{Proofs of results in Section~\ref{sec:transition}}
\label{sec:totpos}

We start by recalling the classical
Tits theorem about reduced words in Coxeter groups.
We call a \emph{$d$-move} a transformation of a reduced word
(for some $w \in W$) that replaces $d$ consecutive entries
$i,j,i,j, \ldots$ by $j,i,j,i,\ldots$,
for some $i$ and $j$ such that $d$ is the order of~$s_i s_j\,$.
Note that, for given $i$ and $j$,
the value of $d$ can be determined from the Cartan matrix as follows:
if $a_{ij}a_{ji} = 0$ (resp.\ $1,2,3$), then $d=2$ (resp.\ $3,4,6$).
The Tits theorem says that every two reduced words for the same element of a Coxeter group
can be obtained from each other by a sequence of $d$-moves.

Applying this to the group $W \times W$,
we conclude that every two double reduced words $\ii, \ii' \in
R(u,v)$ can be obtained from each other by a sequence of the
following operations: \emph{positive} $d$-moves for the alphabet $[1,r]$,
\emph{negative} $d$-moves for the alphabet $[- 1, - r]$, and \emph{mixed} $2$-moves
that interchange two consecutive indices of opposite signs.
Furthermore, if $\ii$ and  $\ii'$ are related by one of these $d$-moves then the
transition map $\tilde R_{\ii}^{\ii'}$ only affects the $d$ entries involved in the move; and the
restriction of $\tilde R_{\ii}^{\ii'}$ to the corresponding segment of $t \in \RR_{> 0}^m$
is the ``local" transition map $\tilde R_{i,j,i, \ldots}^{j,i,j, \ldots}$.
In this section, we give explicit formulas for these local transition maps.

We shall write
$\tilde R_{i,j,i, \ldots}^{j,i,j, \ldots}(t_1, \ldots, t_d) = (p_1, \ldots, p_d)$;
thus the tuples $t_1, \ldots, t_d$ and $p_1, \ldots, p_d$ are related by
\begin{equation}
\label{eq:d-transition}
x_i (t_1) x_j (t_2) x_i (t_3) \cdots =
x_j (p_1) x_i (p_2) x_j (p_3) \cdots \ .
\end{equation}
The transition maps for positive $d$-moves (i.e., for $i, j \in [1,r]$) were computed in
\cite[Theorem~3.1]{bz2}.
For the convenience of the reader, let us reproduce these results here.

\begin{proposition}
\label{pr:positive transition}
Let $i, j \in [1,r]$, and let $d$ be the order of $s_i s_j$ in $W$.
Then the transition map in {\rm (\ref{eq:d-transition})} is given as follows:

\noindent {\rm (1) Type $A_1 \times A_1$:} if $a_{ij} = a_{ji} = 0$ then  $d = 2$, and
\begin{equation*}
p_1 = t_2, \,\, p_2 = t_1  \ .
\end{equation*}

\noindent {\rm (2) Type $A_2$:} if $a_{ij} = a_{ji} = -1$ then  $d = 3$, and
\begin{equation*}
p_1 = \frac{t_2 t_3}{t_1 + t_3}, \,\, p_2 = t_1 + t_3, \,\,
p_3 = \frac{t_1 t_2}{t_1 + t_3} \ .
\end{equation*}

\noindent {\rm (3) Type $B_2$:} if $a_{ij} = -2, \, a_{ji} = -1$ then  $d = 4$, and
\begin{equation*}
p_1 = \frac{t_2 t_3^2 t_4}{\pi_2}, \,\, p_2 = \frac{\pi_2}{\pi_1}, \,\,
p_3 = \frac{\pi_1^2}{\pi_2}, \,\, p_4 = \frac{t_1 t_2 t_3}{\pi_1} \ ,
\end{equation*}
where
$$\pi_1 = t_1 t_2 + (t_1+t_3) t_4, \,\,
\pi_2 = t_1^2 t_2 + (t_1+t_3)^2 t_4 \ . $$

\noindent {\rm (4) Type $G_2$:} if $a_{ij} = -3, \, a_{ji} = -1$ then  $d = 6$, and
\begin{eqnarray*}
\begin{array}{lc}
p_1 = \frac{t_2 t_3^3 t_4^2 t_5^3 t_6}{\pi_3}, \,
p_2 = \frac{\pi_3}{\pi_2}, \,
p_3 = \frac{\pi_2^3}{\pi_3 \pi_4}, \\[.1in]
p_4 = \frac{\pi_4}{\pi_1 \pi_2}, \,
p_5 = \frac{\pi_1^3}{\pi_4}, \,
p_6 = \frac{t_1 t_2 t_3^2 t_4 t_5}{\pi_1} \ ,
\end{array}
\end{eqnarray*}
where
\begin{eqnarray*}
\begin{array}{rcl}
\pi_1 &=& t_1 t_2 t_3^2 t_4 + t_1 t_2 (t_3+t_5)^2 t_6
+ (t_1+t_3)t_4 t_5^2 t_6 \ , \\[.1in]
\pi_2 &=& t_1^2 t_2^2 t_3^3 t_4 + t_1^2 t_2^2 (t_3+t_5)^3 t_6
+(t_1+t_3)^2 t_4^2 t_5^3 t_6 \\[.1in]
&+& t_1 t_2 t_4 t_5^2 t_6
(3t_1 t_3 + 2t_3^2 +  2t_3 t_5 + 2t_1 t_5) \ , \\[.1in]
\pi_3 &=& t_1^3 t_2^2 t_3^3 t_4 + t_1^3 t_2^2 (t_3+t_5)^3 t_6
+(t_1+t_3)^3 t_4^2 t_5^3 t_6 \\[.1in]
&+& t_1^2 t_2 t_4 t_5^2 t_6
(3t_1 t_3 + 3t_3^2+ 3t_3 t_5  + 2t_1 t_5) \ , \\[.1in]
\pi_4 &=& t_1^2 t_2^2 t_3^3 t_4
\bigl(t_1 t_2 t_3^3 t_4 + 2t_1 t_2 (t_3+t_5)^3 t_6 +
(3t_1 t_3 + 3t_3^2+ 3t_3 t_5  + 2t_1 t_5) t_4 t_5^2 t_6 \bigr) \\[.1in]
&+& t_6^2 \bigl(t_1 t_2 (t_3+t_5)^2 + (t_1+t_3) t_4 t_5^2 \bigr)^3 \ .
\end{array}
\end{eqnarray*}

\noindent {\rm (5)} Furthermore, in each of the cases {\rm (1)--(4)} above,
if we interchange $a_{ij}$ with $a_{ji}$ then the corresponding transition
map in {\rm (\ref{eq:d-transition})} is obtained from the given one by
the transformation $p_k \to p_{d+1-k}, \, t_k \to t_{d+1-k}$.
\end{proposition}

The transition maps for mixed $2$-moves are given by the following proposition
which is an immediate consequence of commutation relations
\cite[(2.5), (2.11)]{fz}.

\begin{proposition}
\label{pr:mixed 2-move transition}
For any  $i, j \in [1,r]$, we have $x_j (t_1) x_{-i} (t_2) = x_{-i} (p_1) x_j (p_2)$,
where
\begin{equation*}
p_1 = t_2, \,\, p_2 = t_1  t_2^{a_{ij}}
\end{equation*}
for $i \neq j$, and
\begin{equation*}
\frac{1}{p_1} = t_1 + \frac{1}{t_2}, \,\,
\frac{1}{p_2} = \frac{1}{t_2} \left(1 + \frac{1}{t_1 t_2} \right)
\end{equation*}
for $i = j$.
\end{proposition}

Finally, the transition maps for negative $d$-moves are given as follows.

\begin{proposition}
\label{pr:negative transition}
Let $i, j \in [-1, -r]$, and let $d$ be the order of $s_{|i|} s_{|j|}$ in $W$.
Then the transition map in {\rm (\ref{eq:d-transition})} is given as follows:

\noindent {\rm (1) Type $A_1 \times A_1$:} if $a_{|i|, |j|} = a_{|j|, |i|} = 0$ then  $d = 2$, and
\begin{equation*}
p_1 = t_2, \,\, p_2 = t_1  \ .
\end{equation*}

\noindent {\rm (2) Type $A_2$:} if $a_{|i|, |j|} = a_{|j|, |i|} = -1$ then  $d = 3$, and
\begin{equation*}
\frac{1}{p_1} = \frac{1}{t_3} + \frac{t_1}{t_2}, \,\, p_2 = t_1 t_3, \,\,
p_3 = t_1 + \frac{t_2}{t_3} \ .
\end{equation*}

\noindent {\rm (3) Type $B_2$:} if $a_{|i|, |j|} = -1, \, a_{|j|, |i|} = -2$ then  $d = 4$, and
\begin{eqnarray*}
\begin{array}{lc}
\frac{1}{p_1} = \frac{t_1}{t_2} + \frac{t_2}{t_3} + \frac{1}{t_4}, , \,\,
\frac{1}{p_2} = \frac{1}{t_1} \left(\frac{t_2}{t_3} + \frac{1}{t_4} \right)^2 + \frac{1}{t_3}
\ , \\[.1in]
p_3 = t_2 +  t_1 t_4 + \frac{t_2^2 t_4}{t_3}, \,\,
p_4 = t_1 + t_3 (\frac{t_2}{t_3} + \frac{1}{t_4})^2 \ .
\end{array}
\end{eqnarray*}

\noindent {\rm (4) Type $G_2$:} if $a_{|i|, |j|} = -1, \, a_{|j|, |i|} = -3$ then  $d = 6$, and
\begin{eqnarray*}
\begin{array}{rcl}
\frac{1}{p_1} &=& \frac{t_1}{t_2} + t_3 (\frac{t_2}{t_3} + \frac{1}{t_4})^2 +
\frac{t_4}{t_5} + \frac{1}{t_6},
\\[.1in]
\frac{1}{p_2} &=& \frac{t_1}{t_3} + 2 t_3 \left(\frac{t_2}{t_3} + \frac{1}{t_4}\right)^3 +
\frac{1}{t_1}
\left(t_3 (\frac{t_2}{t_3} + \frac{1}{t_4})^2 + \frac{t_4}{t_5} + \frac{1}{t_6}\right)^3
\\[.1in]
&+& \frac{3 t_2 t_4}{t_3 t_5}  + \frac{3 t_2}{t_3 t_6} + \frac{3}{t_4 t_6} + \frac{2}{t_5},
\\[.1in]
p_5 &=& t_1 t_6 + t_3^2 t_6 (\frac{t_2}{t_3} + \frac{1}{t_4})^3
+ t_4 t_6 (\frac{t_4}{t_5} + \frac{1}{t_6})^2 + 2 t_2 + \frac{2 t_3}{t_4} +
\frac{3 t_2 t_4 t_6}{t_5} + \frac{2 t_3 t_6}{t_5}, \\[.1in]
p_6 &=& t_1 + t_3^2 (\frac{t_2}{t_3} + \frac{1}{t_4})^3
+ t_5 (\frac{t_4}{t_5} + \frac{1}{t_6})^3 +
\frac{3 t_2 t_4}{t_5}  + \frac{3 t_2}{t_6} + \frac{3 t_3}{ t_4 t_6} + \frac{2 t_3}{t_5} \ ;
\end{array}
\end{eqnarray*}
the two middle components $p_3$ and $p_4$ are determined from
two additional relations
$$ p_1 p_3 p_5 = t_2 t_4 t_6 \ , \quad p_2 p_4 p_6 = t_1 t_3 t_5 \ .$$

\noindent {\rm (5)} Furthermore, in each of the cases {\rm (1)--(4)} above,
if we interchange $a_{|i|, |j|}$ with $a_{|j|, |i|}$ then the corresponding transition
map in {\rm (\ref{eq:d-transition})} is obtained from the given one by
the transformation $p_k \to 1/p_{d+1-k}, \, t_k \to 1/t_{d+1-k}$.
\end{proposition}

\proof
Each of the formulas in Proposition~\ref{pr:negative transition}
follows from the corresponding formula in
Proposition~\ref{pr:positive transition} by applying the map
$x \mapsto x^T$ to both sides of (\ref{eq:d-transition}) and using
(\ref{eq:negative-to-positive}).
\endproof

\noindent {\bf Proof of Theorem~\ref{th:geometric lifting}.}
Part (i) of Theorem~\ref{th:geometric lifting} follows from the Tits
theorem and Propositions~\ref{pr:positive transition},
\ref{pr:mixed 2-move transition} and
\ref{pr:negative transition} since all rational expressions appearing there are subtraction-free.
As for part (ii), it is enough to check it for the rank 2 case when the
``geometric" transition maps $\tilde R_{i,j,i, \ldots}^{j,i,j, \ldots}$ are given by
Propositions~\ref{pr:positive transition} and \ref{pr:negative transition}.

The case of $A_1 \times A_1$ is obvious, and the case of $A_2$ is obtained by comparing
the expressions in (\ref{eq:Lustzig A2 transition}) and (\ref{eq:string A2 transition})
with the corresponding expressions in
Propositions~\ref{pr:positive transition} and \ref{pr:negative transition}.

The transition maps for string parametrizations for the types
$B_2$ (or $C_2$) and $G_2$ were found in \cite{lit,naka99}; our theorem
is then proved by direct comparison of these formulas with
the ones given by Proposition~\ref{pr:negative transition}.

For the Lusztig parametrizations, our proof is even less computational
(but still mysterious).
First of all, the statement for types $A_1 \times A_1$ and $A_2$
implies that it is true for any simply-laced type.
The transition maps for Lusztig
parametrizations for the type $B_2$ were found in \cite{lu92} using the
following strategy.
Let $a_{12} = -2$ and $a_{21} = -1$, i.e., $\alpha_2$ is the long simple root.
Lusztig (implicitly) claims that the transition map
$R_{2,1,2,1}^{1,2,1,2}$ for type $B_2$ is obtained from the transition
map $R_{1,3,2,1,3,2}^{2,1,3,2,1,3}$ for type $A_3$
(with the standard numeration of simple roots) by the following procedure:
$$R_{2,1,2,1}^{1,2,1,2}(t_1, t_2, t_3, t_4) = (p_1, p_2, p_3,p_4)$$
$$\Leftrightarrow R_{1,3,2,1,3,2}^{2,1,3,2,1,3}(t_1, t_1, t_2, t_3, t_3, t_4) =
(p_1, p_2,p_2,  p_3,p_4, p_4) \ . $$
Using our statement for type $A_3$, we see that
$R_{1,3,2,1,3,2}^{2,1,3,2,1,3}$ is the tropicalization of the
geometric transition map $\tilde R_{1,3,2,1,3,2}^{2,1,3,2,1,3}$.
But then the equality
$$\tilde R_{1,3,2,1,3,2}^{2,1,3,2,1,3}(t_1, t_1, t_2, t_3, t_3, t_4) =
(p_1, p_2,p_2,  p_3,p_4, p_4)$$
is easily seen to be equivalent to
$\tilde R_{1,2,1,2}^{2,1,2,1}(t_1, t_2, t_3, t_4) = (p_1, p_2, p_3,p_4)$,
where $\tilde R_{1,2,1,2}^{2,1,2,1}$ is the geometric transition map
for $B_2$, with the same convention as above: $\alpha_2$ is the long simple
root.
This proves our statement for the type $B_2$ (notice
that the geometric lifting of $R_{2,1,2,1}^{1,2,1,2}$ is
$\tilde R_{1,2,1,2}^{2,1,2,1} = (\tilde R_{2,1,2,1}^{1,2,1,2})^\vee$,
which explains the necessity of passing to the Langlands dual group).

Now for the type $G_2$ one can use the same argument, with $A_3$ replaced by
$D_4$.
This concludes our proof.
\endproof

\noindent {\bf Proof of Theorem~\ref{th:twist special}.}
Let us denote the left hand side of (\ref{eq:geometric dual}) by
$F_{\ii,\ii'}(t')$.
We consider each $F_{\ii,\ii'}$ as a map from $\ZZ_{\geq 0}^m$ to $\ZZ^m$.
The following properties of these maps are immediate from the definitions:

\smallskip

\noindent (1) $F_{\ii,\ii'}(0, \dots, 0) = (0, \dots, 0)$ for any $\ii, \ii' \in R(\wnot)$.

\smallskip

\noindent (2)
$F_{\ii,\ii''} = R_{-\ii'}^{-\ii} \circ F_{\ii',\ii''} =
F_{\ii,\ii'} \circ R_{\ii''}^{\ii'}$ for any $\ii, \ii', \ii'' \in R(\wnot)$.

\smallskip

The next property is less obvious:

\smallskip

\noindent (3) For any two reduced words $\ii$ and $\ii'$ for $\wnot$ such
that $i'_1 = i_1$, the first component of
$F_{\ii,\ii'}(t'_1, \dots, t'_m)$ is equal to $t'_1$, while all other components
only depend on $t'_2, \dots, t'_m$.

The statement about the first component follows from (\ref{eq:li = ci}) and
Proposition~\ref{pr:Lusztig parametrizations} (ii).
The statement about other components follows from properties of the crystal operators
$\tilde f_i: \BB \to \BB$ in Proposition~\ref{pr:crystal BB} (i).
Indeed, if $i'_1 = i_1$ then $F_{\ii,\ii'}$ commutes with the shift operator
$T_1$ which acts on $\ZZ^m$ by adding $1$ to the first component of a vector.

We now claim that the above properties uniquely determine the family of maps
$F_{\ii,\ii'}$.
More precisely: if a collection of maps $F_{\ii,\ii'}: \ZZ_{\geq 0}^m \to \ZZ^m$
satisfy properties (1)--(3) then $F_{\ii,\ii'}(t') = c_{\ii} (b_{\ii'} (t')^{\rm dual})$
for any $\ii, \ii' \in R(\wnot)$ and any $t' \in \ZZ_{\geq 0}^m$.

First of all, the equality $F_{\ii,\ii''} = F_{\ii,\ii'} \circ R_{\ii''}^{\ii'}$
in (2) implies that there exists a collection of maps $F_{\ii}: \BB \to \ZZ^m$
such that $F_{\ii,\ii'} = F_{\ii} \circ b_{\ii'}$ for any $\ii, \ii' \in R(\wnot)$.
It remains to show that $F_{\ii} (b) = c_{\ii} (b^{\rm dual})$ for any $\ii \in R(\wnot)$
and $b \in \BB$.
By (1), this is true for $b = 1$.
If $b \neq 1$ then, by Proposition~\ref{pr:crystal BB} (ii),
$b = \tilde f_i (b')$ for some $i \in [1,r]$ and $b' \in \BB$.
Pick any $\ii'$ with $i'_1 = i$.
Using induction on the weight of $b$, we can assume that
$F_{\ii'} (b') = c_{\ii'} ({b'}^{\rm dual})$.
By (3), this implies that $F_{\ii'} (b) = c_{\ii'} ({b}^{\rm dual})$.
Finally, the first equality in (2) implies that
$F_{\ii} (b) = c_{\ii} (b^{\rm dual})$ for any $\ii \in R(\wnot)$,
as required.

To complete the proof of Theorem~\ref{th:twist special}, it remains to show that
the functions $\tilde F_{\ii,\ii'}(t)$ given by the right hand side of
(\ref{eq:geometric dual}) satisfy the same properties (1)--(3).
To prove (1) notice that $[Q]_{\rm trop} (0, \dots, 0) = 0$
for any subtraction-free rational expression $Q$.
Property (2) follows from Theorem~\ref{th:geometric lifting} (ii).
Let us prove (3).

Fix two reduced words $\ii$ and $\ii'$ for $\wnot$ such that $i'_1 = i_1$.
Let $t'_1, \dots, t'_m \in \RR_{> 0}$, and let
$$x' = x_{\ii'} (t'_1, \dots, t'_m) \in L^{e,\wnot}_{> 0}, \,\,
x = \eta^{e,\wnot} (x') = x_{-\ii} (t_1, \dots, t_m) \in L^{\wnot,e}_{> 0} \ .$$
We need to show that $t_1 = t'_1$, and, for $k  = 2, \dots, m$, that $t_k$ does not
depend on $t'_1$.
Remembering (\ref{eq:eta1}) and combining Corollary~\ref{cor:factors special} (i)
with (\ref{eq:minors-tau}), we obtain
\begin{equation}
\label{eq:t'-eta-t}
t_k =  \frac{\Delta_{s_{i_1} \cdots s_{i_{k-1}} \omega_{i_k}, \wnot \omega_{i_k}}(x')}
{\Delta_{s_{i_1} \cdots s_{i_{k}} \omega_{i_k}, \wnot \omega_{i_k}}(x')}
\end{equation}
for any $k = 1, \dots, m$.
The equality $t_1 = t'_1$ now follows
by comparing (\ref{eq:t'-eta-t}) for $k=1$ with the first equality in
(\ref{eq:positive t1 special}).
Furthermore, if $k > 1$ then both minors in (\ref{eq:t'-eta-t})
are invariant under the transformation $x' \mapsto x_{i_1} (t) x$
for any $t$ since both elements $s_{i_{k-1}} \cdots s_{i_1}$ and
$s_{i_k} \cdots s_{i_1}$ send $\alpha_{i_1}$ to a negative root;
it follows that $t_k$ does not depend on $t'_1$, and we are done.
\endproof

\noindent {\bf Proof of Theorem~\ref{th:minors-trails}.}
First let us show that (ii) follows from (i).
Indeed, using (\ref{eq:negative minors}) and the first equality in (\ref{eq:minors-tau}),
we see that
$$\Delta_{\gamma, \delta} (x_{-\ii} (t_1, \dots, t_m)) =
\prod_{k=1}^m t_k^{- \delta (\alpha_{i_k}^\vee)} \cdot
\Delta_{-\gamma, -\delta} (x_{\ii} (t'_1, \dots, t'_m)) \ ,$$
where the $t'_k$ are given by (\ref{eq:t-to-t'}).
Computing the minor on the right hand side with the help of (i), we obtain the sum
of monomials corresponding to $\ii$-trails $\pi$ from $-\gamma$ to $-\delta$;
the exponent of $t_k$ in such a monomial is equal to
$$- \delta (\alpha_{i_k}^\vee) + c_k (\pi) + \sum_{l > k} a_{i_k, i_l} c_l (\pi)$$
$$= (- \delta + \sum_{l > k} c_l (\pi) \alpha_{i_l} + c_k (\pi) \frac{\alpha_{i_k}}{2})
(\alpha_{i_k}^\vee) \ .$$
Substituting $c_l (\pi) \alpha_{i_l} = \gamma_{l-1} - \gamma_l$ and
remembering the definition (\ref{eq:dk}), we conclude that the latter exponent is equal to
$d_k (\pi)$, as required.

For the proof of (i) we need a little preparation.
Consider the ring of regular functions $\CC[G]$ as a $G \times G$-representation
under the action $(g_1, g_2)f (x) = f(g_1^T x g_2)$.
We denote by $f \mapsto (u_1, u_2) f$ the corresponding action of $U(\gg) \times U(\gg)$,
where $U(\gg)$ is the universal enveloping algebra of $\gg$.
For every $f \in \CC[G]$, the function $f(x_{\ii} (t_1, \dots, t_m))$ is a polynomial
in $t_1, \dots, t_m$, and the coefficient of each monomial $t_1^{n_1} \cdots t_m^{n_m}$
is equal to $((1, e_{i_1}^{(n_1)} \cdots e_{i_m}^{(n_m)})f) (e)$, where $e$ stands for the identity
element of $G$, and $e_i^{(n)}$ stands for the divided power $e_i^n/n!$
(cf. \cite[Lemma~3.7.5]{bfz}).
If $f$ is (bi)-homogeneous of degree $(\gamma, \gamma')$ then $f(e)$ can be nonzero
only if $\gamma = \gamma'$.
It follows that if degree of $f$ is $(\gamma, \delta)$ then
$f(x_{\ii} (t_1, \dots, t_m))$ contains only monomials
$t_1^{n_1} \cdots t_m^{n_m}$ with $\sum_k n_k \alpha_{i_k} = \gamma - \delta$.

Returning to generalized minors, we notice that $\Delta_{\gamma, \delta}$ has degree
$(\gamma, \delta)$ (see (\ref{eq:minor properties})), and belongs to the submodule
$V_{\omega_i, \omega_i}$ of $\CC[G]$ generated by the highest weight vector
$\Delta_{\omega_i, \omega_i}$.
Furthermore, $\Delta_{\gamma, \delta}$ spans the weight subspace
$V_{\omega_i, \omega_i}(\gamma, \delta)$, and we also have $\Delta_{\gamma, \gamma}(e) = 1$.
It follows that the coefficient $c$ of $t_1^{n_1} \cdots t_m^{n_m}$ in
$\Delta_{\gamma, \delta}(x_{\ii} (t_1, \dots, t_m))$ can be found from the equality
$(1, e_{i_1}^{(n_1)} \cdots e_{i_m}^{(n_m)}) \Delta_{\gamma, \delta} = c \Delta_{\gamma, \gamma}$.
Applying the element $(\overline u^T, e) \in G \times G$ to both sides of this equality, we
see that
$(1, e_{i_1}^{(n_1)} \cdots e_{i_m}^{(n_m)}) \Delta_{\omega_i, \delta}
= c \Delta_{\omega_i, \gamma}$.
Remembering Definition~\ref{def:trails}, we see that
$\Delta_{\gamma, \delta}(x_{\ii} (t_1, \dots, t_m))$ consists precisely of the monomials
$t_1^{c_1 (\pi)} \cdots t_m^{c_m (\pi)}$
for all $\ii$-trails $\pi$ from $\gamma$ to $\delta$ in $V_{\omega_i}$.
It only remains to show that, for every such $\ii$-trail $\pi$, the
corresponding coefficient $c$ is a positive integer.

Let us consider the $U(\gg)$-module structure on $\CC[G]$ given by
$u f = (1,u)f$.
Under this action, $\Delta_{\omega_i, \omega_i}$ is a highest weight vector, and it generates
the submodule isomorphic to $V_{\omega_i}$.
Each minor $\Delta_{\omega_i, \gamma}$ is an extremal vector of weight $\delta$ in this
submodule normalized in such a way that $u \Delta_{\omega_i, \gamma} =
\Delta_{\omega_i, \omega_i}$ for some $u \in U(\gg)$ which is a monomial in the divided
powers $e_j^{(n)}$ (see \cite[Lemma~2.8]{fz}).
We also have $\Delta_{\omega_i, \delta} =
u' \Delta_{\omega_i, \omega_i}$ for some $u' \in U(\gg)$ which is a monomial in the divided
powers $f_j^{(n)}$.
Let us abbreviate $e^{(\pi)} = e_{i_1}^{(c_1 (\pi))} \cdots e_{i_m}^{(c_m (\pi))}$.
The equality $e^{(\pi)} \Delta_{\omega_i, \delta} = c \Delta_{\omega_i, \gamma}$
can be now rewritten as
$u e^{(\pi)} u' \Delta_{\omega_i, \omega_i} = c \Delta_{\omega_i, \omega_i}$.
This shows that part (i) of Theorem~\ref{th:minors-trails} is a consequence
of the following statement.

\begin{lemma}
\label{lem:path positivity}
If $u \in U(\gg)$ is a monomial of degree $0$ in the divided powers of the elements $e_j$ and $f_j$
then $u \Delta_{\omega_i, \omega_i} = c \Delta_{\omega_i, \omega_i}$
for some nonnegative integer $c$.
\end{lemma}

\proof
To see that $c \in \ZZ$ notice that the commutation relations in $U(\gg)$
between divided powers of the $f_j$ and $e_j$ involve
integer coefficients only.
It remains to show that $c \geq 0$.
First consider the case when $\gg$ is simply-laced, i.e., $|a_{ij}| \leq 1$ for $i \neq j$.
Then the nonnegativity of $c$ is a consequence of \cite[Theorem~4.3.13]{lu};
to be more precise, one applies the dual version of this result that says
that each generator $E_j$ or $F_j$ of $U_q(\gg)$ acts in the dual canonical basis in
$V_{\omega_i}$ by a matrix whose entries are nonnegative integer Laurent polynomials in $q$
(cf. \cite{bz93}).

If $\gg$ is not simply-laced, we use a well known embedding of
$\gg$ into a simply-laced complex semisimple Lie algebra $\tilde \gg$.
This embedding can be described as follows: if $\tilde \gg$
has Chevalley generators $f_{\tilde i}, \alpha_{\tilde i}^\vee$, and
$e_{\tilde i}$ for
$\tilde i \in I$ then the Chevalley generators of a subalgebra $\gg$ have the form
$$f_i = \sum_{\tilde i \in I_i} f_{\tilde i}, \,\,
\alpha_i^\vee = \sum_{\tilde i \in I_i} \alpha_{\tilde i}^\vee, \,\,
e_i = \sum_{\tilde i \in I_i} e_{\tilde i} \ ,$$
where the subsets $I_i \subset I$ are disjoint,
and no two indices from the same $I_i$ are adjacent to each other in the Dynkin
diagram of $\tilde \gg$.

The embedding $\gg \subset \tilde \gg$ allows us to identify any fundamental $\gg$-module
$V_{\omega_i}$ with the $\gg$-submodule generated by a highest vector of a fundamental
$\tilde \gg$-module $V_{\omega_{\tilde i}}$ for any $\tilde i \in I_i$.
Since any monomial in the $e_i$ and $f_i$ is a sum of monomials in the
$e_{\tilde i}$ and $f_{\tilde i}$, the desired inequality $c \geq 0$ follows from
the corresponding claim for $\tilde \gg$.
This completes the proofs of Lemma~\ref{lem:path positivity}
and Theorem~\ref{th:minors-trails}.
\endproof

\begin{remark}
\label{rem:step 1 minor}
{\rm In general, we do not know of a nice formula for the
coefficients of the monomials in Theorem~\ref{th:minors-trails}.
In some special cases these coefficients can be found with the
help of the following formulas which are easy consequences of the above proof:
$\Delta_{\gamma, \delta} (x x_i (t)) = \Delta_{\gamma, \delta}(x)$
for any $x \in G$ and $t \in \CC$ whenever $\delta (\alpha_i^\vee) \geq 0$;
and $\Delta_{\gamma, \delta} (x x_i (t)) = \Delta_{\gamma, \delta} (x)
+ t \Delta_{\gamma, s_i \delta} (x)$ whenever $\delta (\alpha_i^\vee) = -1$.}
\end{remark}

\section{Proofs of results in Sections~\ref{sec:multiplicities},
\ref{sec:new stuff on can}, \ref{sec:crystals} and \ref{sec:Plucker models}}
\label{sec:proofs main}


We start with the proof of Theorem~\ref{th:Lusztig-string}; the rest of the results
will follow quite easily.

\noindent {\bf Proof of Theorem~\ref{th:Lusztig-string}.}
In view of Theorem~\ref{th:twist special} and Corollary~\ref{cor:minors-trails-tropical},
part (i) of Theorem~\ref{th:Lusztig-string} is obtained via the tropicalization of (\ref{eq:t'-eta-t})
(with each fundamental weight $\omega_i$ replaced by $\omega_i^\vee$)

To prove part (ii), let us first rewrite conditions (1) -- (4) in Theorem~\ref{th:LR Lusztig-trails}
in terms of the element $x^\vee_\ii (t_1, \dots, t_m)$.

\begin{lemma}
\label{lem:Lusztig trails-Plucker}
Each of the conditions {\rm (1) -- (4)} in Theorem~\ref{th:LR Lusztig-trails}
is equivalent to the corresponding condition in Theorem~\ref{th:LR-Plucker-Lusztig}
with $M_{\gamma, \delta}= [\Delta_{\gamma, \delta}(x^\vee_\ii (t_1, \dots, t_m))]_{\rm trop}$.
\end{lemma}

\proof
For conditions (3) and (4), the claim follows from Corollary~\ref{cor:minors-trails-tropical}.
By the same theorem and the definition of $\ii$-trails,
we have
$$[\Delta_{\omega_i, s_i \omega_i}(x_\ii (t_1, \dots, t_m))]_{\rm trop}
= \min_{i_k = i} t_k \ ,$$
which proves our claim for the condition (1).
As for (2), our claim follows from (\ref{eq:highest minor}) (see also
Corollary~\ref{cor:unique trail from lambda} below).
\endproof

Our second step is to rewrite the same conditions in terms of the element
$\eta^{e,\wnot}(x^\vee_\ii (t_1, \dots, t_m))$.

\begin{lemma}
\label{lem:Plucker-to-Plucker}
The four conditions in Lemma~\ref{lem:Lusztig trails-Plucker}
are equivalent to conditions {\rm (1) -- (4)} in Theorem~\ref{th:LR-Plucker-strings}
with $M_{\gamma, \delta} = [\Delta_{\delta, \gamma}(\eta^{e,\wnot}(x^\vee_\ii (t_1, \dots, t_m)))]_{\rm trop}$.
\end{lemma}

\proof
Let us establish some identities
between generalized minors of an element $x' \in L^{e,\wnot}$
and those of $x = \eta^{e,\wnot} (x') \in L^{\wnot,e}$:
for any $i \in [1,r]$, we have

\begin{equation}
\label{eq:x-x' minors 1}
\Delta_{\omega_i, \wnot \omega_i}(x') =
\frac{1}{\Delta_{\omega_i, \omega_i}(x)} \ ;
\end{equation}

\begin{equation}
\label{eq:x-x' minors 2}
\Delta_{\omega_i, s_i \omega_i}(x') =
\Delta_{\wnot \omega_{i^*}, s_{i^*} \omega_{i^*}}(x) \ ;
\end{equation}

\begin{equation}
\label{eq:x-x' minors 3}
\Delta_{s_i \omega_i, \wnot \omega_i}(x') =
\frac{\Delta_{\wnot s_{i^*} \omega_{i^*}, \omega_{i^*}}(x)}
{\Delta_{\omega_i, \omega_i}(x)} \ ;
\end{equation}

\begin{equation}
\label{eq:x-x' minors 4}
\Delta_{\omega_i, \wnot s_i \omega_i}(x') =
\Delta_{s_i \omega_{i}, \omega_{i}}(x)
\prod_{j \neq i} \Delta_{\omega_j, \omega_j}(x)^{a_{ij}} \ .
\end{equation}

The identities (\ref{eq:x-x' minors 1}) are equivalent to
$[x]_0 = ([x' \overline {\wnot}]_0)^{-1}$, which is an
immediate consequence of (\ref{eq:eta2}).
As for (\ref{eq:x-x' minors 2}) -- (\ref{eq:x-x' minors 4}), they follow
by equating expressions (\ref{eq:negative t1 special}) and
(\ref{eq:positive t1 special}) with the corresponding expressions in
Corollary~\ref{cor:factors special}.
For example, (\ref{eq:x-x' minors 4}) is obtained by equating
the first expression in (\ref{eq:positive t1 special})
with the one in Corollary~\ref{cor:factors special} (ii) for $k = 1$
(note that in these formulas, one has to replace $x$ with $\tau_{\wnot}(x)$).

Our lemma is now obtained by a straightforward calculation (in which one applies
the above formulas to the group ${}^L G$).
For example, (\ref{eq:x-x' minors 1}) shows that conditions (2) in
Theorems~\ref{th:LR-Plucker-Lusztig} and \ref{th:LR-Plucker-strings}
are equivalent to each other.
\endproof

To complete the proof of Theorem~\ref{th:Lusztig-string} (ii),
it remains to show the following: if in Lemma~\ref{lem:Plucker-to-Plucker}
we replace $x^\vee_\ii (t_1, \dots, t_m)$ with $x^\vee_{\ii'} (t'_1, \dots, t'_m)$
and write the element $\eta^{e,\wnot}(x^\vee_{\ii'} (t'_1, \dots, t'_m)$ as
$x^\vee_{-\ii}(t_1, \dots, t_m)$ then the conditions in
Lemma~\ref{lem:Plucker-to-Plucker} become equivalent to the corresponding conditions in
Theorem~\ref{th:LR string-trails}.
This is done in a straightforward way by expanding each expression
$[\Delta_{\delta, \gamma}(x^\vee_{-\ii}(t_1, \dots, t_m))]_{\rm trop}$
with the help of Corollary~\ref{cor:minors-trails-tropical}.
\endproof

\noindent {\bf Proof of Theorem~\ref{th:li}.}
In view of (\ref{eq:li = ci}), the function
$l_i (b_\ii (t_1, \dots, t_m))$ can be computed by using a special case
of (\ref{eq:Lusztig-string}) with $k = 1$ (one also needs to
interchange $\ii$ with $\ii'$, and $t$ with $t'$ there).
It remains to notice that in this situation, the first minimum in
(\ref{eq:Lusztig-string}) becomes
$\sum_k s_{i_1} \cdots s_{i_{k-1}} \alpha_{i_k} (\omega_i^\vee) \ t_k$
by (\ref{eq:highest minor}) (see also
Corollary~\ref{cor:unique trail from lambda} below).
\endproof

\noindent {\bf Proof of Theorem~\ref{th:1st and last Lusztig-through-string}.}
Let us again consider two elements $x$ and $x'$ given by
$$x' = x_{\ii'} (t'_1, \dots, t'_m) \in L^{e,\wnot}_{> 0} \ ,$$
$$x = \eta^{e,\wnot} (x') = x_{-\ii}(t_1, \dots, t_m) \in L^{\wnot,e}_{> 0} \ .$$
Then each $t'_k$ is a subtraction-free rational expression in $t_1, \dots, t_m$,
and these expressions can be found with the help of Corollary~\ref{cor:factors special} (i).
In general, these expressions are a little cumbersome but for $k = 1$ or $k = m$
they can be simplified by using (\ref{eq:positive t1 special}) and
(\ref{eq:x-x' minors 1}) -- (\ref{eq:x-x' minors 4}); this was essentially done
in the above proof of Theorem~\ref{th:Lusztig-string}.
For example, here is the answer for $t'_m$:
\begin{equation}
\label{eq:last geom Lusztig through strings}
(t'_m)^{-1}  = \sum_{k: i_k^* = i'_m} t_k^{-1} \prod_{l > k} t_l^{- a_{i_l, i_k}}  \ .
\end{equation}
Formulas (\ref{eq:1st Lusztig}) and (\ref{eq:last Lusztig})
are obtained by ``tropicalizing" these expressions with the help of
Corollary~\ref{cor:minors-trails-tropical} (and passing from $G$
to ${}^L G$ as usual).
\endproof

\noindent {\bf Proof of Theorem~\ref{th:string cones general}.}
In the course of the proof of Theorem~\ref{th:Lusztig-string}, we have shown
that if $t = c_{\ii} (b_{\ii'} (t')^{\rm dual})$ then the nonnegativity of all
Lusztig parameters $t'_k$ (i.e., condition (1) in Theorem~\ref{th:LR Lusztig-trails})
is equivalent to the fact that the string components $t_k$ satisfy
inequalities (1) in Theorem~\ref{th:LR string-trails}.
This is precisely what we need to show.
\endproof

\begin{remark}
\label{rem:string cone}
{\rm 1. The above argument not only provides an explicit description of the string cones
but actually proves their existence thus providing an independent proof of
Proposition~\ref{pr:strings}.

\noindent 2. The above argument also implies that the string cones $C_\ii$
can be characterized in terms of the transition maps $R_{-\ii}^{-\ii'}$ as follows:
$C_\ii$ consists of all $t \in \RR^m$ such that, for any $\ii' \in R(\wnot)$,
the last component of $R_{-\ii}^{-\ii'}(t)$ is nonnegative.
This characterization was used in \cite{gp}, where the string cones were explicitly described
for type $A_r$.}
\end{remark}


\noindent {\bf Proof of Theorem~\ref{th:LR Lusztig-trails}.}
Taking into account Proposition~\ref{pr:Lusztig parametrizations}
and Corollary~\ref{cor:multiplicity through BB}, our statement
is an easy  consequence of Theorem~\ref{th:li}.
\endproof

\noindent {\bf Proof of Theorem~\ref{th:LR string-trails}.}
This is a consequence of Theorems~\ref{th:LR Lusztig-trails}
and \ref{th:Lusztig-string} (ii).
\endproof

\noindent {\bf Proof of Theorem~\ref{th:geom crystal}.}
Let $x  = x_{-\ii}(t_1, \dots, t_m) \in L^{\wnot,e}_{> 0}$, and
$\tilde x = (\eta^{e,\wnot} \circ F_i^{(c)} \circ
\eta^{\wnot,e})(x) = x_{- \ii}(\tilde t_1, \dots, \tilde t_m)$.
Let $x' = \eta^{\wnot, e} (x) \in L^{e, \wnot}_{> 0}$.
In view of (\ref{eq:geom crystal Lusztig}), we have
$x = \eta^{e, \wnot} (x')$ and
$\tilde x = \eta^{e, \wnot} (x' x_{i^*}(t))$,
where $t = (c-1)t'_m$.
By (\ref{eq:last geom Lusztig through strings}), we have
(in the notation of Theorem~\ref{th:geom crystal}):
\begin{equation}
\label{eq:additive shift}
t^{-1} = (c-1)^{-1} \sum_{l:i_l = i} T_l \ .
\end{equation}

Our goal is to express each $\tilde t_k$ in terms of $t_1, \dots, t_m$.
To do this, we combine (\ref{eq:t'-eta-t}) with its counterpart
for $\tilde x$:
$$\tilde t_k =  \frac{\Delta_{s_{i_1} \cdots s_{i_{k-1}} \omega_{i_k}, \wnot \omega_{i_k}}(x' x_{i^*}(t))}
{\Delta_{s_{i_1} \cdots s_{i_{k}} \omega_{i_k}, \wnot \omega_{i_k}}(x'x_{i^*}(t))}\ .$$
Since $(\wnot \omega_{i_k})(\alpha_{i^*}^\vee) = - \delta_{i_k, i}$,
it follows from Remark~\ref{rem:step 1 minor} that $\tilde t_k = t_k$ unless $i_k = i$;
furthermore, if $i_k = i$ then
$$\tilde t_k = \frac{\Delta_{s_{i_1} \cdots s_{i_{k-1}} \omega_{i}, \wnot
\omega_{i}}(x') + t \Delta_{s_{i_1} \cdots s_{i_{k-1}} \omega_{i}, \wnot s_i
\omega_{i}}(x')}
{\Delta_{s_{i_1} \cdots s_{i_{k}} \omega_{i}, \wnot
\omega_{i}}(x') + t \Delta_{s_{i_1} \cdots s_{i_{k}} \omega_{i}, \wnot s_i
\omega_{i}}(x')}$$
$$= t_k \cdot \frac{t^{-1} + \Delta_{s_{i_1} \cdots s_{i_{k-1}} \omega_{i}, \wnot s_i
\omega_{i}}(x')/\Delta_{s_{i_1} \cdots s_{i_{k-1}} \omega_{i}, \wnot
\omega_{i}}(x')}
{t^{-1} + \Delta_{s_{i_1} \cdots s_{i_{k}} \omega_{i}, \wnot s_i
\omega_{i}}(x')/\Delta_{s_{i_1} \cdots s_{i_{k}} \omega_{i}, \wnot
\omega_{i}}(x')} \ .$$

We claim that
\begin{equation}
\label{eq:two fractions}
\frac{\Delta_{s_{i_1} \cdots s_{i_{k-1}} \omega_{i}, \wnot s_i
\omega_{i}}(x')}{\Delta_{s_{i_1} \cdots s_{i_{k-1}} \omega_{i}, \wnot
\omega_{i}}(x')} = \sum \limits_{l\geq k:i_l=i} T_l, \quad
\frac{\Delta_{s_{i_1} \cdots s_{i_{k}} \omega_{i}, \wnot s_i
\omega_{i}}(x')}{\Delta_{s_{i_1} \cdots s_{i_{k}} \omega_{i}, \wnot
\omega_{i}}(x')} = \sum \limits_{l > k:i_l=i} T_l \ ;
\end{equation}
the desired equality (\ref{eq:geom crystal}) then follows by
plugging these expressions together with the one in (\ref{eq:additive shift})
into the above formula for $\tilde t_k$.

To prove (\ref{eq:two fractions}), we use the identity (\ref{eq:subflag minors})
with each variable $M_{\gamma, \delta}$ replaced by
$\Delta_{\gamma, \delta}(\tau_{\wnot} (x')) =
\Delta_{\wnot \delta, \wnot \gamma} (x')$ (see (\ref{eq:minors-tau})).
We thus obtain
$$\frac{\Delta_{s_{i_1} \cdots s_{i_{k-1}} \omega_{i}, \wnot s_i
\omega_{i}}(x')}{\Delta_{s_{i_1} \cdots s_{i_{k-1}} \omega_{i}, \wnot
\omega_{i}}(x')}$$
$$= \sum \limits_{l \geq k:i_l = i}
\frac{1}{\Delta_{s_{i_1} \cdots s_{i_{l-1}} \omega_{i}, \wnot
\omega_{i}}(x') \Delta_{s_{i_1} \cdots s_{i_{l}} \omega_{i}, \wnot
\omega_{i}}(x')} \prod_{j \neq i} \Delta_{s_{i_1} \cdots s_{i_{l}} \omega_{j}, \wnot
\omega_{j}}(x')^{- a_{ji}} \ .$$
As an easy consequence of (\ref{eq:t'-eta-t}), here the summand
corresponding to each index $l$ is equal to $T_l$; this proves the
first equality in (\ref{eq:two fractions}).
The second equality is proved in the same way.
This completes the proof of Theorem~\ref{th:geom crystal}.
\endproof

\noindent {\bf Proof of Theorem~\ref{th:Plucker model for tp cells}.}
The proof follows that of \cite[Theorem~2.7.1]{bfz}.
Let us sketch the proof of part (i); the proof of (ii) is the same.
Let $\MM^{\wnot} (K)^+$ denote the set of
all tuples $(M_{\omega_i,\gamma}) \in  \MM^{\wnot} (K)$ such that
$M_{\omega_i,\omega_i} = 1$ for all $i$.
It is clear that $p^+$ is a well-defined map
$\LL^{e, \wnot}(K) \to \MM^{\wnot} (K)^+$.
To show that this is a bijection, it suffices to construct the inverse map
in the case when $K = \RR_{> 0}$.
Since the map $\tt \mapsto x_{\ii}(t^\ii)$ is a bijection
between $\LL^{e,\wnot}(\RR_{> 0})$ and $L^{e,\wnot}_{> 0}$
(see Example~\ref{ex:real Luv}), we only need to construct a bijective
correspondence $(M_{\omega_i,\gamma}) \mapsto x$ between
$\MM^{\wnot} (\RR_{> 0})^+$ and $L^{e,\wnot}_{> 0}$.
This is done as follows: pick any $\ii \in R(\wnot)$; define the factorization parameters
$t_k$ as in Corollary~\ref{cor:factors special} (i) with each minor
$\Delta_{\omega_i,\gamma} (\psi^{\wnot,e}(x))$ replaced by $M_{\omega_i,\gamma}$;
form the corresponding product $x' = x_{- i_1} (t_1) \cdots x_{- i_m} (t_m)$;
and finally define $x = \psi^{\wnot,e}(x')$.
\endproof

\noindent {\bf Proof of Theorem~\ref{th:Plucker tropical}.}
The theorem follows by combining Theorem~\ref{th:Plucker model for tp cells}
and Corollary~\ref{cor:minors-trails-tropical} with Examples~\ref{ex:tropical Luv Lusztig}
and \ref{ex:tropical Luv strings} and with Lemmas~\ref{lem:Lusztig trails-Plucker}
and \ref{lem:Plucker-to-Plucker}.
\endproof

\noindent {\bf Proof of Theorems~\ref{th:LR-Plucker-Lusztig}
and \ref{th:LR-Plucker-strings}.}
These theorems follow by combining Theorems~\ref{th:LR Lusztig-trails} and
\ref{th:Plucker tropical} with Lemmas~\ref{lem:Lusztig trails-Plucker}
and \ref{lem:Plucker-to-Plucker}.
\endproof

\section{Proofs of results in Section~\ref{sec:string cones special}}
\label{sec:proofs strings special}

To deduce Theorems~\ref{th:string cone split}, \ref{th:fully commutative factors},
and \ref{th:minuscule} from Theorem~\ref{th:string cones general}, we need to develop
some properties of $\ii$-trails.
Although Theorem~\ref{th:string cones general} involves trails in the fundamental
modules over the Langlands dual Lie algebra ${}^L \gg$, we find it more convenient
to deal with trails in $\gg$-modules (translating from $\gg$ to ${}^L \gg$ is
automatic).

We start with some easy consequences of Definition~{\rm
\ref{def:trails}}.

\begin{proposition}
\label{trails general} {\rm (i)} Condition {\rm (2)} in Definition~{\rm
\ref{def:trails}} is equivalent to the following: $f_{i_l}^{c_l}
\cdots f_{i_1}^{c_1}$ is a non-zero linear map from $V(\gamma)$ to
$V(\delta)$.

\noindent {\rm (ii)} All the weights $\gamma_k$ in a $\ii$-trail are
weights of $V$, and $\gamma = \gamma_0 \geq \gamma_1 \geq \cdots
\geq \gamma_l = \delta$, where $\gamma \geq \delta$ has the usual
meaning that $\gamma - \delta$ is a nonnegative integer linear
combination of simple roots of $\gg$.
\end{proposition}

\proof Part (ii) is trivial; (i) follows from the well-known fact
that there exists a nondegenerate bilinear form $B$ on $V$ such
that $B(xv_1, v_2) = B(v_1, x^T v_2)$ for any $x \in \gg$ and
$v_1, v_2 \in V$.
\endproof

Throughout the rest of this section, we assume that $V$, $\gamma$,
$\delta$, and $\ii$ in Definition~\ref{def:trails} satisfy the
following conditions:

\smallskip

\noindent (1) $V = V_\lambda$, where $\lambda$ is a dominant weight for $\gg$;

\noindent (2) $\gamma$ and $\delta$ are two {\emph extremal} weights in $V_{\lambda}$,
i.e., they belong to the $W$-orbit $W \lambda$;

\noindent (3) $\ii = (i_1, \dots, i_l) \in R(w)$ for some $w \in W$.

\smallskip

Recall that the extremal weights in $V_{\lambda}$ are precisely
the vertices of the weight polytope $P(V_\lambda)$, and the
corresponding weight subspaces $V_{\lambda}(\gamma)$ are
one-dimensional. We call an $\ii$-trail from $\gamma$ to $\delta$
\emph{extremal} if all $\gamma_k$ are extremal weights of
$V_\lambda$. Extremal $\ii$-trails can be described as follows.
Let $W_\lambda$ denote the stabilizer of $\lambda$ in $W$. Every
$\gamma \in W \lambda$ has a unique \emph{minimal presentation}
$\gamma = u \lambda$, where $u \in W$ is the minimal (with respect
to the Bruhat order) representative of its coset $u W_\lambda$.

\begin{proposition}
\label{pr:extremal trails}
Suppose that extremal weights $\gamma$ and $\delta$ have minimal presentations
$\gamma = u \lambda$ and $\delta = v \lambda$.
There exists an extremal $\ii$-trail from $\gamma$ to $\delta$ in $V_\lambda$ if and only if
$\l(u v^{-1}) = \l(v) - \l(u)$, and $u v^{-1} \leq w$ (in the Bruhat order).
Under these conditions, the extremal $\ii$-trails from $\gamma$ to $\delta$
are in a bijection with subwords $(i_{k(1)}, \dots, i_{k(p)})$ of $\ii$ which are reduced
words for $u v^{-1}$; the $\ii$-trail corresponding to a sequence $(k(1) < \cdots < k(p))$
is given by $\gamma_k = s_{i_{k(j)}} \cdots s_{i_{k(1)}} \gamma$, where $j$ is the
maximal index such that $k(j) \leq k$.
\end{proposition}

\proof
First of all, the condition that an $\ii$-trail
$\pi = (\gamma = \gamma_0, \dots, \gamma_l = \delta)$
is extremal can be reformulated as follows:
$\gamma_k \in \{\gamma_{k-1}, s_{i_k} \gamma_{k-1}\}$
for $k = 1, \dots, l$.
Let $k(1)$ be the minimal index such that $\gamma_{k(1)} = s_{i_{k(1)}} \gamma \neq \gamma$.
Since $s_{i_{k(1)}} \gamma  = s_{i_{k(1)}} u \lambda < u \lambda$, it follows easily
that $\l(s_{i_{k(1)}} u) = \l (u) + 1$.
The standard properties of the Bruhat order in $W/ W_\lambda$ then imply that
$\gamma_{k(1)} = s_{i_{k(1)}} u \lambda$ is the minimal presentation of the weight
$\gamma_{k(1)}$.
Replacing $\gamma$ with $\gamma_{k(1)}$ and using induction on the length of $\ii$,
we obtain the desired statement.
\endproof

In particular, if $\lambda$ is \emph{minuscule}
(i.e., the extremal weights $w \lambda$
are the only weights of the $\gg$-module $V_{\lambda}$) then
\emph{all} $\ii$-trails from $\gamma$ to $\delta$ are extremal, and so are given
by Proposition~\ref{pr:extremal trails}.

\noindent {\bf Proof of Theorem~\ref{th:minuscule}.}
Remembering (\ref{eq:dk}) and using Proposition~\ref{pr:extremal trails}, we see
that the inequalities in (\ref{eq:extremal inequalities}) are precisely
the inequalities in Theorem~\ref{th:string cones general} corresponding to
extremal $\ii$-trails.
Thus all these inequalities hold on $C_\ii$.
If $\gg$ is of type $A_r$ then $\gg$ is isomorphic to ${}^L \gg$, and all its
fundamental weights are known to be minuscule.
So all the trails in Theorem~\ref{th:string cones general} are extremal,
and Theorem~\ref{th:minuscule} follows.
\endproof

Our next result gives upper and lower bounds for every $\ii$-trail (between two
extremal weights).
Let us fix $V$, $\gamma$, $\delta$, and $\ii$ as above.
Define weights
${\underline \gamma}_0, \dots, {\underline \gamma}_l$ and
${\overline \gamma}_0, \dots, {\overline \gamma}_l$ by setting
${\underline \gamma}_0 = \gamma$, ${\overline \gamma}_l = \delta$, and
$${\underline \gamma}_k = \min ({\underline \gamma}_{k-1}, s_{i_k} {\underline \gamma}_{k-1}),
\,\, {\overline \gamma}_{k-1} = \max ({\overline \gamma}_{k}, s_{i_k} {\overline \gamma}_{k})$$
for $k = 1, \dots, l$.

\begin{proposition}
\label{pr:bounds}
Any $\ii$-trail $(\gamma_0, \gamma_1, \cdots , \gamma_l)$ from $\gamma$ to $\delta$
in $V_\lambda$ satisfies ${\underline \gamma}_k  \leq \gamma_k \leq {\overline \gamma}_{k}$
for $k = 0, \dots, l$.
\end{proposition}

\proof Let us prove the inequalities $\gamma_k \leq {\overline
\gamma}_{k}$; the remaining ones are proved in a similar way with
the help of Proposition~\ref{trails general}(i). Let $v_\gamma$
denote a nonzero vector in $V_{\lambda}(\gamma)$. We proceed by
induction on $l$ to prove a little stronger statement:

\noindent (*) If a vector $v = e_{i_{k+1}}^{c_{k+1}} \cdots
e_{i_l}^{c_l} v_\delta$ is nonzero for some nonnegative integers
$c_{k+1}, \dots, c_l$ then its weight $\gamma_k = \gamma + c_{k+1}
\alpha_{i_{k+1}} + \cdots + c_l \alpha_{i_l}$ satisfies $\gamma_k
\leq {\overline \gamma}_{k}$.

We can assume that $l \geq 1$, and that our statement holds for any $\ii'$-trail,
where $\ii' = (i_1, \dots, i_{l-1})$.
Let us abbreviate $\overline \delta =  {\overline \gamma}_{l-1}$.
Consider two cases.

\noindent {\bf Case 1:} $\overline \delta = \delta \geq s_{i_l} \delta$.
In this case $\delta + \alpha_{i_l}$ is not a weight of $V_\lambda$.
Therefore, $\delta = \gamma_{l-1} = \overline \delta $, and our statement follows by
induction.

\smallskip

\noindent {\bf Case 2:} $\delta < s_{i_l} \delta = \overline \delta$.
Using the representation theory of $sl_2$, we see that $e_{i_l} v_{\overline \delta} = 0$,
and $v_\delta = f_{i_l}^a v_{\overline \delta}$, where
$a = {\overline \delta}(\alpha_{i_l}^\vee) > 0$.
Thus we have
$$v = e_{i_{k+1}}^{c_{k+1}} \cdots e_{i_l}^{c_l} v_\delta =
e_{i_{k+1}}^{c_{k+1}} \cdots e_{i_l}^{c_l} f_{i_l}^a v_{\overline \delta} \ .$$
Using the commutation relations between $f_{i_l}$ and the elements $e_{i}$ in $U(\gg)$,
we can express $v$ as a linear combination of vectors of the form
$f_{i_l}^{a'} e_{i_{k+1}}^{c'_{k+1}} \cdots e_{i_{l-1}}^{c'_{l-1}} v_{\overline \delta}$.
Hence at least one of these vectors is nonzero, and we conclude by induction that
$$\gamma_k \leq \gamma_k + a' \alpha_{i_l} \leq {\overline \gamma}_{k} \ ,$$
as required.
\endproof

\begin{corollary}
\label{cor:unique trail}
Suppose that
$\gamma \geq s_{i_1} \gamma \geq \cdots \geq s_{i_l} \cdots s_{i_1} \gamma = \delta$.
Then there is a unique $\ii$-trail $\pi$ from $\gamma$ to $\delta$ in $V_\lambda$;
it is given by $\gamma_k = s_{i_k} \cdots s_{i_1} \gamma$,
and it has $c_k (\pi) = \gamma (s_{i_1} \cdots s_{i_{k-1}} \alpha_{i_k}^\vee)$,
and $d_k (\pi) = 0$ for $k = 1, \dots, l$.
\end{corollary}

\proof
The uniqueness of $\pi$ follows from Proposition \ref{pr:bounds} since,
under the present assumptions,
we have ${\underline \gamma}_k = {\overline \gamma}_{k}= s_{i_k} \cdots s_{i_1} \gamma$
for $k = 0, \dots, l$.
The claim about $c_k (\pi)$ and $d_k (\pi)$ follows at once from the definitions
(\ref{eq:ck}) and (\ref{eq:dk}).
\endproof

The following special case of Corollary~\ref{cor:unique trail} extends
\cite[Proposition~3.3]{ber}.

\begin{corollary}
\label{cor:unique trail from lambda}
Suppose that $\gamma = u \lambda$ for some $u \in W$ such that $\l(w^{-1} u) = \l(u) + \l(w)$.
For every $\ii \in R(w)$, there is a unique $\ii$-trail $\pi$ from $\gamma$ to $w^{-1}\gamma$;
it is given by $\gamma_k = s_{i_k} \cdots s_{i_1} \gamma$ for $k = 0, \dots, l$.
\end{corollary}

Our proof of Theorem~\ref{th:string cone split} relies on one more corollary of
Proposition~\ref{pr:bounds}; to formulate it we need the following definition.

\begin{definition}
\label{def:splitting}
{\rm Let $V$, $\gamma$, $\delta$, and $\ii$ have the same meaning as in conditions (1)--(3) above.
An index $k \in [0,l]$ is splitting if it satisfies the following conditions:

\noindent (1) $\gamma \geq s_{i_1} \gamma \geq \cdots \geq s_{i_k} \cdots s_{i_1} \gamma$;

\noindent (2) $\delta \leq s_{i_l} \delta \leq \cdots \leq s_{i_{k+1}} \cdots s_{i_l} \delta$;

\noindent (3) $s_{i_{k+1}} \cdots s_{i_l} \delta - s_{i_k} \cdots s_{i_1} \gamma$
is a simple root for $\gg$.}
\end{definition}

\begin{corollary}
\label{cor:splitting}
If an index $k \in [0,l]$ is splitting then every $\ii$-trail $(\gamma_0, \dots, \gamma_l)$
from $\gamma$ to $\delta$ has either $\gamma_j = s_{i_j} \cdots s_{i_1} \gamma$ for
$0 \leq j \leq k$, or $\gamma_j = s_{i_{j+1}} \cdots s_{i_l} \delta$ for
$k \leq j \leq l$.
\end{corollary}

\proof
Conditions (1) and (2) in Definition~\ref{def:splitting} imply that
${\underline \gamma}_j = s_{i_j} \cdots s_{i_1} \gamma$ for $0 \leq j \leq k$, and
${\overline \gamma}_{j} = s_{i_{j+1}} \cdots s_{i_l} \delta$ for
$k \leq j \leq l$.
Combining condition (3) with Proposition~\ref{pr:bounds}, we conclude that
$\gamma_k$ must be equal to $s_{i_k} \cdots s_{i_1} \gamma$
or $s_{i_{k+1}} \cdots s_{i_l} \delta$.
In the former case, Corollary~\ref{cor:unique trail} guarantees the uniqueness
of an $(i_1, \dots, i_k)$-trail from $\gamma$ to $\gamma_k$, and we conclude that
$\gamma_j = s_{i_j} \cdots s_{i_1} \gamma$ for
$0 \leq j \leq k$.
Similarly, in the latter case, we have $\gamma_j = s_{i_{j+1}} \cdots s_{i_l} \delta$ for
$k \leq j \leq l$, as required.
\endproof

Let $C_\ii (\gamma, \delta)$ denote the cone
of tuples $(t_1, \dots, t_l) \in \RR^l$
such that $\sum_k d_k (\pi) t_k \geq 0$ for any
$\ii$-trail $\pi$ from $\gamma$ to $\delta$.
Our next result is an immediate consequence of Corollary~\ref{cor:splitting}.

\begin{lemma}
\label{lem:cone separation}
Suppose an index $k$ is splitting (for $\ii$, $\gamma$, and $\delta$), and let
$\ii^{(1)} = (i_1, \dots, i_k)$, and $\ii^{(2)} = (i_{k+1}, \dots, i_l)$.
Then we have
$$C_\ii (\gamma, \delta) = C_{\ii^{(1)}} (\gamma, s_{i_{k+1}} \cdots s_{i_l} \delta)
\times C_{\ii^{(2)}} (s_{i_k} \cdots s_{i_1} \gamma, \gamma) \ .$$
\end{lemma}

\noindent {\bf Proof of Theorem~\ref{th:string cone split}.}
It is enough to prove (\ref{eq:string cone split}) for the case when $p= 2$,
i.e., when the flag $\emptyset = I_0 \subset I_1 \subset \cdots \subset I_p = [1,r]$
has only one proper subset $I_1$.
Thus $\ii$ is the concatenation $(\ii^{(1)}, \ii^{(2})$,
where $\ii^{(1)} \in R(\wnot (I_1))$, and
$\ii^{(2)} \in R(\wnot (I_1)^{-1} \wnot)$;
let $l = \l (\wnot)$ be the length of $\ii$, and $k = \l (\wnot (I_1))$ be the length
of $\ii^{(1)}$.

In our present notation, each cone $C_\ii (u,v)$ is the intersection of the cones
$C_{\ii} (u\omega_i^\vee, v s_i \omega_i^\vee)$ for all $i \in [1,r]$; in particular,
the string cone $C_\ii$ in (\ref{eq:string cone split}) is equal to
$$C_\ii = \bigcap_{i \in [1,r]} C_{\ii} (\omega_i^\vee, \wnot s_i \omega_i^\vee) \ .$$
Therefore, it suffices to show that
\begin{equation}
\label{eq:splitting parabolic}
C_{\ii} (\omega_i^\vee, \wnot s_i \omega_i^\vee) =
C_{\ii^{(1)}} (\omega_i^\vee, \wnot (I_1) s_i \omega_i^\vee) \times
C_{\ii^{(2)}} (\wnot (I_1) \omega_i^\vee, \wnot s_i \omega_i^\vee)
\end{equation}
for every $i \in [1,r]$.
Let us distinguish two cases.

\noindent {\bf Case 1:} $i \in I_1$. We claim that in this case the index $k$ is splitting
for $\ii$, $\gamma = \omega_i^\vee$, and $\delta = \wnot s_i \omega_i^\vee$.
Condition (1) in Definition~\ref{def:splitting} is obvious.
To prove (2), we need to show that
$s_{i_{j+1}} \cdots s_{i_l} \wnot s_i \omega_i^\vee (\alpha_{i_{j+1}}) \geq 0$ for
$k \leq j < l$.
This follows from the fact that $s_i \wnot(I_1) s_{i_{k+1}} \cdots s_{i_j} (\alpha_{i_{j+1}})$
is a positive root for $\gg$ (this is clear since $(i_{k+1}, \dots,  i_{j+1})$
is a reduced word for a left factor of $\wnot (I_1)^{-1} \wnot$).
Finally, to prove (3) notice that in the present situation we have
$$s_{i_{k+1}} \cdots s_{i_l} \delta - s_{i_k} \cdots s_{i_1} \gamma =
\wnot (I_1) s_i \omega_i^\vee - \wnot (I_1) \omega_i^\vee = - \wnot (I_1) \alpha_i^\vee \ ,$$
which is a simple root for ${}^L \gg$ whenever $i \in I_1$.

The desired equality (\ref{eq:splitting parabolic}) now follows from
Lemma~\ref{lem:cone separation}.

\noindent {\bf Case 2:} $i \notin I_1$.
Since in this case $s_j \omega_i^\vee = \omega_i^\vee$ for any $j \in I_1$,
it follows that every $\ii^{(1)}$-trail from $\omega_i^\vee$ is trivial, i.e., all its components
are equal to $\omega_i^\vee$.
Thus in this case both sides of (\ref{eq:splitting parabolic}) are equal to
$\RR^k \times C_{\ii^{(2)}} (\omega_i^\vee, \wnot s_i \omega_i^\vee)$.
This concludes the proof of Theorem~\ref{th:string cone split}.
\endproof

\noindent {\bf Proof of Theorem~\ref{th:fully commutative factors}.}
Without loss of generality, we can assume that the sets $I_{j-1}$
and $I_j$ in the formulation are equal to $[1,r-1]$ and $[1,r]$
respectively.
Let us abbreviate $\wnot' = \wnot([1,r-1])$, and let $w = \wnot'
\wnot$, and $\ii = (i_1, \dots, i_l) \in R(w)$.
We need to show that if $w$ is fully commutative then the cone
$$\bigcap_{i \in [1,r]} C_\ii (\wnot' \omega_i^\vee, \wnot s_i
\omega_i^\vee) \subset \RR^l$$
is given by the inequalities in Theorem~\ref{th:fully commutative factors}.

First, let us show that $C_\ii (\wnot' \omega_r^\vee, \wnot
s_r \omega_r^\vee)$ is given by the only inequality $t_l \geq 0$
(this part does not use the fact that $w$ is fully commutative).
Notice that $\wnot' \omega_r^\vee =  \omega_r^\vee$.
Notice also that the coroot
$w \alpha_i^\vee =  \wnot' \wnot \alpha_i^\vee = - \wnot' \alpha_{i^*}^\vee$
is negative for $i = r^*$, and positive otherwise.
It follows that $i_l = r^*$ regardless of the choice of $\ii \in R(w)$.
Since $\wnot s_r \omega_r^\vee (\alpha_{r^*}) = - s_r
\omega_r^\vee(\alpha_r) = -1$, it follows that any $\ii$-trail
$\pi = (\gamma_0, \dots, \gamma_l)$ from $\omega_r^\vee$ to
$\wnot s_r \omega_r^\vee$ must have
$\gamma_{l-1} = \gamma_l = \wnot s_r \omega_r^\vee$.
Using Corollary~\ref{cor:unique trail from lambda}, we see that
$\pi$ is unique, and we have
$\gamma_k = s_{i_k} \cdots s_{i_1} \gamma$ for $k = 0, \dots, l-1$.
It follows that $d_k (\pi) = 0$ for $k = 0, \dots, l-1$, and $d_l (\pi) = 1$.
The corresponding linear inequality defining the cone
$C_\ii (\omega_r^\vee, \wnot s_r \omega_r^\vee)$ is $t_l \geq 0$
as claimed.

Now let us show the following:

\noindent (*) Each of the cones
$C_\ii (\wnot' \omega_{i^*}^\vee, \wnot s_{i^*} \omega_{i^*}^\vee)$
for $i \neq r^*$ is given by some of the inequalities (2)--(4) in
Theorem~\ref{th:fully commutative factors}.

This can be done by analyzing the corresponding
$\ii$-trails but we prefer another method using geometric lifting.
By Corollary~\ref{cor:minors-trails-tropical}, the cone
$C_\ii (\wnot' \omega_{i^*}^\vee, \wnot s_{i^*} \omega_{i^*}^\vee)$
consists of all integer $l$-tuples $(t_1, \dots, t_l)$ satisfying
the inequality
$$[\Delta_{w \omega_{i}^\vee, s_{i}  \omega_{i}^\vee}
(x_{-\ii} (t_1, \dots, t_l))]_{\rm trop} \geq 0 \ .$$
We shall deal with the minor $\Delta_{w \omega_{i}, s_{i}\omega_{i}}$
instead of $\Delta_{w \omega_{i}^\vee, s_{i}  \omega_{i}^\vee}$
(so in the resulting formulas one will have to replace the Cartan matrix with its
transpose).

A calculation in $SL_2$ shows that
$\overline s_i = \lim_{t \to \infty} x_{-i}(t) x_{i}(-t)$.
It follows that the minor in question can be written as
$$\Delta_{w \omega_{i}, s_{i}\omega_{i}} (x_{-\ii} (t_1, \dots, t_l))=
\lim_{t \to \infty} \Delta_{w \omega_{i}, \omega_{i}}
(x_{-\tilde \ii} (t_1, \dots, t_l, t)) \ ,$$
where $\tilde \ii$ is the word $(i_1, \dots, i_l, i)$.
Since $i \neq r^*$, we have $\tilde \ii \in R(w s_i)$.
We also have $w s_i w^{-1} =  \wnot' s_{i^*} {\wnot'}^{\ -1} =
s_{i'}$ for some $i' \in [1,r-1]$.
Therefore, the word $\tilde \ii' = (i', i_1, \dots, i_l)$
is also a reduced word for $w s_i$.
With the help of the transition maps in Proposition~\ref{pr:negative
transition}, we can express the product $x_{-\tilde \ii} (t_1, \dots, t_l, t)$
as $x_{-\tilde \ii'} (p, p_1, \dots, p_l)$, where $p$ and all $p_k$
are subtraction-free rational expressions in $t_1, \dots, t_l, t$.
Since $w^{-1} (\alpha_{i'}) = \alpha_i$, it easily follows from
Propositions~\ref{pr:Luv equations} and \ref{pr:Luv coordinates} that
$$\Delta_{w \omega_{i}, \omega_{i}}(x_{-\tilde \ii} (t_1, \dots, t_l,t))=
\Delta_{w \omega_{i}, \omega_{i}}(x_{-\tilde \ii'} (p, p_1, \dots, p_l))= $$
$$p^{-1} \Delta_{w \omega_{i}, \omega_{i}}(x_{-\ii} (p_1, \dots, p_l)) = p^{-1} \ .$$
Thus it remains to compute $p$ as a subtraction-free rational expression
$p(t_1, \dots, t_l, t)$ and take its limit as $t \to \infty$.

To compute $p$, we shall use a combinatorial lemma valid for
arbitrary Coxeter groups.
It uses the following notation: for any two distinct indices $i$ and $j$ from $[1,r]$,
let $\wnot(i,j)$ denote the longest element
in the parabolic subgroup of $W$ generated by $s_i$ and $s_j$.
Thus $\wnot(i,j) = s_i s_j s_i \cdots = s_j s_i s_j \cdots$ (both products have $d$ factors,
where $d$ is the order of $s_i s_j$ in $W$).

\begin{lemma}
\label{lem:simple-to-simple}
Suppose an element $w$ of an arbitrary Coxeter group $W$, and an
index $i$ are such that $\l (w s_i) = \l(w) + 1 > 1$, and
$w s_i w^{-1} = s_{i'}$ for some $i'$.
Then there exist an index $j \neq i$ and a reduced word $\ii'$ of $w$ such
that $s_i s_j$ has finite order $d$, and $\ii'$ ends with the word
$(\dots, j,i,j) \in R(\wnot(i,j)s_i)$ of length $d-1$.
\end{lemma}

\proof
Pick any reduced word $(i_1, \dots, i_l) \in R(w)$; clearly,
$i_l \neq i$.
Both words $(i_1, \dots, i_l, i)$ and $(i', i_1, \dots, i_l)$ are
reduced words for $w s_i$.
By the Tits theorem, the latter word can be obtained from the
former one by a sequence of $d$-moves.
Consider the first move in this sequence that involves the last letter.
If this move is performed on a word $(\ii', i)$ then the word
$\ii' \in R(w)$ has the desired property.
\endproof

Applying Lemma~\ref{lem:simple-to-simple} several times if
necessary, we obtain a reduced word $\ii' \in R(w)$ with the
following property: $\ii'$ is a concatenation
$\ii^{(1)}, \dots, \ii^{(n)}$ such that each $\ii^{(k)}$
consists of two alternating letters, and one obtains a reduced word
$(\ii', i)$ from $(i', \ii')$ by a sequence of $n$ $d$-moves,
the $k$-th one involving $\ii^{(k)}$ and the preceding index.
Thus if we replace $\ii$ with $\ii'$ then the above rational
function $p(t_1, \dots, t_l, t)$ is easily computed by
repeatedly using the formula for $1/p_1$ in Proposition~\ref{pr:negative
transition}.
Taking its tropicalization, we conclude that the desired property
(*) holds for $\ii'$; more precisely, the
cone $C_{\ii'} (\wnot' \omega_{i^*}^\vee, \wnot s_{i^*} \omega_{i^*}^\vee)$
is given by the inequalities (2)--(4) in Theorem~\ref{th:fully commutative factors}
corresponding to all the intervals $\ii^{(k)}$ that have length $> 1$.

On the other hand, the fact that $\ii$ satisfies (*)
is clearly preserved by switches $(i_k,i_{k+1}) \to (i_{k+1},i_k)$ with
$a_{i_k,i_{k+1}} = 0$.
Therefore if $w$ is fully commutative then (*) holds for any $\ii \in R(w)$.

To complete the proof of Theorem~\ref{th:fully commutative factors},
it remains to show that conversely each of the
inequalities (2)--(4) appears as one of the defining inequalities for
$C_{\ii'} (\wnot' \omega_{i^*}^\vee, \wnot s_{i^*} \omega_{i^*}^\vee)$
with some $i \neq r^*$.
Without loss of generality, we can assume that our inequality
corresponds to a subword
$(i_{k+1}, \dots, i_{k+d-1}) = (j,i,j, \dots)$ of $\ii$, where
$d \in \{3,4,6\}$ is the order of $s_i s_j$.
Let $u = s_{i_1} \cdots s_{i_k}$.
In view of the above argument, to complete the proof
it suffices to show that $\l (u s_i) = \l(u) + 1$, and
$u s_i u^{-1} = s_{i'}$ for some  $i' \in [1,r-1]$.
In other words, it suffices to show that the root $\beta = u \alpha_{i}$ is one of
the simple roots $\alpha_{1}, \dots, \alpha_{r-1}$.

Notice that $i_k \neq i$ (otherwise $\ii$ would contain a
subword $(i_k, \dots, i_{k+d-1}) \in R(\wnot (i,j))$
which is impossible since $w$ is fully commutative).
Therefore, the root $\beta$ is positive.
On the other hand, the word $(i_{k+1}, \dots, i_l)$ is a reduced
word for $u^{-1} w$.
Since $w$ (and hence $u^{-1} w$) is fully commutative, and
$a_{i, i_{k+1}} = a_{ij} \neq 0$,
it follows that no reduced word for $u^{-1} w$ can begin with $i$.
Therefore, the root $w^{-1} u \alpha_{i} = w^{-1} \beta$ is also positive.
Since $w^{-1} \beta = \wnot \wnot' \beta$,
we conclude that the root $\wnot' \beta$ is negative.
Since $\beta$ is a positive root sent to a negative
one by $\wnot'$, it follows that $\beta$ does not contain
$\alpha_r$, i.e., it is a positive integer linear combination
of $\alpha_{1}, \dots, \alpha_{r-1}$.

As a final step in our argument, notice that any reduced word for
$u$ begins with $i_1 = r$ (since $\l (\wnot' s_{i'}) = \l (\wnot') -1$
for any $i' \in [1,r-1]$).
It follows that $u^{-1} \alpha_{i'}$ is a positive root
for any $i' \in [1,r-1]$.
Therefore, $u^{-1} \beta$ can be a simple root only if
$\beta$ is simple.
This completes the proof of Theorem~\ref{th:fully commutative
factors}.
\endproof

\section{Proofs of results in Sections~\ref{sec:LR BCD} and \ref{sec:reduction}}
\label{sec:special multiplicities}

\noindent {\bf Proofs of Theorems~\ref{th:LR string-trails B} and \ref{th:LR string-trails D}.}
To prove Theorem~\ref{th:LR string-trails B}, we apply Theorem~\ref{th:LR string-trails}
to the following reduced word:
$$\ii = (\ii^{(0)}, \dots, \ii^{(r-1)}) \in
R(\wnot), \quad \ii^{(j)} = (j, \dots, 1, 0, 1, \dots, j) \ .$$
We shall rename the variables $t_k$ as follows: the variables
corresponding to each interval $\ii^{(j)}$ will be denoted
$(t_{-j}^{(j)}, \dots, t_{-1}^{(j)}, t_{0}^{(j)}, t_{1}^{(j)},
\dots, t_{j}^{(j)})$.
We only need to show that for this choice of $\ii$, each of the
conditions (1)--(4) in Theorem~\ref{th:LR string-trails}
specializes to the corresponding condition in Theorem~\ref{th:LR string-trails B}.
For conditions (2) and (4) this is straightforward, and for (1)
this is a special case of Theorem~\ref{th:fully commutative factors}
(the corresponding string cone was already found in \cite{lit}).
It remains to analyze condition (3).

Let $\pi = (\gamma_i^{(j)})$ be an $\ii$-trail from $s_j \omega_j^\vee$ to
$\wnot \omega_j^\vee$ in an ${}^L \gg$-module $V_{\omega_j^\vee}$;
here the components $\gamma_k$ of $\pi$ are renamed in the same way as the
corresponding variables $t_k$.
We use Proposition~\ref{pr:bounds} to obtain upper and lower
bounds for the weights $\gamma_i^{(j)}$.
First, since $s_i s_j \omega_j^\vee = s_j \omega_j^\vee$ for $i < j-1$, we have
$\underline \gamma_{i}^{(j')} = s_j \omega_j^\vee$ and $d_{i}^{(j')}(\pi) = 0$
for all $|i| \leq j' < j-1$.
Second, an easy calculation shows that
$\underline \gamma_{j}^{(j)}= \overline \gamma_{j}^{(j)} = s_j
\cdots s_1 s_0 s_1 \cdots s_j \omega_j^\vee$.
It follows that every component $\gamma_{i}^{(j')}$ of $\pi$ such
that $j' > j$ is obtained from the previous component by the
action of $s_{|i|}$; therefore, $d_{i}^{(j')}(\pi) = 0$ whenever $j' > j$.
Thus, the only part of $\pi$ that can contribute to an inequality
in (3) is a $(\ii^{(j-1)}, \ii^{(j)})$-trail from $s_j \omega_j^\vee$ to
$s_j \cdots s_1 s_0 s_1 \cdots s_j \omega_j^\vee$ in $V_{\omega_j^\vee}$.

If $j = 0$ then $(\ii^{(j-1)}, \ii^{(j)}) = (0)$, so there is a
unique trail $\pi$, and the only non-zero number $d_k (\pi)$ is
$d_0^{(0)} (\pi) = -1$.
The corresponding inequality in (3) is $\lambda (\alpha_{0}^\vee)
\geq t_0^{(0)}$, as claimed.

Now let $j > 0$.
Then we can assume without loss of generality that $j = r-1$.
Let us also assume that $\gg = so_{2r+1}$ is of type $B_r$.
Then $V_{\omega_{j}^\vee}$ is the $2r$-dimensional standard
(sometimes also called \emph{vector}) representation of ${}^L \gg = sp_{2r}$.
This ${}^L \gg$-module is minuscule, so the trails in question are extremal.
Therefore, they are given by Proposition~\ref{pr:extremal trails}
with $\lambda = \omega_{j}^\vee$, $\ii = (\ii^{(j-1)}, \ii^{(j)})$,
$u = s_{j}$, and $v = s_{j} \cdots s_1 s_0 s_1 \cdots s_{j}$.
We see that these trails correspond to all occurrences of a
word $(j-1, \dots, 1, 0, 1, \dots, j)$ as a subword of
$\ii = (j-1, \dots, 1, 0, 1, \dots, j-1; j,  \dots, 1, 0, 1, \dots, j)$.
By inspection, there are $4j-1$ such subwords falling into the
following $4$ classes (in each case, we represent a subword by the
list of variables $t_{i'}^{j'}$ corresponding to positions \emph{not} belonging
to this subword):

\smallskip

\noindent (1) $\{t_{i'}^{(j-1)} \, (i' > i); t_{i'}^{(j)} \, (i' \leq
i)\}$ for $0 \leq i < j$;

\smallskip

\noindent (2) $\{t_{i'}^{(j-1)} \, (i' > i); t_{i'}^{(j)} \, (-i - 1 \neq i' \leq
i+1)\}$ for $0 \leq i < j$;

\smallskip

\noindent (3) $\{t_{i'}^{(j-1)} \, (i' \geq  -i); t_{i'}^{(j)} \, (i'
< -i)\}$ for $0 \leq i < j$;

\smallskip

\noindent (4) $\{t_{i'}^{(j-1)} \, (i+1 \neq i' \geq -i-1); t_{i'}^{(j)} \,
(-i - 1 \neq i' \leq i+1)\}$ for $0 \leq i < j-1$.

\smallskip

Computing the coefficients $d_{i'}^{(j')} (\pi)$ for the $4j-1$
extremal $\ii$-trails $\pi$ corresponding to these subwords, we
see that the inequalities in Theorem~\ref{th:LR string-trails B}(3)
are indeed specializations of those in Theorem~\ref{th:LR string-trails}(3).
This completes the proof of Theorem~\ref{th:LR string-trails B}
for $\gg$ of type $B_r$.

Now suppose that $\gg = sp_{2r}$ is of type $C_r$.
The extremal trails and corresponding linear inequalities are
described in the same way as for the type $B_r$.
The only difference is that the ${}^L \gg$-module
$V_{\omega_{j}^\vee}$ for $j = r-1$ is now the ($(2r+1)$-dimensional) vector
representation of $so_{2r+1}$, and it is no longer minuscule.
But it is \emph{quasi-minuscule}, that is, its only
non-extremal weight is the zero weight.
It follows easily that there is a unique non-extremal
$(\ii^{(j-1)}, \ii^{(j)})$-trail $\pi$ from $s_j \omega_j^\vee$ to
$s_j \cdots s_1 s_0 s_1 \cdots s_j \omega_j^\vee$ in
$V_{\omega_j^\vee}$;
and the only non-zero coefficients $d_{i'}^{(j')} (\pi)$ are
$d_{0}^{(j-1)} (\pi) = 1$ and $d_{0}^{(j)} (\pi) = -1$.
The resulting linear inequality in (3) is
$\lambda (\alpha_{j}^\vee) \geq t_0^{(j)} - t_{0}^{(j-1)}$.
But we do not have to include this inequality since it is a
consequence of the two inequalities $\lambda (\alpha_{j}^\vee) \geq
2 t_0^{(j)} - t_1^{(j-1)} - t_{-1}^{(j)}$ and
$\lambda (\alpha_{j}^\vee) \geq t_1^{(j-1)} + t_{-1}^{(j)} - 2t_0^{(j-1)}$
corresponding to extremal trails.
This completes the proof of Theorem~\ref{th:LR string-trails B}.
\endproof

The proof of Theorem~\ref{th:LR string-trails D} is very similar,
and we leave the details to the reader.
Note only that when $\gg = so_{2r}$ is of type $D_r$, we apply
Theorem~\ref{th:LR string-trails} to the following reduced word:
$$\ii = (\ii^{(1)}, \dots, \ii^{(r-1)}) \in
R(\wnot), \quad \ii^{(j)} = (j, \dots, 2, -1, 1, 2 \dots, j) \ .$$
The variables $t_k$ corresponding to each interval $\ii^{(j)}$ are now denoted
$$(t_{-j}^{(j)}, \dots, t_{-1}^{(j)}, t_{1}^{(j)}, \dots, t_{j}^{(j)})\ .$$
The same argument as for the type $B_r$ above shows that all the
$\ii$-trails in Theorem~\ref{th:LR string-trails}(3) are extremal,
and the corresponding linear inequalities are precisely those in
Theorem~\ref{th:LR string-trails D}(3).  \endproof

\noindent {\bf Proofs of Theorems~\ref{th:reduction multiplicity Lusztig} and
\ref{th:reduction multiplicity string}.}
Let us recall a well-known relationship between the reduction
multiplicities and tensor product multiplicities (see e.g.,
\cite{bz88}): the multiplicity $c_{nu, \beta}$ of $V_\beta^{(I)}$ in the reduction of
$V_\nu$ to $\gg (I)$ is equal to $c_{\lambda, \nu}^{\lambda + \beta}$
for any weight $\lambda$ such that $\lambda (\alpha_i^\vee) = 0$
for $i \in I$, and $\lambda (\alpha_i^\vee) \gg 0$ for $i \in [1,r] \setminus I$.
(This follows from the interpretation of $c_{nu, \beta}$ as the dimension
of the subspace of vectors $v \in V_\nu$ of weight $\beta$ such
that $e_i v = 0$ for $i \in I$.)
Thus, an expression for $c_{nu, \beta}$ can be obtained by
computing $c_{\lambda, \nu}^{\lambda + \beta}$ using
either of the theorems \ref{th:LR Lusztig-trails} and
\ref{th:LR string-trails}.
In doing so we choose a reduced word $\tilde \ii \in R(\wnot)$ as
a concatenation $(\ii', \ii)$, where
$\ii' = (i'_1, \dots, i'_{m'})  \in R(\wnot (I))$ and
$\ii = (i_1, \dots, i_n)  \in R(\wnot (I)^{-1} \wnot)$.
Let us write the corresponding variables $t_k$ that appear in
Theorems~\ref{th:LR Lusztig-trails} and \ref{th:LR string-trails}
as $\tilde t = (t'_1, \dots, t'_{m'}, t_1, \dots, t_n)$.

We now claim that conditions (1)--(3) in each of
Theorems~\ref{th:LR Lusztig-trails} and \ref{th:LR string-trails}
imply that $t'_1 = \cdots = t'_{m'} = 0$.
It should be possible to deduce this directly but we prefer another argument.
Let us deal with Theorem~\ref{th:LR Lusztig-trails};
Theorem~\ref{th:LR string-trails} can be treated in the same way.
Recall that condition (3) (combined with (1) and (2)) in
Theorem~\ref{th:LR Lusztig-trails} was obtained
as a reformulation of the following
(see Proposition~\ref{pr:Lusztig parametrizations} (ii) and
Corollary~\ref{cor:multiplicity through BB}):
$l_i (b_{\tilde \ii} (\tilde t)) \leq \lambda (\alpha_i^\vee)$
for all $i \in [1,r]$.
With the choice of $\lambda$ as above, this just means that
$l_i (b_{\tilde \ii} (\tilde t))= 0$ for $i \in I$.
However if $i \in I$ then
$l_i (b_{\tilde \ii} (\tilde t)) = l_i (b_{\ii'}(t'_1, \dots, t'_{m'}))$,
in view of Proposition~\ref{pr:Lusztig parametrizations}(ii);
here $b_{\ii'}(t'_1, \dots, t'_{m'})$ is understood as a canonical
basis vector for the (semisimple part of) the Lie algebra $\gg (I)$
instead of $\gg$.
Remembering the definition of $l_i (b)$, we conclude that
$b_{\ii'}(t'_1, \dots, t'_{m'}) = 1$ hence $t'_1 = \cdots = t'_{m'} =0$,
as claimed.

Now the conditions on $t_1, \dots, t_m$ in
Theorem~\ref{th:reduction multiplicity Lusztig} (resp.
Theorem~\ref{th:reduction multiplicity string}) are easily  seen to be
equivalent to conditions (1), (2) and (4) in
Theorem~\ref{th:LR Lusztig-trails} (resp. Theorem~\ref{th:LR string-trails});
for conditions Theorem~\ref{th:reduction multiplicity string} (1)
and Theorem~\ref{th:reduction multiplicity Lusztig} (3),
this follows from the splitting property (\ref{eq:string cone split})
in Theorem~\ref{th:string cone split}.
\endproof

\noindent {\bf Proofs of Corollaries~\ref{cor:pq reduction Lusztig} and
\ref{cor:pq reduction string}.}
These corollaries are obtained by specializing
Theorems~\ref{th:reduction multiplicity Lusztig} and
\ref{th:reduction multiplicity string} to the following choice
of a reduced word $\ii \in R(\wnot (I)^{-1} \wnot)$:
$$\ii = (p,p+1, \dots, p+q-1; p-1, p, \dots, p+q-2; \dots ; 1, 2,
\dots, q) \ .$$
We rename the corresponding variables $t_k$ as follows:
$$t= (t_{11}, \dots, t_{1q}; t_{21}, \dots, t_{2q}; \dots ; t_{p1}, \dots,
t_{pq}) \ .$$
It is now shown by a direct check that all the conditions in
Corollaries~\ref{cor:pq reduction Lusztig} and
\ref{cor:pq reduction string} are specializations of the
corresponding conditions in Theorems~\ref{th:reduction multiplicity Lusztig} and
\ref{th:reduction multiplicity string} (for conditions
Corollary~\ref{cor:pq reduction string}(1)
and Corollary~\ref{cor:pq reduction Lusztig}(3),
this follows from Theorem~\ref{th:fully commutative factors}).
We leave the details of this check to the reader.
\endproof

\section*{Acknowledgements}
We thank Alex Postnikov for very helpful discussions.
Part of this work was done during Andrei Zelevinsky's visit
to the University of Bielefeld in August 1999; he thanks
Claus Ringel and Steffen K\"onig for their hospitality and
interest in this work.

\end{document}